\documentclass[12pt]{amsart}
\usepackage[hmargin=1in,vmargin=1.25in]{geometry}
\usepackage[pdfencoding=auto,psdextra]{hyperref}
\usepackage{textgreek}
\usepackage{xcolor}
\usepackage{tikz}
\usepackage{tikz-cd}
\usetikzlibrary{arrows.meta}
\usetikzlibrary{shapes,calc}

\usepackage{mathtools}
\usepackage{extarrows}
\usepackage{latexsym}
\usepackage{verbatim}
\usepackage{amsmath,amsthm,amssymb}
\usepackage[mathscr]{euscript}
\usepackage{stmaryrd}
\usepackage{multicol}
\usepackage{bm}
\numberwithin{equation}{section}
\newtheorem{thm}{Theorem}[section]
\newtheorem{prop}[thm]{Proposition}
\newtheorem{lem}[thm]{Lemma}
\newtheorem{cor}[thm]{Corollary}
\newtheorem{claim}{Claim}{\bf}{\it}
{\bf}{\it}
\newtheorem{conjecture}[thm]{Conjecture}
\newtheorem{hypo}[thm]{Hypothesis}
\newtheorem{fthm}{Theorem}{\bf}{\it}
{\bf}{\it}
\newtheorem{fcor}[fthm]{Corollary}{\bf}{\it}
{\bf}{\it}
{\bf}{\it}

\theoremstyle{definition}
\newtheorem{defn}[thm]{Definition}

{\bf}{\rm}

\theoremstyle{remark}
\newtheorem{ex}[thm]{Example}
\newtheorem{rem}[thm]{Remark}
{\bf}{\it}

\newtheorem{definition and corollary}[thm]{Definition and Corollary}

{\it}{\rm}

\DeclareMathOperator{\Hom}{Hom}

\DeclareMathOperator{\Ext}{Ext}
\DeclareMathOperator{\ext}{ext}
\DeclareMathOperator{\cend}{end}

\DeclareMathOperator{\GL}{\mathop{GL}}
\DeclareMathOperator{\SL}{\mathop{SL}}
\DeclareMathOperator{\gl}{\mathfrak{gl}}
\DeclareMathOperator*{\gch}{\mathrm{gch}}
\DeclareMathOperator*{\ugch}{\underline{\mathrm{gch}}}
\DeclareMathOperator*{\bgch}{\mathrm{gch}_{q,t}}
\DeclareMathOperator*{\gdim}{\mathrm{gdim}}
\DeclareMathOperator{\coker}{coker} 
\DeclareMathOperator{\Res}{Res}

\DeclareMathOperator{\ch}{ch}
\DeclareMathOperator{\gr}{gr}
\DeclareMathOperator{\Lie}{Lie}
\DeclareMathOperator{\rk}{rk}

\DeclareMathOperator{\hd}{\mathsf{hd}}
\DeclareMathOperator{\soc}{\mathsf{soc}}

\DeclareMathOperator{\SW}{\mathsf{SW}}
\DeclareMathOperator{\WS}{\mathsf{WS}}

\DeclareMathOperator{\Span}{Span}

\newcommand{\gmod}{\mathchar`-\mathsf{gmod}}

\newcommand{\al}{\alpha}
\newcommand{\af}{\mathrm{af}}

\newcommand{\C}{{\mathbb C}}
\newcommand{\cC}{{\mathscr C}}

\newcommand{\cO}{{\mathcal O}}

\newcommand{\tK}{{}^\sharp K}

\newcommand{\hs}{\mathtt{h}}

\newcommand{\la}{\lambda}
\newcommand{\La}{\Lambda}

\newcommand{\HL}{\mathtt{HL}}

\newcommand{\Sym}{\mathfrak{S}}
\newcommand{\tSym}{\widetilde{\mathfrak{S}}}

\newcommand{\sk}{{\mathtt s}^{(k)}}

\newcommand{\g}{\mathfrak{g}}
\newcommand{\tg}{\widetilde{\mathfrak{g}}}
\newcommand{\gb}{\mathfrak{b}}
\newcommand{\tb}{\widetilde{\mathfrak{b}}}
\newcommand{\tB}{\widetilde{\mathbf{B}}}

\newcommand{\gt}{\mathfrak{t}}
\newcommand{\wgt}{\widetilde{\mathfrak{t}}}

\newcommand{\tI}{\mathtt{I}}
\newcommand{\tJ}{\mathtt{J}}

\newcommand{\gn}{\mathfrak{n}}

\newcommand{\gp}{\mathfrak{p}}
\newcommand{\tp}{\widetilde{\mathfrak{p}}}
\newcommand{\gu}{\mathfrak{u}}

\renewcommand{\P}{\mathbb{P}}
\newcommand{\Par}{\mathtt{Par}}
\newcommand{\Comp}{\mathtt{Comp}}

\newcommand{\sX}{\mathscr{X}}

\newcommand{\po}{\mathrm{pol}}
\newcommand{\ssk}{\mathsf{s}^{(k)}}
\newcommand{\ssc}{\mathsf{s}}
\newcommand{\ws}{\widetilde{s}}
\newcommand{\we}{\widetilde{e}}
\newcommand{\sP}{\mathsf{P}}

\newcommand{\Q}{\mathbb{Q}}
\newcommand{\R}{\mathbb{R}}
\newcommand{\bW}{\mathbb{W}}

\newcommand{\Z}{\mathbb{Z}}

\newcommand{\Gm}{\mathbb G_m}

\newcommand{\CatalanInlineRC}[3]{%
  \begin{tikzpicture}[scale=0.40, baseline={([yshift=-.4ex]current bounding box.center)}]
    \foreach \i in {1,...,5}{
      \foreach \j in {1,...,5}{
        \draw[gray!70, dotted] (\j-1,5-\i) rectangle (\j,6-\i);
      }
    }
    \foreach \r/\c in {#3}{
      \path[fill=red!55] (\c-1,5-\r) rectangle (\c,6-\r);
      \draw[red!75, thick] (\c-1,5-\r) rectangle (\c,6-\r);
    }
    \foreach [count=\r from 1] \w in {#1}{
      \node[font=\scriptsize] at (\r-0.5,5-\r+0.5) {\w};
    }
    \pgfmathsetmacro{\yy}{5-(#2)}
    \draw[black, line width=0.25pt] (0,0) rectangle (5,5);
  \end{tikzpicture}%
}
\begin{document}
\title[Categorification of $k$-Schur functions]{Categorification of $k$-Schur functions\\and refined Macdonald positivity}
\author{Syu Kato}
\address{Department of Mathematics, Kyoto University, Oiwake Kita-Shirakawa Sakyo Kyoto 606-8502 JAPAN}
\email{syuchan@math.kyoto-u.ac.jp}
\date{\today}



\begin{abstract}
We establish a substantial part of the categorification framework for $k$-Schur functions proposed in Chen's Ph.D.\ thesis, written under the supervision of Mark Haiman. More precisely, we realize $k$-Schur functions as the graded characters of a distinguished family of objects in a certain module category and establish module-theoretic analogues of several foundational properties of $k$-Schur functions. A final part of Chen's framework predicts a homological characterization of modules admitting filtrations by these distinguished objects; we establish this characterization under an additional hypothesis.

In addition, we unconditionally prove that modules arising in this framework admit filtrations whose subquotients are affine Demazure modules of level $k$. As a consequence, modified Macdonald polynomials expand positively in the characters of affine Demazure modules. This may be viewed as a homological refinement of Macdonald positivity, arising from an intrinsic $\ext$-orthogonality condition on the corresponding Garsia--Haiman modules. We also prove the combinatorial hypothesis for $m \le 2k$ in the appendix and verify it computationally for $m\leq 19$, yielding $k$-Schur positivity of modified Macdonald polynomials in these cases.

Our approach builds on our previous work on the algebraic and geometric realization of Catalan symmetric functions, a class encompassing both $k$-Schur and Hall--Littlewood functions.
\end{abstract}

\maketitle

\tableofcontents

\section*{Introduction}
The positivity of the coefficients in the Schur function expansion of the modified Macdonald polynomials, often referred to as Macdonald positivity, is a central theme in algebraic combinatorics and representation theory. It was conjectured by Macdonald~\cite{Mac88} and established by Haiman~\cite{Hai01}. Beyond positivity itself, one seeks conceptual explanations for this phenomenon and effective descriptions of these coefficients, often through expansions in more suitable bases. These questions remain the subject of active research; see, for example,~\cite{GH09,Gor12,Los22} for structural approaches and~\cite{HHL05,CMW20,AMM23,AMM24} for explicit formulas. One of the most influential approaches is the study of $k$-Schur functions, introduced by Lapointe, Lascoux, and Morse~\cite{LLM03,GR05} in connection with a conjectural refinement of Macdonald positivity.

The significance of $k$-Schur functions was further reinforced by
Lam's identification~\cite{Lam06} of their specialization at $q=1$
with a basis of the homology of the affine Grassmannian. This result
provides a geometric interpretation of $k$-Schur functions and places
them at the intersection of algebraic combinatorics and geometric
representation theory; see
\cite{LM03a,LM03b,LM04,LM05,LM07,BSZ11,LLMS,AB12,BSS12,LLMS13,CNZ16}.

A major combinatorial development in the study of $k$-Schur functions was achieved by Blasiak, Morse, Pun, and Summers~\cite{BMPS,BMPS2,BMP}, who studied a broader framework of Catalan symmetric functions in which the $k$-Schur functions arise naturally. In particular, they established a number of fundamental combinatorial properties of $k$-Schur functions in this setting.

The Catalan symmetric functions studied in those works were introduced in Chen's Ph.D. thesis~\cite{Che10}, written under the supervision of Mark Haiman. Chen formulated two related frameworks in which $k$-Schur functions arise naturally: an algebraic, or module-theoretic, framework~\cite[\S5.1]{Che10}, and a geometric framework~\cite[\S5.4]{Che10}. In~\cite{Kat23a}, we established the geometric framework by proving the vanishing conjecture posed in~\cite[\S5.4]{Che10}.

The purpose of the present paper is to establish the algebraic realization of $k$-Schur functions predicted by the Chen--Haiman framework. We study the modules arising from Chen's geometric framework by interpreting them as quasicoherent sheaves on the variety $\sX_\Psi$ introduced in~\cite{Kat23a}, rather than as vector bundles on flag varieties. We then use an affine version of Schur--Weyl duality~\cite{FKM19,Fl21} to transport these modules to the corresponding algebraic category. This allows us to realize $k$-Schur functions as the graded characters of a distinguished family of objects in that category and to prove module-theoretic analogues of several foundational properties of $k$-Schur functions. In this sense, the present paper establishes a substantial part of the categorification framework for $k$-Schur functions proposed in Chen's thesis.

To describe our results more concretely, we begin by fixing three positive integers $n > m$ and $k$. We set
\[
\g[z] := \C \cdot \mathrm{Id} \oplus \mathfrak{sl}(n, \C[z]).
\]
We also set
\[
A = A_m := \C \Sym_m \ltimes \C[X_1,\ldots,X_m].
\]
Both $\g[z]$ and $A$ are naturally graded. We consider the category $\g[z]\gmod_m$ of finitely generated graded $\g[z]$-modules on which the diagonal matrix units $E_{ii}$ ($1 \le i \le n$) act semisimply with eigenvalues in $\Z_{\le 0}$ and $\mathrm{Id}$ acts by $-m$, and the category $A\gmod$ of finitely generated graded $A$-modules.

We have the Schur--Weyl duality functor
\[
\SW : \g[z]\gmod_m \xrightarrow{\ \cong\ } A\gmod,
\]
which is an equivalence of categories~\cite{FKM19,Fl21}. The simple objects on both sides are parameterized by the set $\Par_m$ of partitions of $m$. For each $\la\in\Par_m$, we denote by $V_\la^*$ (resp. $L_\la$) the corresponding simple object on the $\g[z]$-side (resp. the $A$-side), normalized to be concentrated in degree~$0$. We write $\mathsf{q}$ for the compatible grading shift functor on both categories. For $\la \in \Par_m$, we write $\la_1$ for the first part of $\la$, and $\la'$ for the conjugate partition.

A Catalan symmetric polynomial $H\!L^\Psi_\la$ is a symmetric polynomial in $n$ variables, depending on a root ideal $\Psi$ of size $n$ (Definition~\ref{defn:ri}) and a partition $\la \in \Par_m$. By~\cite{Kat23a}, each $H\!L^\Psi_\la$ can be interpreted as the graded character of a module $\HL^\Psi_\la$ in the category $\g[z]\gmod_m$:
\[
H\!L^\Psi_\la = \gch \HL^\Psi_\la \in \bigoplus_{\mu \in \Par_m} \Z[q] s_\mu,
\]
where $s_\mu$ denotes the Schur polynomial associated with $\mu \in \Par_m$. We denote by $\HL_\la$ the module corresponding to the case where $\Psi$ is the maximal root ideal.

Let $\Par_m^{(k)}$ denote the set of $k$-bounded partitions of $m$, that is, those $\la \in \Par_m$ satisfying $\la_1 \le k$. There is a notion of $k$-conjugation $\la \mapsto \la^{\omega_k}$, which defines an involution on $\Par_m^{(k)}$~\eqref{eqn:defkinv}. As shown in~\cite{BMPS}, the $k$-Schur polynomial $s^{(k)}_\la$ can be realized as the Catalan symmetric polynomial $H\!L^{\Psi[\la,k]}_\la$, where $\Psi[\la,k]$ is a root ideal determined by $k$ and $\la \in \Par_m^{(k)}$. We set
\[
\sk_\la := \HL^{\Psi[\la,k]}_\la.
\]
In~\cite{Kat23a}, we proved that $\HL^\Psi_\la$ admits a unique simple quotient for any root ideal $\Psi$. However, the arguments there suggest that $\HL^\Psi_\la$ may admit more than one simple submodule in general. In this paper, we provide a general sufficient condition under which the module $\HL^\Psi_\la$ has a unique simple submodule (Theorem~\ref{thm:HLss}). Consequently, we obtain:

\begin{fthm}[$\doteq$ Theorem~\ref{thm:Cat-sh} $+$ Theorem~\ref{thm:HLss}]\label{fthm:ks}
For each $\la \in \Par^{(k)}_m$, the module $\sk_\la$ has a unique simple quotient $V_\la^*$ and a unique simple submodule.
\end{fthm}

For each $\la \in \Par^{(k)}_m$, we set $\ssc_\la^{(k)} := \SW ( \sk_\la )$. We have an autoequivalence $\circledast$ on $A\gmod$ induced by tensoring with $\mathsf{sgn}$. From Theorem~\ref{fthm:ks}, we derive the following result, originally conjectured in~\cite[Conjecture~4.4.3]{Che10}:

\begin{fthm}[$\doteq$ Theorem~\ref{thm:CHid} $+$ Corollary~\ref{cor:kdual}]\label{fthm:char}
For each $\la \in \Par^{(k)}_m$, we have
$$\ssc_\la^{(k)} \cong \mathsf q^{d_k (\la)} \bigl( ( \ssc_{\la^{\omega_k}}^{(k)} )^{*} \bigr)^{\circledast},$$
for some $d_k(\la) \in \Z_{\ge 0}$ specified in~\eqref{eqn:defkinv}. Moreover, $\ssc_\la^{(k)}$ is precisely the module $M_{\la,(\la^{\omega_k})'}$ in~\cite[\S4.2]{Che10}.
\end{fthm}

We further analyze our geometric realization of $\HL^\Psi_\la$, with particular emphasis on $\sk_\la$, in Sections~\ref{sec:gpr}--\ref{sec:kbr}. This yields geometric recursions between modules (Theorem~\ref{thm:assgr}, Corollary~\ref{cor:kPieri}, and Theorem~\ref{thm:Fdp}), together with new vanishing theorems for sheaf cohomology on $T^*_\Psi X$. These geometric results imply the following module-theoretic refinement of~\cite[Theorem~2.6]{BMPS} and~\cite[Theorem~4.1]{BMPS2}:

\begin{fthm}[$\doteq$ Corollary~\ref{cor:sskbranch}
$+$ Corollary~\ref{cor:HLfilt}]\label{fthm:filt}
For each $\la \in \Par^{(k)}_m$, the following assertions hold:
\begin{itemize}
\item $\sk_\la$ admits a finite filtration by the grading shifts of the modules
$\{\mathtt{s}_\mu^{(k+1)}\}_{\mu \in \Par_m^{(k+1)}}$.
\item $\HL_\la$ admits a finite filtration by the grading shifts of the modules
$\{\sk_\mu\}_{\mu \in \Par_m^{(k)}}$.
\end{itemize}
\end{fthm}

We denote by $\ext^{\bullet}_{A} ( \bullet, \bullet )$ the graded extension groups in $A \gmod$. Applying $\SW$ to Theorem~\ref{fthm:filt} and combining this with the results\footnote{Alternatively, one can import the $\ext$-orthogonality results from~\cite{CL06,CI15,Kle14} via $\SW$ and apply the twist $\circledast$ as in Remark~\ref{rem:Atw}.} recalled in Section~\ref{subsec:Amod}, we derive:

\begin{fcor}[$\doteq$ Corollary~\ref{cor:korth}]\label{fcor:orth}
For each $\la \in \Par^{(k)}_m$, we have
\begin{equation}
\ext^\bullet_{A} ( \ssc_\la^{(k)}, L_{\mu} ) \equiv 0 \qquad (\mu \in \Par_m \setminus \Par_m^{(k)}).\label{eqn:Intro-extorth}
\end{equation}
\end{fcor}

The preceding results give an algebraic realization of $k$-Schur functions through a distinguished family of modules, establish their structural properties, and show that they satisfy the $\ext$-orthogonality condition~\eqref{eqn:Intro-extorth}.

We next study the consequences of the $\ext$-orthogonality condition~\eqref{eqn:Intro-extorth}.

\begin{fthm}[$\doteq$ Theorem~\ref{thm:Dstr}]\label{fthm:main}
Let $M \in A \gmod$ be a finite-dimensional module. If $M$ satisfies
\begin{equation}
\ext^{\bullet} _{A} ( M, L_\mu ) = 0 \qquad (\mu \notin \Par_m^{(k)})\label{eqn:fextLk},
\end{equation}
then $M$ admits a finite filtration whose subquotients are grading
shifts of modules $\mathsf{W}_\la^{(k)}$ obtained, via $\SW$, from
affine Demazure modules of level $k$; see
Section~\ref{subsec:level-k}.
\end{fthm}

Taking the graded characters of particular modules appearing in Theorem~\ref{fthm:main}, we obtain:

\begin{fcor}[affine Demazure positivity of Macdonald polynomials $\doteq$ Corollary~\ref{cor:dMac}]\label{fcor:dmp}
Let $\la \in \Par_m^{(k)}$ and let $\mathsf{GH}_\la$ be the Garsia--Haiman module~\cite{GH93}. Here we use the convention
$\bgch \mathsf{GH}_\la = \widetilde{H}_\la(t,q)$~\cite{Hai01},
where $\bgch$ denotes the bigraded Frobenius characteristic. Then
\[
\bgch \mathsf{GH}_\la^{\circledast} = \omega \widetilde{H}_\la(t,q) \in \sum_{\mu \in \Par_m} \Z_{\ge 0}[q^{\pm 1},t] \gch \mathsf{W}_\mu^{(k)},
\]
where $\omega$ is a $\Q(q,t)$-linear involution satisfying $\omega(s_\mu)=s_{\mu'}$ for $\mu \in \Par_m$.
\end{fcor}

For a fixed partition $\la$, the character $\bgch\mathsf{GH}_\la$ is independent of $k$, and every $k\ge\la_1$ is admissible. Thus Corollary~\ref{fcor:dmp} provides a nontrivial refinement of Macdonald positivity for a range of levels $k$.

In view of~\cite[Theorem~4.1]{Kat23a}, Theorem~\ref{fthm:main} and Corollary~\ref{fcor:dmp} may be regarded as homological refinements of Macdonald positivity and as weaker forms of the $k$-Schur positivity predicted by Chen's algebraic framework. We next turn to the final part of that framework, which predicts a homological characterization of modules admitting $k$-Schur filtrations. We reduce the proof of the original statement to a filtration-inheritance property formulated as Conjecture~\ref{conj:ker}:

\begin{fthm}[$\doteq$ Theorem~\ref{thm:str}]\label{fthm:cmain}
Assume Conjecture~\ref{conj:ker}. Let $M \in A\gmod$ be a finite-dimensional module. Then $M$ admits a finite filtration whose subquotients are grading shifts of the modules
$\{ \ssc_\la^{(k)}\}_{\la\in\Par_m^{(k)}}$ if and only if
\[
\ext_A^\bullet(M,L_\mu)=0
\qquad
(\mu\notin\Par_m^{(k)}).
\]
\end{fthm}

We further formulate an explicit combinatorial sufficient condition
for Conjecture~\ref{conj:ker}, stated as
Hypothesis~\ref{hypo:cycle}. Although Hypothesis~\ref{hypo:cycle} fails in full generality, we
have verified it computationally for $m\leq19$. In the appendix, we prove the hypothesis whenever $m\leq 2k$.

As an application of Theorem~\ref{fthm:cmain}, we obtain:

\begin{fcor}[refined Macdonald positivity $\doteq$ Corollary~\ref{cor:rMac}]\label{fcor:rmp}
Assume Conjecture~\ref{conj:ker}. For each $\la \in \Par_m^{(k)}$, we have
\[
\omega \widetilde{H}_\la(t,q) \in \sum_{\mu \in \Par_m^{(k)}} \Z_{\ge 0}[q,t] s_\mu^{(k)}.
\]
\end{fcor}

Thus Theorem~\ref{fthm:cmain} gives the refined Macdonald positivity statement in Chen's formulation~\cite{Che10}, conditional on Conjecture~\ref{conj:ker}. Its relation to the original formulation of Lapointe--Lascoux--Morse~\cite[Conjecture~8]{LLM03} is discussed in Remark~\ref{rem:ksch}.

This paper is organized as follows. In Section~\ref{sec:prelim}, we fix notation and recall the background material used throughout, including results from~\cite{Kat15,Kat26b,Kat23a} and the representation theory of current algebras~\cite{CI15,FKM19}. In Section~\ref{sec:mtp}, we introduce our basic modules, and establish their fundamental properties.

In Section~\ref{sec:as}, we recall further constructions from~\cite{Kat23a} and prove that the relevant classes of $\g[z]$-modules, including the modules $\sk_\la$, have simple socles. Section~\ref{sec:ksch} refines the resulting characterization of $\sk_\la$, proves Theorem~\ref{fthm:char}, and develops the module-theoretic input required in the subsequent sections.

Sections~\ref{sec:tvb},~\ref{sec:gpr}, and~\ref{sec:kst} establish geometric counterparts of the Pieri rule and of the $k$-Schur straightening rule. In Section~\ref{sec:kbr}, we derive the module-theoretic $k$-Schur branching rules summarized in Theorem~\ref{fthm:filt}. Together with the results recalled in Section~\ref{subsec:Amod}, these yield the $\ext$-orthogonality properties of the families $\{\ssc_\la^{(k)}\}_\la$ and $\{\sk_\la\}_\la$.

In Section~\ref{sec:main}, we construct filtered resolutions and use them, together with Lemma~\ref{lem:fres} and Theorem~\ref{fthm:char}, to derive criteria for the existence of $\ssc^{(k)}$-filtrations, including the filtration criterion in Theorem~\ref{thm:kfilt}. In Section~\ref{sec:comb}, we formulate the combinatorial Hypothesis~\ref{hypo:cycle} and prove that it implies the filtration-inheritance property stated in Conjecture~\ref{conj:ker}.

Finally, in Section~\ref{sec:CH}, we combine these results to prove Theorem~\ref{fthm:main} and Corollary~\ref{fcor:dmp}. Under Conjecture~\ref{conj:ker}, the same argument yields the stronger Theorem~\ref{fthm:cmain} and Corollary~\ref{fcor:rmp}.

Earlier versions of this paper claimed that
Theorem~\ref{fthm:cmain} holds without additional assumptions.
The argument given there was incomplete: it did not establish that
the kernel of a surjection between $\ssc^{(k)}$-filtered modules is
again $\ssc^{(k)}$-filtered. In the present version, this inheritance
property is isolated as Conjecture~\ref{conj:ker}, and a sufficient
combinatorial condition is formulated in
Hypothesis~\ref{hypo:cycle}. The author apologizes for the error and
for any confusion caused by the earlier versions.

\section{Preparatory materials}\label{sec:prelim}

\subsection{Conventions}\label{subsec:conv}

We always work over the field of complex numbers $\C$. By a graded vector space $M$, we mean a $\Z$-graded vector space $M = \bigoplus_{i \in \Z} M_i$. The graded dimension of $M$ is
$$\gdim M := \sum_{i \in \Z} q^i \dim M_i.$$

Let $\mathsf{gVec}$ denote the category of graded vector spaces. We define the restricted dual of $M \in \mathsf{gVec}$ as $M^{\vee} := \bigoplus_{i \in \Z} ( M_i )^*$, with the convention that $( M_i )^*$ has degree $-i$. Note that when $M$ is finite-dimensional, we have $M^{\vee} = M^*$. A graded algebra $E$ is a $\C$-algebra whose underlying vector space is a graded vector space satisfying the condition $E_i \cdot E_j \subset E_{i+j}$ for all $i,j\in \Z$. Similarly, a graded Lie algebra is a Lie algebra $\mathfrak e=\bigoplus_i\mathfrak e_i$ satisfying $[\mathfrak e_i,\mathfrak e_j]\subset\mathfrak e_{i+j}$. This grading induces a graded algebra structure on $U(\mathfrak e)$. A graded module $M = \oplus_{i \in \Z} M_i$ over a graded algebra $E$ is an $E$-module such that $E_i \cdot M_j \subset M_{i+j}$ for all $i,j\in \Z$. The category of graded $E$-modules, denoted by $E\mathchar`-\mathsf{gMod}$, and the category of finitely generated graded $E$-modules, denoted by $E\gmod$, are both equipped with autoequivalences $\mathsf q$ and $\mathsf q^{-1}$, which act on a graded module $M$ as follows:
$$\left( \mathsf{q} ( M ) \right) _i = M_{i-1}, \quad \text{and} \quad \left( \mathsf{q}^{-1} ( M ) \right)_i = M_{i+1} \qquad (i \in \Z).$$

We define the graded $\Ext$-group as
$$\ext^i_{E} ( M, N ) := \bigoplus_{m \in \Z}\Ext^i_{E\mathchar`-\mathsf{gMod}} ( \mathsf{q}^m M, N )$$
for each $M,N \in E\mathchar`-\mathsf{gMod}$ and $i \in \Z$, and regard it as a graded vector space. More generally, for an abelian category $\cC$ equipped with grading shift functors $\mathsf{q},\mathsf{q}^{-1}$, we define
$$\ext^i_{\cC} ( M, N ) := \bigoplus_{m \in \Z}\Ext^i_{\cC} ( \mathsf{q}^m M, N ).$$
Here $\ext^0_E(M,N)$ and $\ext^0_{\cC}(M,N)$ are also denoted by $\hom_E(M,N)$ and $\hom_{\cC}(M,N)$, respectively. We also write $\cend_E(M):=\hom_E(M,M)$. For brevity, we may omit the grading shift functor $\mathsf{q}^{\bullet}$ when the precise grading shifts are not important.

An increasing exhaustive filtration of a module (or a graded module) $M$ is an increasing sequence
$$0 = F_0 M\subset F_1 M \subset F_2 M \subset \cdots$$
of submodules (or graded submodules) that satisfies the condition
$$\bigcup_i F_i M = M.$$
It is finite if we have $F_i M = M$ for all sufficiently large $i$. Its associated graded is defined as
$$\gr_F M = \bigoplus_{i \ge 0}F_{i+1} M / F_{i} M.$$

Dually, a decreasing separable filtration of a module $M$ is a decreasing sequence
$$M = F^0 M \supset F^1 M \supset F^2 M \supset \cdots$$
of submodules such that $\bigcap_i F^i M = \{0\}$. We define its associated graded as
$$\gr^F M = \bigoplus_{i \ge 0} F^{i} M/ F^{i+1} M.$$

For a collection $\{Z_\la\}_\la$ of graded modules, we say that a graded module $M$ admits a filtration (resp. a finite filtration) by $\{Z_\la\}_\la$ if $M$ admits a decreasing separable filtration (resp. a finite filtration) such that its associated graded $\gr M$ is a direct sum of grading shifts of modules in $\{Z_\la\}_\la$.

\subsection{Partitions and symmetric functions}\label{subsec:part}

Standard references for this subsection include~\cite{Mac95,LLM03,LM05}.

For a nonnegative integer $m$, a partition $\la$ of $m$ is defined as a non-increasing sequence of nonnegative integers
$$\la = ( \la_1, \la_2,\la_3,\ldots),\quad \la_1 \ge \la_2 \ge \la_3 \ge \cdots$$
such that $\sum_{i \ge 1} \la_i = m$. We refer to $m$ as the size of $\la$ and denote it by $|\la|$. We set $\ell ( \la ) := \max \bigl(\{ i \mid \la_i > 0 \}\cup\{0\}\bigr)$ and call it the length of $\la$. Let $\Par_m$ be the set of partitions of $m$. The conjugate $\la' \in \Par_m$ of a partition $\la \in \Par_m$ is defined by
$$( \la' )_i := \# \{ j \mid \la_j \ge i \}.$$
We define
\[
\mathsf{m} ( \la ) := \sum_{i\ge 1} \frac{\la_i ( \la_i - 1)}{2} = \sum_{i\ge 1} \sum_{j=0}^{\la_i-1} j, \quad \text{and} \quad \mathsf{n} ( \la ) := \mathsf{m} ( \la' ) = \sum_{i \ge 1} (i-1) \la_i \qquad ( \la \in \Par_m).
\]
Set $\Par := \bigsqcup_{m \ge 0} \Par_m$. An element $(x,y) \in ( \Z_{\ge 0} )^2$ is called a box of a partition $\la \in \Par$ if and only if $0 \le x < \la_{y+1}$. The hook length of $\la$ at a box $(x,y)$ is defined as
$$\la_{y+1} - x + (\la')_{x+1} - y - 1.$$
For a positive integer $r$, a partition $\la \in \Par$ is called an $r$-core if there is no box of $\la$ whose hook length is $r$. For $\la, \mu \in \Par$, we write $\la \subset \mu$ if $\la_i \le \mu_i$ for each $i \ge 1$. We refer to a pair $(\la,\mu)$ of partitions such that $\la \subset \mu$ as a skew partition and denote it by $\mu/\la$.

Let ${\bm \La}$ denote the ring of symmetric functions over $\Z$, and let ${\bm \La}_q$ and ${\bm \La}_{q,t}$ denote the rings of symmetric functions over $\Q(\!(q)\!)$ and $\Q(\!(q,t)\!)$, respectively. The set of Schur functions $\{\ws_\la\}_{\la \in \Par}$, the set of elementary symmetric functions $\{\we_\la\}_{\la \in \Par}$, and the set of power-sum symmetric functions $\{\widetilde{p}_\la\}_{\la \in \Par}$ form bases of ${\bm \La}_q$. For each $f \in {\bm \La}$, we have its adjoint operator $f^{\perp} \in \mathrm{End} ( {\bm \La} )$ (\cite[Chap~I~\S5~Example~3]{Mac95}). We may extend $f^{\perp}$ to a $\Q(\!(q)\!)$-linear endomorphism on ${\bm \La}_q$.

For $k \in \Z_{\ge 1}$, define
$$\Par^{(k)} := \{\la \in \Par \mid \la_1 \le k \}$$
and call its elements $k$-bounded partitions. Note that $\Par^{(0)} = \{\emptyset\}$.

\begin{thm}[{\cite[\S3]{LM05}}]\label{thm:kcore}
For $\la \in \Par^{(k)}$, there exists a unique $(k+1)$-core $\mu \in \Par$ and $\nu \subset \mu$ such that a box of $\mu$ lies in $\nu$ if and only if it has hook length $> (k+1)$ and
\[
\la_i = \mu_i - \nu_i \qquad (i \ge 1).
\]
\end{thm}

Using the skew partition $\mu/\nu$ associated with $\la \in \Par^{(k)}$ through Theorem~\ref{thm:kcore}, we define $d_k (\la) \in \Z_{\ge 0}$ and $\la^{\omega_k} \in \Par$ as:
\begin{equation}
d_k (\la) := |\nu|, \qquad (\la^{\omega_k})_i := ( \mu' )_i - ( \nu' )_i \qquad (i \ge 1).\label{eqn:defkinv}
\end{equation}

\begin{ex}\label{ex:n7conj}
Let $\la =(6,5^2,3,1^2) \in \Par^{(7)}$. The skew partition $\mu/\nu$ obtained from $\la$ in Theorem~\ref{thm:kcore} can be computed by the procedure described in~\cite[Definition~32]{LLM03}, which yields
$$\nu = (6,1^2) \subset (12,6^2,3,1^2) =\mu.$$
Thus, we have $\nu' = (3,1^5)$, $\mu' = (6,4^2,3^3,1^6)$, and
$$\la^{\omega_7} = (3^3,2^3,1^6).$$
\end{ex}

\begin{thm}[{\cite[\S7]{LLM03}}]\label{thm:kinv}
For each $\la \in \Par^{(k)}$, we have $\la^{\omega_k} \in \Par^{(k)}$. In addition, $\omega_k$ defines an involution on $\Par^{(k)}$.
\end{thm}

\begin{cor}\label{cor:dktr}
In the setting of Theorem~\ref{thm:kinv}, we have $d_k ( \la ) = d_k ( \la^{\omega_k} )$.
\end{cor}

Note that we have $d_k ( \la ) = 0$ and $\la^{\omega_k} = \la'$ whenever $|\la| \le k$ (see~\cite[P138]{LLM03}).

\subsection{Algebraic groups and root systems}\label{subsec:roots}

For the basic material contained here, we refer to~\cite{CG97,Kum02}.

Let $n \in \Z_{>0}$ and $m \in \Z_{\ge 0}$. We set $G := \GL (n,\C)$. Let $T$ be the subgroup of diagonal matrices of $G$, and let $B$ be the subgroup of upper-triangular matrices of $G$. We follow the standard convention of denoting Lie algebras by the corresponding fraktur letters (e.g., $\g = \mathrm{Lie} \, G, \gt = \mathrm{Lie} \, T,\ldots$). We set $\g_\circ := \mathfrak{sl} (n,\C) \subset \g$, $\gt_\circ := ( \gt \cap \g_\circ )$, and $\gb_\circ := ( \gb \cap \g_\circ )$. We also set
\begin{eqnarray*}
\g [z] & := & \C \mathrm{Id} \oplus ( \g_\circ \otimes \C [z] ) \subset \mathfrak{gl} (n,\C [z])\quad \text{and}\\
G [\![z]\!] & := & ( \C^{\times} \mathrm{Id} ) \cdot \SL (n, \C [\![z]\!] ) \subset \GL (n, \C [\![z]\!] ).
\end{eqnarray*}

Note that $\gt \subset \g[z]$, $T \subset G[\![z]\!]$, and $\g[z]$ is a $\Z$-graded Lie algebra in which $\g_\circ \otimes z^s$ has degree $s$. We also set $\widehat{G}[\![z]\!] := \Gm \ltimes G[\![z]\!]$, where $\Gm$ acts on $G[\![z]\!]$ by the dilation of $z$ (so-called loop rotation action). We may add $n$ as a subscript in order to clarify the dependence of these objects on the choice of $n$; for instance, we write $G_n = G = \GL(n,\C)$, $\g_n[z] = \g[z]$, $G_n[\![z]\!] = G[\![z]\!]$, etc.

We fix the primitive characters $\epsilon_1,\ldots,\epsilon_n$ of $T$, where each $\epsilon_i$ extracts the $i$-th diagonal entry. We set
\begin{align*}
\sP := & \bigoplus_{i=1}^n \Z \epsilon_i \supset \sP^+ := \{\sum_i \la_i \epsilon_i \in \sP \mid \la_1\ge \la_2 \ge \cdots \ge \la_n \},\\
\Comp ( n ) & := \{\sum_i \la_i \epsilon_i \in \sP \mid \la_i \ge 0 \ \text{for all} \ i \}, \qquad \Par ( n ) := \Comp ( n ) \cap \sP^+.
\end{align*}

The elements of $\Par(n)$ are called dominant polynomial weights of $\g$. The symmetric group $\Sym_n$ acts on $\sP$ and we have a non-degenerate $\Sym_n$-invariant pairing $\langle \bullet, \bullet \rangle$ on $\sP$ defined as $\langle \epsilon_i,\epsilon_j \rangle = \delta_{ij}$. For $\la = \sum_{i} \la_i \epsilon_i \in \sP$, we set $|\la| := \sum_i \la_i$. We have $|\la| \ge 0$ if $\la \in \Comp (n)$. Unless otherwise specified, we understand a weight $\la \in \sP$ to be expressed as $\sum_i \la_i \epsilon_i$. We might also represent this as $\la = (\la_1,\ldots,\la_n)$. For $\la, \mu \in \Par(n)$, we say $\la \subset \mu$ if $\la_i \le \mu_i$ for $1 \le i \le n$. We set
\begin{align*}
\Comp_m ( n ) & := \{\la \in \Comp (n) \mid |\la| = m\}, & \Par_m ( n ) := \{\la \in \Par (n) \mid |\la| = m\},\\
\Par^{(k)}_m ( n ) & := \{\la \in \Par_m(n) \mid \la_1 \le k\}. &
\end{align*}

We freely identify $\Par_m(n)$ with the subset of $\Par_m$ consisting of partitions of $m$ of length at most $n$. When $n \ge m$, we may omit $(n)$ in $\Comp_m(n)$, $\Par(n)$, $\Par_m(n)$, and $\Par^{(k)}_m(n)$, since padding with additional zeros identifies some of these objects with those introduced in Section~\ref{subsec:part}. In particular, this identifies $\Par(n)$ with a subset of $\Par$, in a way compatible with the calculation of $|\la|$.

We set $\varpi_i := \sum_{j=1}^i \epsilon_j \in \Par(n)$ for $0 \le i \le n$ (i.e., $\varpi_0 = 0$). We set $\al_i := (\epsilon_i - \epsilon_{i+1}) \in \sP$ for $1\le i < n$, $\al_{ij} := (\epsilon_i - \epsilon_j) \in \sP$ for $1 \le i, j \le n$, and
$$\Delta^+ := \{\al_{ij} \mid 1 \le i < j \le n\} \subset \sP.$$

Let $E_{ij} \in \g$ be the matrix unit corresponding to $1 \le i, j \le n$ (whose $T$-weight is $\al_{ij}$). We set $\gn := \bigoplus_{1\le i<j\le n}\C E_{ij}=[\gb,\gb]$. We introduce the dominance order $\le$ on $\sP$ by
$$\la \le \mu \qquad \Leftrightarrow \qquad \mu \in \la + \sum_{\beta \in \Delta^+} \Z_{\ge 0} \beta.$$

For subsets of $\sP$, we may consider the partial order induced by the dominance order on $\sP$, which we will also refer to as the dominance order for simplicity. Note that $\la, \mu \in \sP$ are comparable with respect to $\le$ only if $|\la| = |\mu|$.

\begin{lem}\label{lem:interval}
If $\la, \mu \in \sP^+$ satisfy $\la \ge \mu$ and $\la \in \Par(n)$, then we have $\mu \in \Par(n)$. If $\la \in \Par_m^{(k)} ( n )$ and $\mu \in \Par_m ( n )$ satisfies $\la \ge \mu$, then we have $\mu \in \Par_m^{(k)} ( n )$.
\end{lem}

\begin{proof}
We prove the first assertion. Note that $\la \in \Par ( n )$ implies $\la_n \ge 0$. This condition is preserved under subtraction of an element in $\Delta^+$. Therefore, the first assertion follows. For the second assertion, it suffices to note that $\mu_1 \le \la_1 \le k$.
\end{proof}

We denote the flag variety of $G$ by $X := G / B$. We have the canonical projection
$$\pi : T^* X \to X.$$

For each $\la \in \sP$, we consider the line bundle
$$\mathcal L_{X} ( \la ) := G \times ^{B} \C_{-\la}$$
over $X$ whose associated sheaf of sections is denoted by $\cO_{X} ( \la )$. We set $\cO_{T^* X} ( \la ) := \pi^* \cO_{X} ( \la )$.

For each $1 \le i < n$, we consider the parabolic subgroup $B \subset P_i \subset G$ and its Levi subgroup $\SL^{+}(2,i)$ stable under the adjoint $T$-action such that $\gp_i = \C F_i \oplus \gb$ (as vector spaces), where we set $F_i := E_{i+1,i}$. We also set $B_i := \SL^{+}(2,i) \cap B$.

For a semisimple $T$-module $V$, we define its character as
$$\ch V := \sum_{\la \in \sP} e^\la \dim \Hom_{T} (\C_{-\la}, V).$$

We set
$$V_\la := H^0 ( X, \cO_{X} ( \la ))^* \quad \text{and} \quad s_\la := \ch V_\la^* \qquad (\la \in \sP^+).$$
We call $V_\la^*$ ($\la \in \Par(n)$) a polynomial representation of $G$. This terminology is justified by the fact that $s_\la \in \C [X_1,\ldots,X_n]$, where $X_i = e^{\epsilon_i}$ ($1 \le i \le n$), if and only if $\la \in \Par(n)$. In this case, $s_\la$ is the image of $\ws_\la \in {\bm \La}_{q}$ under the truncation map, i.e., we have a graded ring map
\begin{equation}
\mathsf{tr}_n: {\bm \La}_q \longrightarrow \C (\!(q)\!) [X_1,\ldots,X_n]^{\Sym_n}\label{eqn:trmap}
\end{equation}
such that $\mathsf{tr}_n ( \ws_\la ) = s_\la$. In a similar manner, $e_\la := \mathsf{tr}_n ( \we_\la )$ is the elementary symmetric polynomial in $n$ variables.

The vector space $V_\la$ acquires the structure of a $G$-module, and $\{V_\la\}_{\la \in \sP^+}$ is a complete set of representatives of the isomorphism classes of irreducible finite-dimensional rational $G$-modules. $V_\la$ is also a finite-dimensional irreducible $\g$-module via differentiation. Hence $\{V_\la\}_{\la \in \sP^+}$ also forms a complete set of representatives of the isomorphism classes of irreducible finite-dimensional $\g$-modules whose $\gt$-weights extend to $T$-weights by integration (or equivalently belong to $\sP \subset \gt^*$). Here $V_\la$ has a unique $B$-eigenline of $T$-weight $\la$.

We have an automorphism $\overline{\bullet} : \sP^+ \rightarrow \sP^+$ given by
$$\overline{\la} = ( - \la_n,\ldots,-\la_1 ) = - w_0 \la \qquad \text{when} \qquad \la = (\la_1,\ldots,\la_n),$$
where $w_0 \in \Sym_n$ is the longest element. We have an isomorphism $( V_\la )^* \cong V_{\overline{\la}}$ for $\la \in \sP^+$.

\subsection{Affine Demazure functors}\label{subsec:AffDem}

For preliminaries relevant to this subsection, we refer to~\cite{Kum02}.

Let $\tg_{\circ}$ be the untwisted affine Kac-Moody algebra corresponding to $\mathfrak{sl} (n,\C)$. We set
$$\tg := \C \mathrm{Id} \oplus \tg_{\circ},$$
where $\mathrm{Id}$ is a central element that records the action of the center of $\mathfrak{gl}(n,\C)$. We index the Dynkin nodes of $\mathfrak{sl} (n,\C)$ and $\tg_{\circ}$ as $\tI := \{1,2,\ldots,(n-1)\}$ and $\tI_\af = \{0\} \cup \tI$, respectively.

Let $\wgt$ be the Cartan subalgebra of $\tg$ which contains the Cartan subalgebra $\gt$ of $\g$. Let $\wp$ be the level one fundamental weight of $\wgt$ and set
$$\La_i := \varpi_i + \wp \qquad (0 \le i \le n).$$
We have an imaginary root $\delta \in \wgt^*$ such that
\begin{align*}
\wgt^* \cong \gt^* \oplus \C \wp \oplus \C \delta, \qquad & \wgt = \gt \oplus \C K \oplus \C d, \quad \text{and}\\
\wp ( K ) = 1, \wp ( \gt ) = 0 = \wp ( d ), & \quad \delta ( d ) = 1.
\end{align*}

Let $\tSym_n$ be the (affine) Weyl group of $\tg$ generated by $\{s_i\}_{i=0}^{n-1}$. We set
$$\sP_\af := \left( \bigoplus_{i \in \tI_\af} \Z \La_i \right) \oplus \Z \varpi_n \oplus \Z \delta \supset \sP_\af^+ := \left( \sum_{i \in \tI_\af} \Z_{\ge 0} \La_i \right) \oplus \Z \varpi_n \oplus \Z \delta \subset \wgt^*$$
and
$$\sP_\af^\po := \{\La \mid \La = w \La_+, \ w \in \tSym_n, \ \La_+ \in \sP_\af^+, \ \La \! \mid_{\gt} \in \Comp(n) \}.$$

We regard $\al_i \in \gt^* \subset \wgt^*$ for $1 \le i < n$. Together with $\al_0 = \al_{n1} + \delta \in \wgt^*$, these form the set of simple roots of $\tg_{\circ} \subset \tg$. We define the affine dominance order $\lhd$ on $\sP_\af$ by
$$\La \lhd \La' \quad \Leftrightarrow \quad \La' \in \La + \sum_{i \in \tI_\af} \Z_{\ge 0} \al_i.$$

For each $i \in \tI_\af$, we have a positive Kac-Moody generator $E_i \in \tg_{\circ}$ and a negative Kac-Moody generator $F_i \in \tg_{\circ}$ with $\wgt$-weight $\al_i$ and $-\al_i$, respectively. Together with $\wgt$, the set $\{E_i, F_i\}_{i \in \tI_\af}$ generates $\tg$. (We choose $F_i = E_{i+1,i}$ for $1 \le i < n$, consistent with our earlier notation.) We have a linear action of $\tSym_n$ on $\wgt^*$ given by
$$s_i ( \La _j + l \varpi_n + m \delta ) := \La_j - \delta_{ij} \al_i + l \varpi_n + m \delta \qquad (0 \le i,j < n, \quad l, m \in \Z).$$
Restriction to $\gt$ gives the projection
$$\mathsf{pr} : \sP_\af \ni \La_i \mapsto \varpi_i \in \sP \qquad (0 \le i \le n).$$

We fix an Iwahori subalgebra $\tb \subset \tg$ generated by $\wgt$ and the set of positive Kac-Moody generators $\{ E_i \}_{i \in \tI_\af}$, and the minimal parabolic subalgebra
$$\tb \subset \tp_i \subset \tg \qquad (0 \le i < n)$$
such that $\tp_i = \C F_i \oplus \tb$. Note that $\dim \tp_i/\tb = 1$. We have a unique subalgebra $\mathfrak{sl}(2,i) \subset \tp_i$ that is normalized by $\wgt$ and is isomorphic to $\mathfrak{sl}(2,\C)$ as a Lie algebra. We have $\mathfrak{sl}(2,i) = \Span \{E_i,[E_i,F_i],F_i\}$, and these generators form a $\mathfrak{sl}(2)$-triple. In addition, we have $( \mathfrak{sl}(2,i) + \gt ) = \Lie \SL^+ (2,i)$ for $1 \le i < n$. We have a direct sum decomposition
$$\tp_i = \bigl( \mathfrak{sl} (2,i) + \wgt \bigr) \oplus \gu_i \qquad (0 \le i < n),$$
and both summands are $\wgt$-stable Lie subalgebras of $\tg$.

For a $\wgt$-semisimple $\tb$-module $M$ on which each $E_j$
($j\in\tI_\af$) acts locally nilpotently, let $\mathscr D_i ( M )$ be the maximal $\tp_i$-module quotient of
$$U ( \tp_i ) \otimes_{U(\tb)} M$$
whose restriction to $\mathfrak{sl}(2,i)$ is a direct sum of finite-dimensional modules.

By construction, each $\mathscr D_i$ defines a right exact functor and we have a composition map
$$M \xlongrightarrow{1 \otimes \mathrm{id}} U ( \tp_i ) \otimes_{U(\tb)} M \longrightarrow \mathscr D_i ( M ).$$
We denote the resulting natural transformation by $\imath_i$. For $w\in \tSym_n$, we fix a reduced expression
$$w = s_{i_1} s_{i_2} \cdots s_{i_l} \quad i_1,\ldots,i_l \in \tI_\af \quad \text{and set}\quad \mathscr D_w := \mathscr D_{i_1} \circ \mathscr D_{i_2} \circ \cdots \circ \mathscr D_{i_l}.$$
We define
$$\imath_w := \imath_{i_1} \circ \cdots \circ \imath_{i_l} : M \longrightarrow \mathscr D_w ( M ).$$

\begin{thm}[Joseph~{\cite[2.15]{Jos85}}; see also~\cite{Kum02}]\label{thm:Jos}
For each $w \in \tSym_n$, the functor $\mathscr D_w$ does not depend on the choice of a reduced expression of $w$.
\end{thm}

By Theorem~\ref{thm:Jos}, the functor $\mathscr D_w$ and the natural transformation $\imath_w$ do not depend on the choice of the reduced expression.

\begin{thm}[Joseph~\cite{Jos85}]\label{thm:Jos2}
For each $i \in \tI_\af$, the following hold:
\begin{enumerate}
\item There is a natural isomorphism of functors $\mathscr D_i \xrightarrow{\ \cong\ } \mathscr D_i \circ \mathscr D_i$.
\item For any finite-dimensional $\tp_i$-module $M$, there is an isomorphism of functors
$$\mathscr D_i ( M \otimes \bullet ) \cong M \otimes \mathscr D_i ( \bullet ).$$
\end{enumerate}
Moreover, the functor $\mathscr D_i$ maps finite-dimensional $\tb$-modules to finite-dimensional $\tp_i$-modules, which may be regarded as $\tb$-modules via restriction.
\end{thm}

\subsection{The algebra $A$ and its representations}\label{subsec:Amod}

Fix $m>0$. Define
\[
A=A_m=\C\Sym_m\ltimes \C[X_1,\ldots,X_m]=\C\Sym_m\ltimes \C[X],
\]
where the relations are
\[
X_iX_j=X_jX_i,\qquad wX_i=X_{w(i)}w
\quad (1\le i,j\le m,\; w\in \Sym_m).
\]

By adjoining another set of variables $Y_1,\ldots,Y_m$, we define a larger algebra
\[
B=B_m=\C\Sym_m\ltimes \C[X_1,\ldots,X_m,Y_1,\ldots,Y_m]
=\C\Sym_m\ltimes \C[X,Y].
\]
We naturally have $A\subset B$, with the additional relations given by
\[
Y_iX_j=X_jY_i,\qquad Y_iY_j=Y_jY_i,\qquad wY_i=Y_{w(i)}w
\quad (1\le i,j\le m,\; w\in \Sym_m).
\]

The algebra $A$ is graded by
$$\deg w = 0 \qquad w \in \Sym_m, \quad \text{and} \quad \deg X_i = 1 \qquad (1 \le i \le m).$$

The anti-involution of $A$ defined by $X_i \mapsto X_i$ ($1 \le i \le m$) and $w \mapsto w^{-1}$ ($w \in \Sym_m$) equips $M^{\vee}$ ($M \in A \gmod$) with the structure of a graded $A$-module. Note that $M^{\vee} = M^*$ and $M^{\vee} \in A \gmod$ for $M \in A \gmod$ if and only if $\dim M < \infty$.

For each $\la =\{\la_i\}_i \in \Par_m$, we set
$$\Sym_\la := \prod_{i = 1}^{\ell(\la)} \Sym_{\la_i} \subset \Sym_m.$$

For each $\la \in \Par_m$, we have an irreducible $\Sym_m$-module $L_\la$ such that
$$L_{m} = \mathsf{triv}, \quad L_{1^m} = \mathsf{sgn}, \quad L_\la \otimes\mathsf{sgn} \cong L_{\la'}, \quad \text{and} \quad \Hom_{\Sym_\la} ( \C, L_\la ) \neq 0.$$

We may regard $L_\la$ as a simple graded $A$-module concentrated in degree zero. We have $L_\la \cong L_\la^{\vee}$ for each $\la \in \Par_m$ since $\Sym_m$ is a real reflection group. Let $e_\la \in \C \Sym_m$ be an idempotent such that $L_\la \cong \C \Sym_m e_\la$. We set
\[
P_\la := A e_\la \cong \C[X_1,\ldots,X_m] \otimes L_\la \qquad \la \in \Par_m.
\]

\begin{lem}[See, e.g.,~\cite{Kat15}]\label{lem:Asp}
The family $\{L_\la \}_{\la \in \Par_m}$ is a complete set of representatives of simple modules in $A \gmod$ up to grading shift. For each $\la \in \Par_m$, the graded $A$-module $P_\la$ is the projective cover of $L_\la$ in $A \gmod$.
\end{lem}

Note that tensoring with $\mathsf{sgn}$ over $\C$ yields an involution $M \mapsto M^{\circledast}$ on the category $A \gmod$ such that $( L_\la )^{\circledast} = L_{\la'}$ and $( P_\la )^{\circledast} = P_{\la'}$ for each $\la \in \Par_m$. We use the same notation for tensoring bigraded $B_m$-modules with $\mathsf{sgn}$ in bidegree $(0,0)$.

For each $\la \in \Par_m$, we define
\begin{equation}
\tK_\la := P_\la / \bigl( \sum_{f \in \hom_A ( P_\mu, P_\la ), \ \mu \not\ge \la} \mathrm{Im} \, f \bigr), \quad \text{and} \quad K_\la := \tK_\la / \bigl( \sum_{f \in \hom_A ( P_\la, \tK_\la )_{>0}} \mathrm{Im} \, f \bigr),\label{eqn:Kdef}
\end{equation}
where the dominance order $\le$ on $\Par_m$ is induced via the identification $\Par_m = \Par_m ( m )$.

For $M \in A \gmod$, we set
\begin{align*}
[M:L_\la]_q & := \gdim \hom_A ( P_\la, M ) &  (\la \in \Par_m), \quad \text{and} \\
\gch M & := \sum_{\mu \in \Par_m} [M:L_\mu]_q \cdot \ws_\mu & \in \bigoplus_{\mu \in \Par_m} \Q(\!(q)\!) \ws_\mu.
\end{align*}

\begin{thm}[See~{\cite[7.5.6]{MR01}}, cf.~{\cite[Theorem~1.7]{Kat22b}}]\label{thm:gdimA}
The global dimension of $A$ is finite. In particular, every $M \in A \gmod$ admits a finite projective resolution.
\end{thm}

\begin{thm}[\cite{Kat15,Kat17,Kat22b}]\label{thm:Kostka}
In the above setting, the following hold for each $\la,\mu \in \Par_m$:
\begin{enumerate}
\item $\ext^i_A ( \tK_\la, ( K_{\mu} )^{\vee} ) = \C^{\delta_{i,0}\delta_{\la,\mu}}$.
\item $\ext^\bullet_A ( \tK_\la, K_{\mu} ) = 0 = \ext^\bullet_A ( \tK_\la, L_{\mu} )$ when $\mu \not\le \la$.
\item $P_\la$ admits a finite filtration by $\{ \tK_{\mu} \}_{\mu \le \la}$.
\item $\tK_{\la}$ admits a decreasing separable filtration whose associated graded is a direct sum of grading shifts of $K_{\la}$.
\item $\mathrm{end}_{A} ( \tK_\la )$ is a polynomial ring with strictly positive degree generators.
\item If $(P_\la:\tK_\mu)_q$ denotes the graded multiplicity of $\tK_\mu$ in such a filtration of $P_\la$, then
\[
(P_\la:\tK_\mu)_q=[K_\mu:L_\la]_q=K_{\la,\mu}(q),
\]
where the rightmost term is the Kostka polynomial.
\end{enumerate}
\end{thm}

\begin{proof}
Statements~\textup{(1)}, \textup{(3)}--\textup{(6)} are consequences of Corollary~C in \cite{Kat17},
applied to the geometric realization of $A$ by the action of $G$ on the nilpotent cone of $\mathfrak g_{\circ}$.
More precisely,

\begin{itemize}
  \item[(1)] follows from the second part of~\cite[Corollary~C(2)]{Kat17};
  \item[(3)] follows from~\cite[Corollary~C(4)]{Kat17};
  \item[(4)] follows from the first part of~\cite[Corollary~C(1)]{Kat17};
  \item[(5)] follows from the second part of~\cite[Corollary~C(1)]{Kat17}
        (or~\cite[Theorem~2.3(1)]{Kat22b});
  \item[(6)] follows from the third part of~\cite[Corollary~C]{Kat17}.
\end{itemize}

The first equality in~\textup{(2)} is exactly given by~\cite[Corollary~3.9(1)]{Kat15}. 
The second equality in~\textup{(2)} follows either from
~\cite[Corollary~2.10]{Kat17} or by applying d\'evissage to the Jordan--H\"older series of $K_\mu$, inducting on the dominance order from the case $K_m=L_m$ using~\eqref{eqn:Ktri}.
\end{proof}

By Theorem~\ref{thm:Kostka}\textup{(4)(6)}, we have
\begin{equation}
[\tK_\la:L_\mu]_q = 0 = [K_\la:L_\mu]_q \qquad (\la \not\le \mu).\label{eqn:Ktri}
\end{equation}

\begin{cor}[see~{\cite[Proposition~3.3]{Kat17}}]\label{cor:Kostka}
For each $\la, \mu \in \Par_m$ such that $\mu \not\le \la$, we have
$$\ext^\bullet_A ( K_\la, K_{\mu} ) = 0 = \ext^\bullet_A ( K_\la, L_{\mu} ).$$
\end{cor}

\begin{rem}\label{rem:Atw}
We can twist the dominance order on $\Par_m$ to obtain another partial order $\le'$ such that
$$\la \le' \mu \Leftrightarrow \la' \le \mu' \Leftrightarrow \la \ge \mu \qquad (\la,\mu \in \Par_{m}).$$
We replace $\le$ with $\le'$ in the definitions of $K_\la$ and $\tK_\la$ ($\la \in \Par_m$) and denote the resulting modules by $K_\la'$ and $\tK_\la'$ respectively. We have
$$K_\la' = ( K_{\la'} )^{\circledast} \quad \text{and} \quad \tK_\la' = ( \tK_{\la'} ) ^{\circledast} \qquad (( L_{\la'} )^{\circledast} = L_\la, ( P_{\la'} )^{\circledast} = P_\la).$$
In particular, the conclusions of Theorem~\ref{thm:Kostka} and Corollary~\ref{cor:Kostka} hold for these twisted modules with respect to $\le'$.
\end{rem}

\begin{thm}[Hotta--Springer]\label{thm:GP}
For each $\la \in \Par_m$, the module $K_\la$ has a simple socle $\mathsf{q}^{\mathsf{n} ( \la )} L_{m}$ as a graded $A$-module. In addition, we have an inclusion $\mathsf{q}^{-\mathsf{n} ( \la )} K_\la \hookrightarrow \mathsf{q}^{-\mathsf{n} ( 1^m )} K_{1^m}$.
\end{thm}

\begin{proof}
By~\cite{Kat15}, the module $( K_\la )^{\vee}$ is identified with the Garsia--Procesi module (\cite{DP81,Tan82,GP92}), which in turn is realized as a quotient of the coinvariant ring by~\cite[Theorem~1.1]{HS77}. Together with $[K_\la:L_m]_q = q^{\mathsf{n}(\la)}$, which follows from Theorem~\ref{thm:Kostka}\textup{(6)}, this proves both assertions.
\end{proof}

\begin{ex}[$m=3$ case taken from~\cite{Kat17}]
We have
\begin{align*}
\gch K_{1^3} & = \ws_{1^3} + (q+q^2) \ws_{21} + q^3 \ws_3, \quad \gch K_{21} = \ws_{21} + q \ws_3, \quad  \gch K_{3} = \ws_3,\\
\gch \tK_{1^3} & = \frac{\gch K_{1^3}}{(1-q)(1-q^2)(1-q^3)}, \quad \gch \tK_{21} = \frac{\gch K_{21}}{(1-q)^2}, \quad \text{and} \quad \gch \tK_{3} = \frac{\gch K_{3}}{(1-q)}.
\end{align*}
Here the fractions are interpreted as formal power series in $q$.
\end{ex}

\subsection{Representation theory of current algebras}

Let $\g[z] \mathchar`-\mathsf{gMod}$ denote the category of graded $\g[z]$-modules $M$ such that
$$\dim M_i < \infty \qquad \text{for all} \quad i \in \Z,$$
and $M$ decomposes as a $\g$-module into a direct sum of $\{V_{\la}^*\}_{\la \in \sP^+}$. For each $\mu \in \sP^+$, we define $\g[z] \gmod^{\le \mu}$ to be the full subcategory of $\g[z] \mathchar`-\mathsf{gMod}$ whose modules are finitely generated as $U(\g[z])$-modules and, as $\g$-modules, are isomorphic to direct sums of $\{V_{\la}^*\}_{\mu \ge \la \in \sP^+}$. We also set
$$\g[z] \gmod_{m}:= \g[z] \gmod^{\le m \epsilon_1} \quad \text{and} \quad \g[z] \gmod := \bigoplus_{m \ge 0} \g[z] \gmod_{m} \subset \g[z] \mathchar`-\mathsf{gMod}.$$
These categories are full subcategories of $\g[z]\mathchar`-\mathsf{gMod}$ that inherit the endofunctors $\mathsf q^{\pm 1}$.

The module $V_\la^*$ can be regarded as a graded simple module of $\g [z]$ concentrated in degree $0$, via the evaluation map
$$\mathtt{ev}_0 : \g [z] \longrightarrow \g \qquad z \mapsto 0.$$

For a graded $\g[z]$-module $M$, we set
$$[M:V_\la^*]_q:= \gdim \hom_{\g} ( V_\la^*, M ) = \sum_{i \in \Z} q^i \dim \hom_{\g} ( V_\la^*, M_i ) \qquad (\la \in \sP^+).$$
We also define the graded characters of $M \in \g[z] \mathchar`-\mathsf{gMod}$ and $M' \in \g[z]\gmod_m$ as:
\begin{align*}
\ugch M & = \sum_{\la \in \sP^+} [M:V_{\la}^*]_q \cdot s_\la\\
\gch M' & = \sum_{\la \in \Par(n)} [M':V_{\la}^*]_q \cdot \ws_\la \in {\bm \La}_q.
\end{align*}
The first expression is interpreted as a formal sum, degree by degree.
By definition, we have $\mathsf{tr}_n ( \gch M' ) = \ugch M'$.
These record the graded multiplicities of $V_{\la}^*$. We also set
$$[M:V_{\la}^*] := \Bigl.[M:V_{\la}^*]_q \Bigr| _{q=1} \qquad (\la \in \sP^+)$$
when the right-hand side defines an integer. Otherwise, we set $[M:V_{\la}^*] := \infty$.

Although not explicitly stated in this exact form, the following result is a direct consequence of the cited works:

\begin{thm}[Chari--Greenstein~\cite{CG07}, Kleshchev~\cite{Kle14}]
Let $\mu \in \sP^+$. For each $\la \in \sP^+$ with $\la \le \mu$, we have the projective cover $\P _\la^{\le \mu}$ of $V_{\la}^*$ in $\g[z] \gmod^{\le \mu}$ and the projective cover $\P _\la$ of $V_{\la}^*$ in $\g[z] \mathchar`-\mathsf{gMod}$.
\end{thm}

For $\la \in \sP^+$, we set $\bW_\la := \P _\la^{\le \la}$.

\begin{thm}[Chari--Loktev~\cite{CL06}, Feigin--Khoroshkin--Makedonskyi~\cite{FKM19}]
For $\la \in \Par(n)$, we have
$$\mathrm{End}_{\g[z]} ( \bW_\la ) \cong \bigotimes_{i=1}^{n-1} \C [x_{i,1},\ldots,x_{i,m_i}]^{\Sym_{m_i}},$$
where $m_i = \la_i - \la_{i+1}$ for $1 \le i < n$. Moreover, $\bW_\la$ is free over this endomorphism ring. Specializing all $x$-variables in $\mathrm{End}_{\g[z]} ( \bW_\la )$ to $0$ yields a finite-dimensional graded $\g[z]$-module $W_\la$.
\end{thm}

\begin{ex}
When $n \ge 3$, we have
$$\gch W_{3} = \ws_3 + (q + q^2) \ws_{21} + q^3 \ws_{1^3}, \quad \gch W_{21} = \ws_{21} + q \ws_{1^3}, \quad \gch W_{1^3} = \ws_{1^3}.$$
\end{ex}

\begin{thm}[Feigin--Khoroshkin--Makedonskyi~\cite{FKM19}, Flicker~\cite{Fl21}]\label{thm:FKM}
Assume $n > m$. We have
$$\bW_{\varpi_1} \cong V_{1}^* \otimes \C [z]= (\C^n)^* \otimes \C [z],$$
whose $\g[z]$-action is the scalar extension of the $\g$-module $V_1^* = V_{\varpi_1}^*$. We have
\begin{equation}
A = A_m \cong \mathrm{End}_{\g [z]} ( \bW_{\varpi_1} ^{\otimes m} ),	\label{eqn:progen}
\end{equation}
where $X_i$ corresponds to the action of $z$ on the $i$-th tensor factor of $\bW_{\varpi_1}^{\otimes m}$.

In addition, $\bW_{\varpi_1}^{\otimes m}$ is a projective $A$-module. This establishes mutually quasi-inverse equivalences
\begin{align*}
\WS : \, & A \gmod \ni M \mapsto ( \bW_{\varpi_1} ^{\otimes m} ) \otimes_{A} M \in \g[z] \gmod_m,\quad\text{and}\\
\SW : \, & \g[z] \gmod_m \ni M \mapsto \hom_{U(\g [z])}( \bW_{\varpi_1} ^{\otimes m}, M ) \in A \gmod
\end{align*}
such that
$$\WS ( L_\la ) \cong V_\la^* \qquad (\la \in \Par_m).$$
\end{thm}

\begin{rem}[twisted Schur--Weyl functor]\label{rem:twSW}
We can also consider the twisted version $\WS^{*}$ of $\WS$, which is the original formulation in~\cite{FKM19,Fl21} and uses
$$\bW_{-\epsilon_n} \cong \C^n \otimes \C [z]$$
in place of $\bW_{\varpi_1}$.
\end{rem}

\subsection{A geometric realization of Catalan functions}

We recall the following results from~\cite{Kat23a}:

\begin{defn}[Root ideal]\label{defn:ri}
A subset $\Psi \subset \Delta^+$ is called a root ideal (of $\g$) if and only if
 $$( \Psi + \Delta^+ )\cap \Delta^+ \subset \Psi.$$
A corner of a root ideal $\Psi$ is a pair $(i,j) \in \{1,2,\ldots,n\}^2$ (or $\al_{ij}$) such that
$$\al_{ij} \in \Psi, \quad \text{and} \quad \{ \al_{i+1,j},\al_{i,j-1} \} \cap \Psi = \emptyset.$$
Here, the latter condition is imposed only when the corresponding roots are well-defined (i.e., when $i < n$ or $j > 1$). Note that for a root ideal $\Psi$, a corner $(i,j)$ is uniquely determined by its row index $i$ (resp. column index $j$). We thus refer to such a corner as the corner of the form $(i,\bullet)$ (resp. $(\bullet,j)$).

Removing a corner $\al_{ij}$ from $\Psi$ defines a root ideal that we denote by $\Psi\setminus(i,j)$. We also write $\Psi = ( \Psi\setminus(i,j) ) \cup (i,j)$.
\end{defn}

\begin{thm}[Cellini~\cite{Cel00}]
The set of $B$-submodules of $\gn$ is in bijection with the set of root ideals of $\Delta^+$.
\end{thm}

We define the $B$-submodule $\gn ( \Psi ) \subset \gn$ corresponding to the root ideal $\Psi$ as
$$\gn (\Psi ) := \bigoplus_{\al_{ij} \in \Psi} \C E_{ij} \subset \gn.$$
By abuse of notation, we sometimes identify $\Psi$ with $\gn ( \Psi )$. Under this identification, the containment relation of root ideals corresponds to the inclusion of subspaces. That is, for root ideals $\Psi,\Psi'$, we have $\Psi \subset \Psi'$ if and only if $\gn ( \Psi ) \subset \gn ( \Psi' )$. This defines a partial order $\subset$ on the set of root ideals. We set $\Psi_0 := \emptyset \subset \Delta^+$. We define
$$T^*_\Psi X := G \times^B \gn( \Psi) \subset G \times^B \gn = T^* X.$$
We have $X = T^*_{\Psi_0} X \subset T^*_\Psi X \subset T^*_{\Delta^+} X = T^* X$.

\begin{thm}[{\cite[Theorem~A]{Kat23a}}]\label{thm:Cat}
Let $\Psi$ be a root ideal. There exists a connected smooth projective algebraic variety $\sX_\Psi$ equipped with an action of $\widehat{G}[\![z]\!]$ with the following properties:
\begin{enumerate}
\item We have an inclusion $T^*_\Psi X \subset \sX_\Psi$, whose image is a $(\Gm \times G)$-stable Zariski open dense subset.
\item $T^*_\Psi X$ is the attracting set of $X \subset \sX_\Psi$ with respect to the $\Gm$-action.
\item For each $\la \in \sP$, there exists a $\widehat{G}[\![z]\!]$-equivariant line bundle $\cO_{\sX_\Psi} ( \la )$ such that its restriction to $T^*_\Psi X$ is $\cO_{T^*_\Psi X} ( \la )$.
\item For each $\la \in \Par(n)$, the higher cohomology vanishes
\[
H^{>0} ( \sX_\Psi, \cO_{\sX_\Psi} ( \la )) = 0\quad\text{and}\quad H^{>0} ( T^*_\Psi X, \cO_{T^*_\Psi X} ( \la ))=0.
\]
Moreover, 
\begin{equation}
 [H^0 ( \sX_\Psi, \cO_{\sX_\Psi} ( \la )):V_{\mu}^*]_q = [H^0 ( T^*_\Psi X, \cO_{T^*_\Psi X} ( \la )):V_{\mu}^*]_q \qquad (\mu \in \Par(n)),\label{eqn:HPtrunc}
\end{equation}
and the left-hand side is zero for $\mu \in \sP^+ \setminus \Par(n)$.
\end{enumerate}
\end{thm}

Note that the vanishing on the open total space in Theorem~\ref{thm:Cat}\textup{(4)} holds for every $\la \in \sP^+$. Indeed, for a fixed root ideal $\Phi$, choose $N\ge0$ such that $\la+N\varpi_n\in\Par(n)$. Since the line bundle of weight $\varpi_n$ on $T^*_\Phi X$ is a determinant character twist, we have
\[
H^\bullet\bigl(T^*_\Phi X,\cO_{T^*_\Phi X}(\la)\bigr)
\cong H^\bullet\bigl(T^*_\Phi X,\cO_{T^*_\Phi X}(\la+N\varpi_n)\bigr)\otimes\C_{N\varpi_n}.
\]
We use this extension of the vanishing statement in the proof of Corollary~\ref{cor:kPieri} and Theorem~\ref{thm:Fdp}.

\begin{thm}[\cite{Kat23a} Theorem~C]\label{thm:Cat-sh}
Let $\Psi$ be a root ideal. For each $\la \in \Par(n)$, differentiating the action of $\widehat{G}[\![z]\!]$ induces a graded $\g[z]$-module structure on $H^0 ( \sX_\Psi, \cO_{\sX_\Psi} ( \la ))$ which has a simple head $V_{\la}^*$.
\end{thm}

The following records special cases of~\cite[Proposition~5.12 and Corollary~5.13]{Kat23a} that follow from Theorem~\ref{thm:Cat-sh}:

\begin{cor}\label{cor:Cat-sg}
Let $\Psi$ be a root ideal. For each $\la \in \Par(n)$, the restriction induces surjections of graded $\g[z]$-modules
\[
H^0 ( \sX_{\Delta^+}, \cO_{\sX_{\Delta^+}} ( \la )) \twoheadrightarrow H^0 ( \sX_\Psi, \cO_{\sX_\Psi} ( \la )), \quad \text{and}\quad
H^0 ( T^* X, \cO_{T^* X} ( \la )) \twoheadrightarrow H^0 ( T^*_{\Psi} X, \cO_{T^*_{\Psi} X} ( \la )).
\]
\end{cor}

\begin{rem}
The variety $\sX_\Psi$ in~\cite[Introduction]{Kat23a} refers to the variety $\sX_\Psi (w_0)$ in the main body of~\cite{Kat23a}. Since we are mainly interested in varieties equipped with the $G$-action, we adopt the shorthand $\sX_\Psi$ here.
\end{rem}

\begin{cor}\label{cor:Cat-fd}
Let $\Psi$ be a root ideal, and let $\la \in \sP$. We have
$$\dim H^0 ( \sX_\Psi, \cO_{\sX_\Psi} ( \la )) < \infty.$$
When $\la \in \Par_m(n)$, we have
\[
H^0 ( \sX_\Psi, \cO_{\sX_\Psi} ( \la )) \in \g[z]\gmod_m.
\]
\end{cor}
\begin{proof}
The first assertion follows from the projectivity of $\sX_\Psi$ together with~\cite[I\!I\!I Theorem~8.8]{Har77}. The second assertion follows from~\eqref{eqn:HPtrunc} and the fact that the action of $\mathrm{Id}$ is constant on $H^0 ( T^*_\Psi X, \cO_{T^*_\Psi X} ( \la ))$.
\end{proof}

\begin{cor}[{\cite[Corollary~4.7 $+$ Corollary~5.14]{Kat23a}}]\label{cor:inclCat}
Let $\Psi \subset \Psi'$ be two root ideals. We have the following commutative diagram of embeddings:
\begin{equation}
\begin{tikzcd}
    \sX_{\Psi} \arrow[r, hook] & \sX_{\Psi'} \\
    T^*_{\Psi} X \arrow[r, hook] \arrow[u, hook] & T^*_{\Psi'} X \arrow[u, hook]
\end{tikzcd}
\label{eqn:cemb}.
\end{equation}
In addition, the restriction of a line bundle $\cO_{\sX_{\Psi'}} ( \la )$ yields the line bundle with the same label for each $\la \in \Par(n)$ and each variety in the diagram~\eqref{eqn:cemb}.
\end{cor}

\subsection{Level $k$ affine Demazure modules}\label{subsec:level-k}
We fix $k \in \Z_{>0}$ and $n  > m \in \Z_{> 0}$.

\begin{lem}
For $\la\in\Par_m(n)$, we have
\[
\la+k\wp\in\sP_{\af}^{+}
\quad\Longleftrightarrow\quad
\la\in\Par_m^{(k)}(n).
\]
\end{lem}

\begin{proof}
The assumption $n > m$ implies $\la_n=0$. We have
\[
\la + k \wp = (k-\la_1) \wp + \sum_{i=1}^{n-1}(\la_i-\la_{i+1})\La_i.
\]
Because $\la$ is a partition, the coefficients
$\la_i-\la_{i+1}$ are nonnegative. Hence
$\la+k\wp$ is affine dominant if and only if
$k-\la_1\geq 0$, or equivalently, if and only if
$\la\in\Par_m^{(k)}(n)$.
\end{proof}

By~\cite[Corollary~10.1]{Kac}, for each $\la \in \sP^+$, there is a unique weight $[\la]_k \in \sP$ and a unique integer $e_k ( \la )$ such that
\begin{equation}
w ( w_0 \la + k \wp ) = [\la]_k + k \wp + e_k ( \la ) \delta \in \sP^+_\af\label{eqn:afdom}
\end{equation}
for some $w \in \tSym_n$. We therefore define the level $k$-Demazure module associated with $\la \in \Par_m (n)$ by
\[
W^{(k)}_{\la} := \mathscr D_{w^{-1}} ( \C_{[\la]_k + k \wp + e_k ( \la ) \delta} ).
\]
The module $W^{(k)}_{\la}$ is independent of the choice of $w$ satisfying~\eqref{eqn:afdom} by Theorem~\ref{thm:Jos2}. Since $w_0 \la$ is anti-dominant, we deduce that $W^{(k)}_\la$ is a graded $\g[z]$-module (see~\cite[Lemma~2.7]{Kat26b}).

The following proposition differs slightly from the original statement because of our character convention:

\begin{prop}[{\cite[Proposition~1.20]{Kat23a}}]\label{prop:D-dual}
If $n > m \in \Z_{>0}$ and $\la \in \Par_m (n)$, then $( W^{(k)}_\la )^*\in \g[z]\gmod_m$.
\end{prop}

\begin{thm}\label{thm:D-filt}
Let $\Psi$ be a root ideal. For each $\la \in \Par^{(k)}_m (n)$, the module
\[
H^0 ( \sX_\Psi, \cO_{\sX_\Psi} ( \la ))
\]
admits a finite filtration by $\{( W^{(k)}_\mu )^*\}_{\mu \in \Par_m(n)}$.
\end{thm}

\begin{proof}
By~\cite[Theorem~4.1(3)]{Kat23a},
$H^0(\sX_\Psi,\cO_{\sX_\Psi}(\la))$ admits a finite filtration by
the duals of the $\g$-stable affine Demazure modules of level
$\la_1$. These factors are of the form $(W^{(\la_1)}_\mu)^*$ ($\mu\in \sP^+$).

Since $k\geq\la_1$,~\cite[Theorem~A(4)]{Kat26b} implies that each
$(W^{(\la_1)}_\mu)^*$ admits a finite filtration by modules of the
form $\mathsf{q}^a (W^{(k)}_\nu)^*$ ($a \in \Z, \nu \in \sP^+$). Refining the original filtration therefore
gives a finite filtration of
$H^0(\sX_\Psi,\cO_{\sX_\Psi}(\la))$ by such modules.

Finally, every filtration factor satisfies the weight condition defining
$\g[z]\gmod_m$. Hence only indices
$\nu\in\Par_m(n)$ can occur.
\end{proof}

The following result is a slight variant of~\cite[Lemma~5.5]{Kat26b}:

\begin{prop}\label{prop:D-ker}
Suppose that $M, N \in \g[z]\gmod_m$ are two modules which admit finite filtrations by $\{( W^{(k)}_\la )^*\}_{\la \in \Par_m(n)}$. For an injection $f : N \hookrightarrow M$ of graded $\g[z]$-modules, $\coker f$ admits a finite filtration by $\{( W^{(k)}_\la )^*\}_{\la \in \Par_m(n)}$.
\end{prop}

\begin{proof}
By~\cite[Lemma~5.5]{Kat26b} applied to $W^{(k)}$, we find that the kernel of the dual map
\[
f^{*} : M^* \twoheadrightarrow N^*
\]
admits a finite filtration by $\{W^{(k)}_\la\}_{\la \in \sP^+}$. Since every weight of $\ker f^{*}$ belongs to $\Comp_m(n) \subset \sP$, every filtration factor $W^{(k)}_\la$ which can contribute to $\ker f^{*}$ satisfies $\la \in \Par_m(n)$.

Since $\coker f \cong (\ker f^*)^*$, dualizing this filtration yields the assertion.
\end{proof}

\section{Module-theoretic preparations}\label{sec:mtp}

Keep the setting of the previous sections.

\subsection{Kostka systems and Hall--Littlewood modules}\label{subsec:KSHL}

We define the modules
$$\HL _\la := H^0 ( \sX_{\Delta^+}, \cO_{\sX_{\Delta^+}} ( \la )) \in \g[z]\gmod_m \qquad (\la \in \Par_m (n)).$$

\begin{prop}\label{prop:Kostrans}
Assume $n > m$. For each $\la \in \Par_m$, we have
$$\HL_\la \cong \WS ( K_{\la} ).$$
In particular, $\HL_\la$ is the largest module in $\g[z]\gmod_m$ whose head is $V_\la^*$ and satisfies
\[
[ \HL_\la : V_{\mu}^* ]_q =
\begin{cases}
1 & (\la = \mu)\\
0 & (\la \not\le \mu)
\end{cases}.
\]
\end{prop}

\begin{proof}
By~\eqref{eqn:HPtrunc}, together with~\cite{Gup87} and Theorem~\ref{thm:Kostka}\textup{(6)}, we obtain
\begin{equation}
[\HL_\la : V_{\mu}^*]_q = K_{\mu,\la}(q) = [K_{\la} : L_{\mu}]_q, \qquad (\la,\mu \in \Par_m),
\label{eqn:stdmult}
\end{equation}
where $K_{\mu,\la}(q)$ denotes the Kostka polynomial. (We note that the first equation holds for $\la,\mu \in \Par_m (n)$ for every $n$.) 
Recall from~\eqref{eqn:Kdef} that $K_\la$ is the maximal object in $A_m \gmod$ generated by $L_\la$ subject to the multiplicity condition
$$
[K_\la:L_{\mu}]_q =
\begin{cases}
1 & (\la = \mu)\\
0 & (\la \not\le \mu)
\end{cases}.
$$
Applying the category equivalence $\WS$, we see that $\WS ( K_\la )$ is the maximal object in $\g [z]\gmod_m$ generated by $V_\la^*$ such that
\[
[\WS ( K_\la ) : V_{\mu}^* ]_q =
\begin{cases}
1 & (\la = \mu)\\
0 & (\la \not\le \mu)
\end{cases}.
\]
The module $\HL_\la$ contains $V_\la^*$ with multiplicity one, and for any $\mu$ with $\la \not\le \mu$, the multiplicity of $V_\mu^*$ is zero. Therefore, the maximality property of $\WS( K_\la )$ read off from the definition~\eqref{eqn:Kdef} yields a nonzero morphism
\[
\WS ( K_\la ) \longrightarrow \HL_\la
\]
in $\g[z]\gmod_m$ that induces the identity on the head $V_\la^*$.  
By Theorem~\ref{thm:Cat-sh}, the module $\HL_\la$ has a simple head. 
It follows that the morphism $\WS ( K_\la ) \to \HL_\la$ is necessarily surjective in $\g [z]\gmod_m$.  
Comparing graded characters in~\eqref{eqn:stdmult}, we conclude that this surjection must be an isomorphism.

This proves both assertions.
\end{proof}

For each root ideal $\Psi$ and $\la \in \Par_m(n)$, we have a graded $\g[z]$-module
$$\HL_\la^{\Psi} := H^{0} ( \sX_\Psi, \cO_{\sX_\Psi} ( \la )) \in \g[z]\gmod_m$$
afforded by Theorem~\ref{thm:Cat-sh}.

\subsection{Structural results on some filtrations}

\begin{lem}\label{lem:Ksplit}
Let $\la, \mu \in \Par_m$ and $d \in \Z$. If $\la \not> \mu$, then every short exact sequence of graded $A_m$-modules
$$0 \rightarrow \mathsf q^d ( \tK_{\mu} ) \rightarrow E \rightarrow \tK_\la \rightarrow 0$$
splits.
\end{lem}

\begin{proof}
In view of Theorem~\ref{thm:Kostka}\textup{(4)}, this follows from the $\ext^1$-part of Theorem~\ref{thm:Kostka}\textup{(2)} when $\la \neq \mu$, and the fact that $\tK_\la$ is projective in the full subcategory of modules $N\in A_m\gmod$ satisfying $[N:L_\nu]_q=0$ for $\nu\not\ge\la$. Indeed, by~\eqref{eqn:Kdef}, every map $P_\la\to N$ factors through $\tK_\la$, so $\hom_{A_m}(\tK_\la,N)\cong e_\la N$ is exact on this subcategory.
\end{proof}

\begin{cor}\label{cor:filtproj2}
Let $\la \in \Par_m$. We enumerate the set
$$\{ \mu \in \Par_m \mid \mu \le \la \} = \{\la^{(1)},\la^{(2)},\ldots,\la^{(l)}\}$$
such that $\la^{(i)} < \la^{(j)}$ in the dominance order implies $i < j$.
Then there exists a finite filtration
$$0 = F_0 \subset F_1 \subset F_2 \subset \cdots \subset F_l = P_\la$$
of graded $A_m$-modules such that each subquotient $F_{i} / F_{i-1}$ $(0 <i \le l)$ is isomorphic to a direct sum of grading shifts of $\tK_{\la^{(i)}}$.
\end{cor}

\begin{proof}
By Theorem~\ref{thm:Kostka}\textup{(3)}, there exists a finite filtration
$$0 = G_0 \subsetneq G_1 \subsetneq G_2 \subsetneq \cdots \subsetneq G_r = P_\la$$
of $P_\la$ such that each subquotient $G_{i}/G_{i-1}$ $(1 \le i \le r)$ is a direct sum of grading shifts of $\tK_{\la^{(s(i))}}$ for some $1 \le s(i) \le l$. Consider the short exact sequence
$$0 \longrightarrow G_{i} / G_{i-1} \longrightarrow G_{i+1} / G_{i-1} \longrightarrow G_{i+1}/G_{i} \longrightarrow 0.$$
Whenever $s(i) \ge s(i+1)$, the condition $\la^{(s(i+1))} \not> \la^{(s(i))}$ holds. Thus, by Lemma~\ref{lem:Ksplit}, the sequence splits, and we have
\[G_{i+1} / G_{i-1} \cong G_{i+1}/G_i \oplus G_{i}/G_{i-1}.\]
This allows us to modify the filtration: if $s(i) > s(i+1)$, we can swap the order of factors by replacing $G_i$ with
$$\ker ( G_{i+1} \longrightarrow G_i/G_{i-1} ).$$
If $s(i)=s(i+1)$, we can merge the steps by removing $G_i$ from the filtration. In either case, we obtain a new filtration satisfying the same condition as the original one. Applying this procedure strictly decreases either the number of inversions (pairs $i<j$ with $s(i)>s(j)$) or the length $r$ of the filtration. Therefore, the process terminates, and we obtain the desired filtration (possibly by inserting trivial steps to match the length $l$).
\end{proof}

\section{Analysis of the socle}\label{sec:as}

Keep the setting of the previous section, particularly that of \S\ref{subsec:AffDem}.

For $\tJ \subsetneq \tI_\af$, we have a Lie subalgebra $\mathfrak{l}_{\tJ}$ of $\tg$ generated by $E_i,F_i,[E_i,F_i]$ for $i \in \tJ$, and the Lie subalgebra $\tp_{\tJ}$ of $\tg$ generated by $\tb$ and $\{F_i\}_{i\in\tJ}$. We have a decomposition
\[
\tp_\tJ = \mathfrak l_\tJ \oplus \wgt_\tJ^{\perp} \oplus \gu_{\tJ},
\]
where $\wgt_\tJ^{\perp} \subset \wgt$ is an abelian Lie subalgebra and $\gu_\tJ$ is a Lie subalgebra of $[\tb,\tb]$ normalized by the action of $( \mathfrak{l}_\tJ \oplus \wgt_{\tJ}^{\perp} )$. We also set $\Sym_{\tJ} := \left< s_j \mid j \in \tJ \right> \subset \tSym_n$. Let $w_0^{\tJ} \in \Sym_{\tJ}$ be its longest element. If $\tJ$ defines a connected subdiagram of the Dynkin diagram of $\tg_{\circ}$, then $\Sym_{\tJ} \cong \Sym_{|\tJ|+1}$ as Coxeter groups.

As explained in~\cite[\S1.1]{Kat23a}, we have two subgroups $L_{\tJ} \subset \widehat{P}_{\tJ} \subset \widehat{G}(\!(z)\!)$ whose Lie algebras are $\mathfrak{l}_{\tJ}$ and $\tp_{\tJ} / \C K$, respectively. We set $B_{\tJ} := ( L_{\tJ} \cap \widehat{P}_{\emptyset} )$, which is a Borel subgroup of the semisimple group $L_{\tJ}$. We define the flag manifold of $\tJ$ by
$$X_{\tJ} := \widehat{P}_{\tJ} / \widehat{P}_{\emptyset} \cong L_{\tJ} / B_{\tJ}.$$
This embeds into the flag variety $\mathbf{X}$ of $\tg_{\circ}$ as the $\widehat{P}_{\tJ}$-orbit of its unique $\widehat{P}_{\emptyset}$-fixed point in $\mathbf{X}$, namely, the Schubert variety of $\mathbf{X}$ corresponding to $w_0^{\tJ} \in \tSym_n$.

\begin{thm}[Demazure]\label{thm:kum}
Let $\tJ \subsetneq \tI_\af$ be a subset. Then $\mathscr D_{w_0^{\tJ}} ( M )$ is a $\tp_{\tJ}$-module for a $\wgt$-semisimple $\tb$-module $M$.
\end{thm}

\begin{proof}
The functor $\mathscr D_{w_0^{\tJ}}$ coincides with the (dual form of the) induction defined using the Bott--Samelson--Demazure--Hansen resolution~\cite[8.1.17 Proposition]{Kum02}. It follows that this induction, together with the Bott--Samelson--Demazure--Hansen resolution used to define it, is the same whether one starts from a $\tb$-module or from a $( \tb \cap \mathfrak{l}_{\tJ})$-module. In view of~\cite[8.1.26]{Kum02}, $\mathscr D_{w_0^{\tJ}}$ is isomorphic to the induction functor defined by $X_{\tJ} \subset \mathbf{X}$. Thus, in addition to the $\tb$-action built into the construction of Demazure functors, it is stable under the action of $\mathfrak{l}_{\tJ}$. For an algebraic explanation, see, e.g.,~\cite[Corollary~1.15]{Kat26b}.
\end{proof}

\begin{cor}\label{cor:Henv}
Keep the setting of Theorem~\ref{thm:kum}. We have a quotient map
$$U ( \tp_{\tJ} ) \otimes_{U ( \tb )} M \longrightarrow \mathscr D_{w_0^{\tJ}} ( M ).$$
\end{cor}

\begin{proof}
Set $D:=\mathscr D_{w_0^{\tJ}}(M)$. By Theorem~\ref{thm:kum}
and $\mathscr D_iD\cong D$ for $i\in\tJ$, the module $D$ is
$\mathfrak l_{\tJ}$-integrable.
For every $\mathfrak l_{\tJ}$-integrable $\tp_{\tJ}$-module $N$,
rank-one adjunction, iterated along a reduced expression, gives
\[
 \Hom_{\tb}(D,N)\cong\Hom_{\tb}(M,N).
\]
Indeed, $\mathscr D_i$ is the maximal integrable quotient of
induction, and a $\tb$-linear map between
$\mathfrak{sl}(2,i)$-integrable modules is automatically $\tp_i$-linear.
The displayed isomorphism is induced by the natural map $M\to D$;
this map need not be injective or nonzero.
Let $D'$ be the $\tp_{\tJ}$-submodule generated by its image.
The quotient $N=D/D'$ is again $\mathfrak l_{\tJ}$-integrable.
The quotient map $D\to N$ is zero on $M$, so the displayed
isomorphism forces it to be zero. Thus $D'=D$, proving the
asserted surjection, including the case $D=0$.
\end{proof}

\begin{defn}[Cyclic interval of $\tI_\af$ and associated Demazure functors]
A cyclic interval $[i,j]$ of $\tI_\af$ is a proper subset
$$\{i,(i+1),\ldots,(n-1),0,\ldots,j\}$$
consisting of a sequence of consecutive elements in $\tI_\af$ whose indices are counted modulo $n$. We also set $(i,j] := [i+1,j]$ and $[i,j) := [i,j-1]$, where we understand $(i,j] = \emptyset = [i,j)$ whenever $i \equiv j \mod n$. For a cyclic interval $[i,j] \subsetneq \tI_\af$, we define
\begin{align*}
\mathscr D_{[i,j]} & := \mathscr D_{j} \circ \cdots \circ \mathscr D_{0} \circ \mathscr D_{n-1} \circ \cdots \circ \mathscr D_{i+1} \circ \mathscr D_i,\\
\mathscr D_{(i,j]} & := \mathscr D_{j} \circ \cdots \circ \mathscr D_{0} \circ \mathscr D_{n-1} \circ \cdots \circ \mathscr D_{i+2} \circ \mathscr D_{i+1}, \quad \text{and}\\
\mathscr D_{[i,j)} & := \mathscr D_{j-1} \circ \cdots \circ \mathscr D_{0} \circ \mathscr D_{n-1} \circ \cdots \circ \mathscr D_{i+1} \circ \mathscr D_i.
\end{align*}
We have a natural transformation
$$\imath_{[i,j]} := \imath_{j} \circ \cdots \circ \imath_0 \circ \imath_{n-1} \circ \cdots \circ \imath_{i+1} \circ \imath_{i} : \bullet \longrightarrow \mathscr D_{[i,j]} ( \bullet )$$
in the category of $\wgt$-semisimple $\tb$-modules.
\end{defn}

\begin{lem}\label{lem:Levi}
Let $[i,j] \subsetneq \tI_\af$ be a cyclic interval. For a finite-dimensional $\wgt$-semisimple $\tp_{(i,j]}$-module $M$,
\begin{equation}
\mathscr D_{[i,j]} ( M ) \cong \mathscr D_{w_0^{[i,j]}} (M) \label{eqn:inflM}
\end{equation}
is a finite-dimensional $\wgt$-semisimple $\tp_{[i,j]}$-module.
\end{lem}

\begin{proof}
Note that $M$ is a finite-dimensional $\mathfrak{sl} (2,l)$-module for each $l \in (i,j]$ by restriction. It follows that $\mathscr D_l ( M )\cong M$ for $l \in (i,j]$ by Theorem~\ref{thm:Jos2}\textup{(2)}. Note that $s_{j} s_{j-1} \cdots s_{i+1} s_i$ extends to the longest element $w_0^{[i,j]}$ of $\Sym_{[i,j]}$ by appending
$$( s_{j} \cdots s_{i+1} )( s_{j} \cdots s_{i+2} ) \cdots (s_{j}s_{j-1})s_{j},$$
and this latter element belongs to $\Sym_{(i,j]}$. Thus, the left-hand side of~\eqref{eqn:inflM} coincides with $\mathscr D_{w_0^{[i,j]}} (M)$. Therefore, Theorem~\ref{thm:kum} equips~\eqref{eqn:inflM} with the structure of a $\tp_{[i,j]}$-module that extends the $\tb$-module structure.
\end{proof}

The next proposition is the basic device that allows us to propagate cyclic generators through a Demazure step.

\begin{prop}\label{prop:gen}
Let $[i,j] \subsetneq \tI_\af$ be a cyclic interval. Assume that a finite-dimensional $\wgt$-semisimple $\tp_{(i,j]}$-module $M$ satisfies the following properties:
\begin{itemize}
\item $M$ is generated by a $\wgt$-eigenvector $v$ of $\wgt$-weight $\La$ as a $\tb$-module;
\item the natural map $\imath_{[i,j]} : M \to \mathscr D_{[i,j]} ( M )$ is injective;
\item the span $M^{\top} := U ( \mathfrak l_{(i,j]} ) v \subset M$ is an irreducible $\mathfrak l_{(i,j]}$-module;
\item we have $M^{\top} \subset \mathscr D_{[i,j]} ( M^{\top} )$ as $(\tp_{(i,j]} \cap \mathfrak l_{[i,j]})$-modules.
\end{itemize}
Then $\mathscr D_{[i,j]} ( M )$ is generated by $U ( \mathfrak l_{[i,j]} ) v$ as a $\tb$-module, and the latter span is isomorphic to $\mathscr D_{[i,j]} ( M^{\top} )$.
\end{prop}

\begin{proof}
We set $\tJ :=[i,j]$ and $\tJ' := (i,j]$. We first show that $v$ generates $\mathscr D_{\tJ}(M)$ as a $\tp_{\tJ}$-module.

Note that the $\wgt$-weight $\La$-part of $M$ is one-dimensional since the $\tb$-action preserves $\C v$ or raises its weight with respect to $\lhd$. Using the commutative diagram
\[
\begin{tikzcd}
    & U(\mathfrak{l}_{\tJ}) \otimes_{U(\gb_{\tJ})} M \arrow[dr, two heads] & \\
    M \arrow[ru, hook] \arrow[r, hook] & \mathscr{D}_{\tJ}(M) \arrow[r, equal] & \mathscr{D}_{w_0^{\tJ}}(M)
\end{tikzcd}
\]
provided by~\eqref{eqn:inflM} and Corollary~\ref{cor:Henv}, we find that $M$ generates $\mathscr D_{\tJ} ( M )$ as a $\tp_{\tJ}$-module. Since the category of finite-dimensional $\gt$-semisimple $\mathfrak{l}_{\tJ}$-modules is semisimple, if we view $\mathscr D_{\tJ} ( M )$ as a direct sum of irreducible $\mathfrak{l}_{\tJ}$-modules, then their extensions as $\tp_{\tJ}$-modules are realized by the action of $\gu_{\tJ}$. Here $v$ generates $M$ by the $\tb$-action. It follows that $v$ generates $\mathscr D_{\tJ} ( M )$ as a $\tp_{\tJ}$-module. Thus $v \in M^\top$ generates $\mathscr D_{\tJ} ( M )$ as a $\tp_{\tJ}$-module.

Finally, we reduce the generation statement from $\tp_{\tJ}$ to $\tb$. In view of Lemma~\ref{lem:Levi} and the assumptions, $M^\top$ is the $\mathfrak l_{\tJ'}$-module induced from its highest weight line. This line is the $(\tb \cap \mathfrak{l}_{\tJ})$-eigenspace $\C_{\La'}$. By the Demazure character formula~\cite[5.6]{Jos85}, $\mathscr D_{\tJ} ( M^{\top} )$ is nonzero if and only if $\La'$ defines a dominant $\mathfrak{l}_\tJ$-weight. In this case, $\mathscr D_{\tJ} ( M^{\top} )$ is an irreducible $\mathfrak{l}_{\tJ}$-module with extremal weight $\La$.  The inclusion $M^\top\hookrightarrow\mathscr D_\tJ(M)$ induces, by the adjunction for Demazure induction over $\mathfrak l_\tJ$, an $\mathfrak l_\tJ$-module map
\[
\mathscr D_\tJ(M^\top)\longrightarrow\mathscr D_\tJ(M).
\]
Its restriction to $M^\top$ is the given inclusion, so it is nonzero. Since its source is irreducible, it is injective and its image is $U(\mathfrak l_\tJ)v$. Thus $\mathscr D_\tJ(M^\top)\cong U(\mathfrak l_\tJ)v\subset\mathscr D_\tJ(M)$. Applying $\tp_{\tJ}$ to the $\mathfrak l_{\tJ}$-lowest weight vector of $\mathscr D_{\tJ}(M^\top)$, we find that $F_i$ ($i \in \tJ$) acts by zero on this vector. Hence the action of every $\wgt$-eigenvector in $\tp_{\tJ}\setminus\tb$ on this vector is zero. Therefore, by the PBW theorem, the $\tp_{\tJ}$-submodule of $\mathscr D_{\tJ}(M)$ generated by $\mathscr D_{\tJ}(M^\top)$ is already generated by $\mathscr D_{\tJ}(M^\top)$ as a $\tb$-module. It follows that $\mathscr D_{\tJ}(M^\top)$ generates $\mathscr D_{\tJ}(M)$ as a $\tb$-module.
\end{proof}

For a root ideal $\Psi$ and each $1 \le i \le n$, we set
$$\hs_i ( \Psi ) := \max \left( \{l \mid E_{li} \in \gn( \Psi )\} \cup \{ 0 \} \right),$$
which is a natural extension of the quantity with the same name in~\cite[\S1.2]{Kat23a}.
We have $0 \le \hs_i (\Psi) < i$ for each $1 \le i \le n$. We define the right exact functor
\begin{equation}
\mathscr C_j ^{\Psi} ( m; \bullet) :=\mathscr D_{[j,\hs_j(\Psi))} (\C_{m \La_j} \otimes \bullet) \qquad (1 \le j\le n, m \in \Z)\label{eqn:Cjpdef}
\end{equation}
on the category of finite-dimensional $\wgt$-semisimple $\tb$-modules.

\begin{lem}\label{lem:Dstab}
For a root ideal $\Psi$ and a sequence $m_j,m_{j+1},\ldots,m_n \in \Z$, the $\tb$-module
$$\left( \mathscr C_{j}^\Psi ( m_j; \bullet ) \circ \mathscr C_{j+1}^\Psi ( m_{j+1}; \bullet ) \circ \cdots \circ \mathscr C_n^\Psi ( m_n; \bullet ) \right) ( \C )$$
extends to a $\tp_{[j,\hs_j(\Psi))}$-module.
\end{lem}

\begin{proof}
We first observe that
\[
(n,\hs_n(\Psi)) \subset [1,n), \qquad
(i,\hs_i(\Psi)) \subset [i+1,\hs_{i+1}(\Psi)) \quad (1 \le i < n).
\]
Since tensoring with $\C_{m_i\La_i}$ commutes with $\mathscr D_j$ for $j \neq i$, the claim follows by repeated application of Lemma~\ref{lem:Levi}.
\end{proof}

The following identification is essentially a reformulation of~\cite[Theorem~4.1]{Kat23a} (which is stated for $\la \in \Par(n)$ and its generalization to $\la \in \sP$ follows by an interpretation of the construction of $\cO_{\sX_\Psi} ( \la )$ in~\cite[Lemma~3.6 and Proposition~3.7]{Kat23a}) and will be the starting point of the socle analysis.

\begin{thm}\label{thm:catalan}
For a root ideal $\Psi$, $r \in \Z$ and $\la =\sum_{i=1}^n m_i \varpi_i \in \sP$, we have
$$H^{r} ( \sX_\Psi, \cO_{\sX_\Psi} ( \la ))^* \cong \mathbb L^{- r} \left( \mathscr C_1^\Psi ( m_1 ; \bullet ) \circ \mathscr C_2^\Psi (m_2; \bullet) \circ \cdots \circ \mathscr C_n^\Psi (m_n; \bullet) \right) ( \C ).$$
\end{thm}

\begin{proof}
For each $1 \le d < n$, we have
$$w_0 w_0^{[d,n)} = (s_{n-1} \cdots s_2 s_1 )(s_{n-1} \cdots s_3 s_2) \cdots (s_{n-1}\cdots s_d s_{d-1}).$$
Thus, we have
$$\mathscr D_{w_0 w_0^{[d,n)}} = \mathscr D_{[1,n)} \circ \mathscr D_{[2,n)} \circ \cdots \circ \mathscr D_{[d-1, n)}.$$
Comparing this with~\cite[Lemma~2.12]{Kat23a}, we deduce
$$N^{\Psi}_{w_0} ( \la ) \cong \left( \mathscr C_1^\Psi ( m_1 ; \bullet ) \circ \mathscr C_2^\Psi (m_2; \bullet) \circ \cdots \circ \mathscr C_n^\Psi (m_n; \bullet) \right) ( \C ),$$
where the left-hand side is defined in~\cite[(2.2)]{Kat23a}. Thus, the assertion follows from~\cite[Theorem~4.1]{Kat23a}.
\end{proof}

\begin{lem}\label{lem:Dincl}
Let $i_1,i_2,\ldots,i_{l} \in \tI_\af$, and let $m_1,\ldots,m_l \in \Z_{\ge 0}$. We define a graded $\tb$-module $M_r$ for $1 \le r \le l$ by setting $M_{l+1} := \C$ and
$$M_r := \mathscr D_{i_r} ( \C _{m_r \La_{i_r}} \otimes M_{r+1} ) \qquad (1 \le r \le l)$$
by induction. Then we have
$$\C_{m_r \La_{i_r}} \otimes M_{r+1} \subset M_r \qquad (1 \le r \le l).$$
\end{lem}

\begin{proof}
The module $\C = \C_0$ is an affine Demazure module corresponding to the zero-weight in the sense of~\cite[\S1.4]{Kat23a}. In view of this, each step follows from~\cite[Corollary~1.18]{Kat23a}.
\end{proof}

\begin{cor}\label{cor:Cincl}
For a root ideal $\Psi$ and $\la = \sum_{j=1}^n m_j \varpi_j \in \Par(n)$, we have
\begin{align*}
\C_{m_j \La_j} \otimes & \, \left( \mathscr C_{j+1}^\Psi (m_{j+1}; \bullet ) \circ \mathscr C_{j+2}^\Psi (m_{j+2}; \bullet ) \circ \cdots \circ \mathscr C_n^\Psi (m_{n}; \bullet ) \right) ( \C )\\
&  \subset \left( \mathscr C_{j}^\Psi (m_{j}; \bullet ) \circ \mathscr C_{j+1}^\Psi (m_{j+1}; \bullet ) \circ \cdots \circ \mathscr C_n^\Psi (m_{n}; \bullet ) \right) ( \C )
\end{align*}
for each $1 \le j \le n$.
\end{cor}

\begin{proof}
This follows from a successive application of Lemma~\ref{lem:Dincl} to the definition.
\end{proof}

\begin{defn}[Shallow ideal]
A root ideal $\Psi$ is said to be shallow if and only if $E_{ij} \in \mathfrak{n} (\Psi)$ and $E_{i,j-1} \notin \gn(\Psi)$ implies $E_{i+1,j}\notin \gn (\Psi)$ for each $1 \le i < j \le n$.
\end{defn}

\begin{thm}\label{thm:HLss}
Let $\Psi$ be a shallow root ideal. For each $\la \in \Par(n)$, the module $\HL^\Psi_\la \ (\cong H^0 ( \sX_\Psi, \cO_{\sX_\Psi} ( \la )) )$ has a simple socle as a $\g[z]$-module.
\end{thm}

\begin{ex}[$n=6,k=7$]\label{ex:7schur}
We consider the case $\la=(6,5^2,3,1^2) \in \Par_{21}(6)$ (i.e., $\la = \varpi_1+2\varpi_3+2\varpi_4+\varpi_6$) and
$$\Psi = \{\al_{1,3},\al_{1,4},\al_{1,5},\al_{1,6},\al_{2,5},\al_{2,6},\al_{3,6}\},$$
which is shallow and is depicted as follows:
\begin{center}
\begin{tikzpicture}[scale=0.5, thick, >=Stealth]
  \def\s{1.0}

  \draw[thick] (0,0) rectangle (6*\s,6*\s);

  \foreach \i in {0,...,5} {
    \pgfmathsetmacro{\x}{6 - \i - 1}
    \pgfmathsetmacro{\y}{\i}
    \draw (\x*\s,\y*\s) rectangle ++(\s,\s);
  }

  \foreach \x/\y in {
    2/5,
    4/4, 5/4,
    3/5, 4/5, 5/5,
    5/3} {
    \fill[red!55] (\x*\s,\y*\s) rectangle ++(\s,\s);
    \draw[thick, dash pattern=on 2pt off 5pt] (\x*\s,\y*\s) rectangle ++(\s,\s);
  }

  \node at (3.5*\s,5.6*\s) {$\Psi$};

  \draw[<->, shorten >=5pt, shorten <=5pt]
    (5.5*\s,3*\s) -- (5.5*\s,6*\s)
    node[midway, fill=white, inner sep=1pt] {$\mathtt{h}_6$};

  \draw[<->, shorten >=5pt, shorten <=5pt]
    (4.5*\s,4*\s) -- (4.5*\s,6*\s)
    node[midway, fill=white, inner sep=1pt] {$\mathtt{h}_5$};

\end{tikzpicture}
\end{center}

This realizes a module whose graded character is $\ws^{(7)}_\la$ in~\S\ref{sec:ksch} through Theorem~\ref{thm:catalan}. We have
$$\hs_{1} (\Psi) = 0, \hs_{2} (\Psi) = 0, \hs_{3} (\Psi) = 1, \hs_{4} (\Psi) = 1, \hs_{5} (\Psi) = 2, \hs_{6} (\Psi) =3.$$
The $\wgt$-weights of the unique simple quotients as $\tb$-modules of
$$\left( \mathscr C_{j}^\Psi (m_{j}; \bullet ) \circ \cdots \circ \mathscr C_n^\Psi (m_{n}; \bullet ) \right) ( \C ) \qquad (1 \le j \le 6)$$
presented in the proof of Theorem~\ref{thm:HLss} are, modulo $\C\delta$,
\begin{align*}
&\La^{(6)}= (1^2,2,1^2,0) + \wp, \, \La^{(5)} = (1,2^2,1,0^2) + \wp, \, \La^{(4)} = (6,4^2,0^3) + 3 \wp,\\
&\La^{(3)}= (11,6,0^3,3) + 5 \wp, \, \La^{(2)} = (11,0^3,3,6) + 5 \wp, \, \La^{(1)} = (0^3,3,6,12) + 6 \wp.
\end{align*}
The $\Sym_6$-orbit of $\mathsf{pr}(\La^{(1)})$ in $\Comp(6)$ contains the partition $(12,6,3)$, whose conjugate is
\[
(12,6,3)' = (3^3,2^3,1^6) = \la^{\omega_7}
\]
by Example~\ref{ex:n7conj}.
\end{ex}

\begin{proof}[Proof of Theorem~\ref{thm:HLss}]
We prove that the dual module has a simple head as a $\tb$-module. We adopt the notation from Theorem~\ref{thm:catalan}, omitting $\Psi$ from the notation for brevity. The proof proceeds by a downward induction on $i$, and the main task is to construct, at each stage, a cyclic vector whose Levi span remains irreducible and nonzero.

Since $\la = \sum_j m_j \varpi_j$, we have $m_j \in \Z_{\ge 0}$ by $\la \in \Par(n)$. Set
\[
m^{(j)} := \sum_{i=j}^n m_i \quad (1 \le j \le n) \quad \text{and} \quad m^{(n+1)} := 0.
\]

Let $N_{n+1} = \C$ as a $\tb$-module, $\La^{(n+1)} = 0 \in \sP_\af$ be its $\wgt$-weight and $v_{n+1} = 1 \in N_{n+1}$. We define $\wgt$-weights
\[
\La^{(i)} = \la^{(i)} + d_i \delta + m^{(i)} \wp, \quad \text{with } \la^{(i)} \in \sP,
\]
inductively as
$$\La^{(i)} := ( s_{\hs_i -1} \cdots s_1 s_0 s_{n-1} s_{n-2} \cdots s_i ) ( \La^{(i+1)} + m_i \La_i ) \qquad (1 \le i \le n)$$
from the initial term $\la^{(n+1)} = (0^n)$. (If $\hs_i = 0$, then the sequence starts by $s_{n-1} s_{n-2}\cdots$.)

The goal is to construct $\tb$-cyclic vectors $v_i$ in $\tb$-modules $N_i$ of the form
\[
N_i := \bigl( \mathscr C_{i} (m_i; \bullet ) \circ \cdots \circ \mathscr C_n (m_n; \bullet ) \bigr) (\C),
\]
such that the vector $v_i$, together with its $\wgt$-weight $\La^{(i)}$, satisfies the hypotheses of Proposition~\ref{prop:gen}, assuming that there exists a $\tb$-cyclic vector $v_{i+1} \in N_{i+1}$ with $\wgt$-weight $\La^{(i+1)}$. The case $i+1=n+1$ is immediate from the definition.

\begin{claim}\label{claim:weightrec}
The sequence $\la^{(i)}$ satisfies the following recurrence depending on $\hs_i$:
\begin{itemize}
\item If $\hs_i > 0$, then $\la^{(i)}$ is given by
\begin{equation}
\la^{(i)}_j = \begin{cases}
\la^{(i+1)}_{j} + m_i & (\hs_i < j < i)\\
\la^{(i+1)}_{j+1} + m_i & (1 \le j < \hs_i)\\
\la^{(i+1)}_i + m_i+m^{(i)} & (j = \hs_i)\\
\la^{(i+1)}_{j+1} & (i \le j < n)\\
\la^{(i+1)}_1 - m^{(i+1)} & (j = n)
\end{cases}\label{eqn:recur>0}.
\end{equation}
  \item If $\hs_i = 0$, then $\la^{(i)}$ is given by
\begin{equation}
\la^{(i)}_j = \begin{cases}
\la^{(i+1)}_{j} + m_i & (1 \le j < i)\\
\la^{(i+1)}_{j+1} & (i \le j < n)\\
\la^{(i+1)}_i + m_i & (j = n)\\
\end{cases}.\label{eqn:recur=0}
\end{equation}
\end{itemize}
\end{claim}

\begin{proof}
This follows by a direct computation from the definition of $\La^{(i)}$.
\end{proof}

\begin{claim}\label{claim:mono}
Assume $\Psi$ is shallow, so that
\begin{equation}
\hs_{i} \le \hs_{i+1} \le \hs_i+1 \qquad (1 \le i < n).\label{eqn:monidx}
\end{equation}
If $\hs_i > 0$ for some $1 \le i \le n$, then $\la^{(i)}$ satisfies
\begin{equation}
\la^{(i)}_1 \le \cdots \le \la^{(i)}_{\hs_i}
\ge
\la^{(i)}_{\hs_i+1} \ge \cdots \ge \la^{(i)}_{i-1}
\quad\text{and}\quad
\la^{(i)}_{\hs_i} \le \la^{(i)}_{i-1} + m^{(i-1)}.
\label{eqn:recur3}
\end{equation}
\end{claim}

\begin{proof}
We argue by downward induction on $i$, starting from $i=n$.

\emph{Base case $i=n$.}
From~\eqref{eqn:recur>0} with $i=n$ and $\la^{(n+1)}\equiv 0$, $m^{(n+1)}=0$, we obtain
\[
\la^{(n)}_j=
\begin{cases}
m_n & (1\le j<\hs_n)\\
2m_n & (j=\hs_n)\\
m_n & (\hs_n<j<n)\\
0 & (j=n)
\end{cases}.
\]
Hence $\la^{(n)}_1\le\cdots\le\la^{(n)}_{\hs_n}\ge\cdots\ge\la^{(n)}_{n-1}\ge\la^{(n)}_n$, and
$\la^{(n)}_{\hs_n}=2m_n\le m_n+(m_{n-1}+m_n)\le\la^{(n)}_{n-1}+m^{(n-1)}$, proving~\eqref{eqn:recur3} for $i=n$.

\emph{Induction step.}
Assume~\eqref{eqn:recur3} holds at stage $i+1$. By~\eqref{eqn:monidx}, $\hs_i\in\{\hs_{i+1},\,\hs_{i+1}-1\}$.
Using~\eqref{eqn:recur>0}, for $1\le j < i$ we have
\[
\la^{(i)}_j=
\begin{cases}
\la^{(i+1)}_{j+1}+m_i & (1\le j<\hs_i)\\
\la^{(i+1)}_{i}+m_i+m^{(i)} & (j=\hs_i)\\
\la^{(i+1)}_{j}+m_i & (\hs_i<j<i)
\end{cases},
\qquad 
\la^{(i)}_{n}=\la^{(i+1)}_{1}-m^{(i+1)}.
\]
Since the blocks $[1,\hs_i-1]$ and $[\hs_i+1,i-1]$ are obtained from monotone blocks of $\la^{(i+1)}$ by adding the constant $m_i$, their monotonicity is preserved. For the peak, by the induction hypothesis at stage $i+1$,
\[
\la^{(i)}_{\hs_i-1}=\la^{(i+1)}_{\hs_i}+m_i
\ \le\
\la^{(i+1)}_{\hs_{i+1}}+m_i
\ \le\
\la^{(i+1)}_{i}+m^{(i)}+m_i
=\la^{(i)}_{\hs_i},
\]
when $\hs_i > 1$, and similarly
\(
\la^{(i)}_{\hs_i+1}=\la^{(i+1)}_{\hs_i+1}+m_i
\le \la^{(i+1)}_{\hs_{i+1}}+m_i
\le \la^{(i)}_{\hs_i}
\)
when $\hs_i + 1 < i$.

Finally,
\[
\la^{(i)}_{\hs_i}-m^{(i)}
= \la^{(i+1)}_{i}+m_i
\ \le\ \la^{(i+1)}_{i-1}+m_i
= \la^{(i)}_{i-1}
\ \le\ \la^{(i)}_{i-1}+m_{i-1}
\]
when $\hs_i<i-1$, and $\la_{\hs_i}^{(i)}=\la_{i-1}^{(i)}$ when $\hs_i=i-1$, which yields the rightmost bound in~\eqref{eqn:recur3}.
\end{proof}

\begin{claim}\label{claim:monotone}
Keep the assumptions of Claim~\ref{claim:mono}, and let $2\le i\le n$. Then:
\begin{itemize}
\item If $\hs_i>0$, we have
\begin{equation}
\la_{i}^{(i)} \le \cdots \le \la_{n}^{(i)} \le \bigl(\la_{1}^{(i)}-m^{(i)}\bigr)
\le \cdots \le \bigl(\la_{\hs_{i-1}}^{(i)}-m^{(i)}\bigr) \le \la_{i-1}^{(i)}+m_{i-1}.
\label{eqn:recur4}
\end{equation}
When $\hs_{i-1}=0$, we understand the middle segment
$\bigl(\la^{(i)}_{1}-m^{(i)}\bigr) \le \cdots \le \bigl(\la^{(i)}_{\hs_{i-1}}-m^{(i)}\bigr)$
to be empty, so that~\eqref{eqn:recur4} reads
$\la^{(i)}_{i} \le \cdots \le \la^{(i)}_{n} \le \la^{(i)}_{i-1}+m_{i-1}$.
\item If $\hs_i=0$, we have
\begin{equation}
\la_i^{(i)} \le \cdots \le \la_{n}^{(i)} \le \la_{i-1}^{(i)} \le \la_{i-2}^{(i)} \le \cdots \le \la_{1}^{(i)}.
\label{eqn:extr}
\end{equation}
\end{itemize}
\end{claim}

\begin{proof}
We proceed by downward induction on $i$ from the case $i=n$. The case $i=n$ is immediate by Claim~\ref{claim:mono} and its proof.

\emph{Case $\hs_i>0$.}
For $i\le j<n$, $\la_{j}^{(i)}=\la_{j+1}^{(i+1)}$ and $\la_n^{(i)}=\la_1^{(i+1)}-m^{(i+1)}$, so
$\la_{i}^{(i)}\le\cdots\le \la_{n}^{(i)}$ follows from the induction hypothesis at stage $i+1$.
To compare $\la_n^{(i)}$ with $\la_1^{(i)}-m^{(i)}$, distinguish two subcases by~\eqref{eqn:recur>0}:
\[
\la_1^{(i)}-m^{(i)}=
\begin{cases}
\la_2^{(i+1)}-m^{(i+1)} & (\hs_i\ge 2)\\
\la_i^{(i+1)}+m_i & (\hs_i=1)
\end{cases}.
\]
In the first subcase, $\la_1^{(i+1)}-m^{(i+1)}\le \la_2^{(i+1)}-m^{(i+1)}$ since the initial segment of $\la^{(i+1)}$ is nondecreasing. In the second subcase,
\[
\la_1^{(i+1)}-m^{(i+1)}\ \le\ \la_{\hs_{i+1}}^{(i+1)}-m^{(i+1)}
\ \le\ \la_i^{(i+1)}+m_i
=\la_1^{(i)}-m^{(i)},
\]
where we used Claim~\ref{claim:mono} at stage $i+1$ for the last inequality.
Thus $\la_n^{(i)}\le \la_1^{(i)}-m^{(i)}$.
Finally,~\eqref{eqn:recur3} gives
\(
\la_{1}^{(i)}-m^{(i)} \le \cdots \le \la_{\hs_{i-1}}^{(i)}-m^{(i)} \le \la_{i-1}^{(i)}+m_{i-1},
\)
which completes~\eqref{eqn:recur4}.

\emph{Case $\hs_i=0$.}
From~\eqref{eqn:recur=0} we have $\la_{j}^{(i)}=\la_{j+1}^{(i+1)}$ for $i\le j<n$ and
$\la_n^{(i)}=\la_i^{(i+1)}+m_i$. By~\eqref{eqn:monidx} we have $\hs_{i+1}\le \hs_i+1=1$, so the
induction hypothesis at stage $i+1$ reads
\[
\la_{i+1}^{(i+1)} \le \cdots \le \la_{n}^{(i+1)} \le \la_{i}^{(i+1)}+m_i
\quad\text{and}\quad
\la_{i}^{(i+1)} \le \cdots \le \la_1^{(i+1)} :
\]
the first chain is (possibly a part of)~\eqref{eqn:recur4} if $\hs_{i+1}=1$, and~\eqref{eqn:extr}
if $\hs_{i+1}=0$, while the second is~\eqref{eqn:recur3} if $\hs_{i+1}=1$ and again~\eqref{eqn:extr}
if $\hs_{i+1}=0$. Adding $m_i$ to the second chain, we obtain
\[
\la_{i+1}^{(i+1)} \le \cdots \le \la_{n}^{(i+1)} \le \la_{i}^{(i+1)}+m_i \le \cdots \le \la_1^{(i+1)}+m_i,
\]
which is~\eqref{eqn:extr} at stage $i$.
\end{proof}

The next claim isolates the representation-theoretic input needed at each induction step.

\begin{claim}\label{claim:demazure-irr}
Let $N_{i+1}^\top$ denote the $\mathfrak l_{(i,\hs_{i})}$-span of $v_{i+1} \in N_{i+1}$. 
Suppose that $N_{i+1}^\top$ is irreducible as an $\mathfrak l_{(i,\hs_{i})}$-module. 
Assume furthermore that the tuple
\[
( \la_{i}^{(i+1)} + m_{i}, \la_{\hs_{i}}^{(i+1)} - m^{(i+1)}, \ldots, \la_1^{(i+1)} - m^{(i+1)}, \la_{n}^{(i+1)}, \ldots, \la_{i+1}^{(i+1)} )
\]
is weakly decreasing if $\hs_i > 0$, or that
\[
( \la^{(i+1)}_{i}+m_{i}, \la^{(i+1)}_{n}, \la^{(i+1)}_{n-1}, \ldots, \la^{(i+1)}_{i+1} )
\]
is weakly decreasing if $\hs_i = 0$. 
Then $\mathscr D_{[i,\hs_i)}$ applied to $\C_{m_i \La_i} \otimes N_{i+1}^{\top}$ is an irreducible and nonzero $\mathfrak l_{[i,\hs_i)}$-module. 
In this case, we have
\[
\C_{m_i \La_i} \otimes N_{i+1}^{\top}\subset
\mathscr D_{[i,\hs_i)} \bigl( \C_{m_i \La_i} \otimes N_{i+1}^{\top} \bigr),
\]
and the $\mathfrak l_{[i,\hs_i)}$-lowest weight vector of 
$\mathscr D_{[i,\hs_i)} \bigl( \C_{m_i \La_i} \otimes N_{i+1}^{\top} \bigr)$ 
has $\wgt$-weight $\La^{(i)}$.
\end{claim}

\begin{proof}
We identify $\mathfrak l_{(i,\hs_{i})} \subset \mathfrak l_{[i,\hs_{i})}$ with $\mathfrak{sl} (n-i+\hs_{i},\C) \subset \mathfrak{sl} (n-i+1+\hs_{i},\C)$, respectively. By Claim~\ref{claim:monotone}, the highest weight (as $\mathfrak{sl} (n-i+1+\hs_{i},\C)$-weight) of the $\mathfrak l_{(i,\hs_{i})}$-span of $\C_{m_i \La_i} \otimes v_{i+1}$ is represented by the $(n-i+1+\hs_{i})$-tuple
\begin{equation}
( \la_{i}^{(i+1)} + m_{i},\la_{\hs_{i}}^{(i+1)}-m^{(i+1)}, \ldots,\la_1^{(i+1)} -m^{(i+1)}, \la_{n}^{(i+1)},\ldots,\la_{i+1}^{(i+1)} )\label{eqn:Lwt}
\end{equation}
if $\hs_{i} > 0$, and $(\la_{i}^{(i+1)}+m_{i},\la_n^{(i+1)},\ldots,\la_{i+1}^{(i+1)})$ if $\hs_{i} = 0$. Hence the above weight is dominant as a weight of $\mathfrak{sl} (n-i+1+\hs_{i},\C)$ if and only if the corresponding sequence is weakly decreasing. Let $v_{i}^+$ denote the corresponding highest weight vector. Then we have
$$\mathscr D_{[i,\hs_i)} ( \C _{m_i \La_i} \otimes N_{i+1}^{\top} ) = \mathscr D_{w_0^{[i,\hs_i)}} ( \C v_i^+ ),$$
which is $\mathfrak{l}_{[i,\hs_i)}$-irreducible if and only if the weakly decreasing condition is satisfied and is zero otherwise. This implies the first assertion.

The second assertion follows by examining the extremal weight vector provided by this irreducible $\mathfrak{l}_{[i,\hs_i)}$-module.
\end{proof}

We now complete the proof of Theorem~\ref{thm:HLss}. Assuming the induction hypothesis at stage $i+1$, Corollary~\ref{cor:Cincl} and Claim~\ref{claim:demazure-irr} show that $N_i$ admits a $\tb$-cyclic vector of weight $\La^{(i)}$ satisfying the hypotheses of Proposition~\ref{prop:gen}. Since $\Psi$ is shallow, Claims~\ref{claim:mono} and~\ref{claim:monotone} ensure that the required dominance conditions remain valid at each step, so the resulting Levi-module remains nonzero and irreducible throughout the induction.

Thus the downward induction closes. Therefore the $\tb$-module
\[
N_{w_0}^\Psi(\la) \cong H^0(\sX_\Psi,\cO_{\sX_\Psi}(\la))^*
\]
has a simple head. Equivalently, $\HL_\la^\Psi$ has a simple socle. This proves the theorem.
\end{proof}

\section{$k$-Schur modules and Chen--Haiman modules}\label{sec:ksch}
We retain the setting of the previous section. Fix $k \in \Z_{> 0}$. For $\la\in\sP$ set
\begin{equation}
\mathsf o^k_\la(i):=\max\{i, i+k-\la_i\}\qquad(1\le i\le n).\label{eqn:okdef}
\end{equation}
When it exists, define the root ideal $\Psi[\la,k]$ by
\[
\gn(\Psi[\la,k])=\bigoplus_{1 \le i \le \mathsf o^k_\la(i) < j \le n} \C E_{ij} \subset \gn.
\]

\begin{lem}\label{lem:Psi}
For $k \in \Z_{\ge 1}$ and $\la \in \sP$, the root ideal $\Psi[\la,k]$ exists whenever $\la_i +1 \ge \la_{i+1}$ for all $1 \le i < n$. In particular, this condition is satisfied for every $\la \in \sP^+$.
\end{lem}

\begin{proof}
The verification is straightforward.
\end{proof}

For $\la \in \Par_m (n)$, define
\[
\sk_\la := \HL_\la^{\Psi[\la,k]} \in \g[z]\gmod_m,
\]
and refer to $\sk_\la$ as the $k$-Schur module. For each $\la \in \Par_m^{(k)}$, we view $\la$ in $\Par_m^{(k)}(m+1)$ and set
\[
\ssc^{(k)}_{\la} := \SW \bigl(\sk_\la\bigr) \in A_m\gmod,
\]
where $\SW$ denotes the Schur--Weyl duality functor for $\g=\mathfrak{gl}(m+1)$. We call $\ssc^{(k)}_{\la}$ the Chen--Haiman module. As we shall see in Corollary~\ref{cor:CHindep}, alternative definitions of Chen--Haiman modules produce isomorphic objects.

\begin{rem}
As shown in Theorem~\ref{thm:ksch}, the graded $\g[z]$-module $\sk_\la$
realizes the $k$-Schur function (in the sense of~\cite{BMPS}) only when
$\la\in\Par_m^{(k)}(n)$.
On the other hand, we will need $\sk_\la$ for general
$\la\in\Par_m(n)$ in order to state Theorem~\ref{thm:skfilt}.
\end{rem}

Recall that the basic terminology concerning filtrations was fixed in \S\ref{subsec:conv}, where the notion of ``filtration by $\{Z_\la\}_\la$" was introduced.

\begin{defn}
A module $M \in \g[z]\gmod_m$ admits an $\sk$-filtration (resp.\ an $\HL$-filtration) if there exists a filtration of $M$ by $\{\sk_\la\}_{\la \in \Par_m^{(k)}(n)}$ (resp.\ by $\{\HL_\la\}_{\la \in \Par_m(n)}$).

Similarly, $M \in A_m\gmod$ admits an $\ssc^{(k)}$-filtration (resp.\ a $K$-filtration, $\tK$-filtration) if $M$ admits a filtration by $\{\ssc^{(k)}_\la\}_{\la \in \Par^{(k)}_m}$ (resp.\ by $\{K_\la\}_{\la \in \Par_m}$,  $\{\tK_\la\}_{\la \in \Par_m}$, respectively).

When $\dim M<\infty$, these are equivalent to the existence of an increasing exhaustive filtration (or a finite filtration) with the same associated graded type.
\end{defn}

\begin{thm}\label{thm:ksch}
For every $\la \in \Par_m^{(k)} (n)$,
\[
\ws_\la^{(k)} = \sum_{\mu\in\Par_m(n)} [\sk_\la : V_{\mu}^*]_q \,\ws_{\mu},
\]
where $\ws_\la^{(k)}$ denotes the $k$-Schur function in the $q$-parameter convention of~\cite[(2.4)]{BMPS}. Moreover,
\begin{equation}
[\sk_\la : V_{\mu}^*] =
\begin{cases}
1 & \text{if }\mu=\la\\
1 & \text{if }\mu=(\la^{\omega_k})'\\
0 & \text{if }\mu \not\ge \la \ \text{or}\ \mu \not\le (\la^{\omega_k})'
\end{cases}\label{eqn:ksch-tri}.
\end{equation}
\end{thm}

\begin{proof}
By~\cite[Theorem~2.4]{BMPS}, the stated expansion is a special case of Theorem~\ref{thm:Cat}\textup{(4)}. By~\cite[\S6(27)]{LM07}, we have
\[
[\sk_\la : V_{\mu}^*] =
\begin{cases}
1 & \text{if }\mu=\la\\
0 & \text{if }\mu \not\ge \la
\end{cases}.
\]
In addition, the $k$-conjugation symmetry for $k$-Schur functions at $q=1$~\cite[(8)]{LM07} implies
\[
[\sk_\la : V_{\mu}^*] = [\sk_{\la^{\omega_k}} : V_{\mu'}^*]
\qquad (\la \in \Par_m^{(k)},\ \mu \in \Par_m).
\]
Combining the usual triangularity with this $k$-conjugation symmetry yields~\eqref{eqn:ksch-tri}.
\end{proof}

\begin{rem}
In~\cite[\S10]{BMPS2}, the authors suggest different root ideals that give rise to the same $\sk_\la$. In the case of Example~\ref{ex:7schur}, this implies that we can also employ
$$\Psi' := \Psi \cup \{\al_{3,5}\} = \Psi [\la,7] \cup \{\al_{3,5}\}$$
instead of $\Psi$. In this case, the natural restriction map
$$\HL_\la^{\Psi'} \longrightarrow \HL_\la^{\Psi}$$
is an isomorphism. The same holds in the setting of~\cite[Example~10.5]{BMPS}.
\end{rem}

\begin{lem}\label{lem:gchcomp}
For each $\la \in \Par_m^{(k)} (n)$ and $\mu \in \Par_m$,
\[
[\sk_\la : V_{\mu}^*]_q = [\ssc^{(k)}_\la : L_{\mu}]_q.
\]
In particular, when $|\la| \le k$, we have
\[
\sk_\la \cong V_\la^*, \quad \text{and} \quad \ssc^{(k)}_\la \cong L_\la.
\]
\end{lem}

\begin{proof}
Since $\ws_\la^{(k)}$ does not depend on $n$ as long as $\la_{n+1} =0$, we can assume $n=m+1$ by Theorem~\ref{thm:ksch}. Then we apply Theorem~\ref{thm:FKM} to conclude the first assertion. The second assertion follows from Theorem~\ref{thm:ksch}, the first assertion, and $\ws^{(k)}_\la =\ws_\la$ when $|\la| \le k$ (\cite[\S4.1]{LLMSSZ}).
\end{proof}

\begin{cor}\label{cor:unitri}
For each $\la \in \Par_m^{(k)}$ and $\mu \in \Par_m$, we have
\[
[\ssc^{(k)}_\la : L_{\mu}]_q = \begin{cases} 1 & (\la = \mu)\\ 0 & (\la \not\le \mu) \end{cases}.
\]
\end{cor}

\begin{proof}
Combine Theorem~\ref{thm:ksch} and Lemma~\ref{lem:gchcomp}. (Alternatively, compare~\eqref{eqn:Ktri} with the surjection $K_\la \twoheadrightarrow \ssc^{(k)}_\la$ afforded by Corollary~\ref{cor:Cat-sg} through $\SW$.)
\end{proof}

\begin{thm}\label{thm:tri}
Let $\Psi$ be a root ideal and $\la \in \sP^+$. For each $\mu \in \sP^+$,
\[
\bigl[ H^0\bigl(T^*_{\Psi} X, \cO_{T^*_{\Psi} X}(\la)\bigr) : V_\mu^* \bigr]\neq 0
\]
implies $\mu \ge \la$. Consequently, if $\la \notin \Par(n)$, the module $H^0\bigl(T^*_{\Psi} X, \cO_{T^*_{\Psi} X}(\la)\bigr)$ contains no nonzero irreducible polynomial $G$-subrepresentation.
\end{thm}

\begin{proof}
By Corollary~\ref{cor:Cat-sg}, we have
\[
\bigl[ H^0\bigl(T^*_{\Psi} X, \cO_{T^*_{\Psi} X}(\la)\bigr) : V_\mu^* \bigr] \le \bigl[ H^0\bigl(T^* X, \cO_{T^* X}(\la)\bigr) : V_\mu^* \bigr].
\]
By~\cite[Theorem~5.1]{Kat23a} and~\eqref{eqn:stdmult}, we have
\begin{align*}
\bigl[ H^0\bigl(T^* X, \cO_{T^* X}(\la)\bigr) : V_\mu^* \bigr] & = \lim_{s \to \infty} \dim \Hom_G ( V_{\mu}^*, H^0(\sX_{\Delta^+}, \cO_{\sX_{\Delta^+}}(\la+s\varpi_n))\otimes \C_{s\varpi_n} )\\
& = \lim_{s \to \infty} \dim \Hom_G ( V_{\mu+s\varpi_n}^*, H^0(\sX_{\Delta^+}, \cO_{\sX_{\Delta^+}}(\la+s\varpi_n)) )\\
& = \lim_{s \to \infty} K_{\mu+s\varpi_n, \la+s\varpi_n} ( 1 ),
\end{align*}
and the last term is nonzero only if $\mu \ge \la$.

If $\la \notin \Par(n)$, then $\mu \ge \la$ forces $\mu \notin \Par(n)$ by Lemma~\ref{lem:interval}, so no dominant polynomial weight can occur. This completes the proof.
\end{proof}

\begin{lem}\label{lem:ksch}
For every $\la \in \Par_m (n)$ such that $\la_i \notin \{1,2,\ldots,k-1\}$ for every $1 \le i \le n$, we have
\[
\sk_\la = \HL_\la.
\]
In particular, we have $\mathtt{s}^{(1)}_\la = \HL_\la$.
\end{lem}

\begin{proof}
By assumption, we have $\la_i \ge k$ or $\la_i =0$. This implies $\Psi [\la,k] = \Delta^+$ up to $\ell ( \la )$-th row. Recall that we have
\[
\HL_0^\Psi = V_0^*
\]
as $\mathfrak{gl} (n-\ell(\la))$-modules regardless of the root ideal $\Psi$ of size $n-\ell ( \la )$ by~\eqref{eqn:stdmult} and Corollary~\ref{cor:Cat-sg}. Furthermore $\gn ( \Psi [\la,k])$ contains the upper-right block of size $\ell(\la) \times (n-\ell(\la))$ in $\gn$. It follows that 
\[
\gch \HL_\la^{\Psi [\la,k]} \equiv \gch \HL_\la^{\Delta^+} \mod (q - 1) 
\]
by~\cite[Corollary~3.4]{BMPS2}. The natural surjective map $\HL_\la^{\Delta^+} \twoheadrightarrow \HL_\la^{\Psi [\la,k]}$ obtained from Corollary~\ref{cor:Cat-sg} must be an isomorphism by the above comparison of characters. Consequently, we have $\sk_\la = \HL_\la$.
\end{proof}

\begin{thm}\label{thm:dla}
Let $k \in \Z_{>0}$. For each $\la \in \Par^{(k)}_m (n)$, let $d \in \Z_{\ge 0}$ be the maximal degree such that $ ( \sk_\la )_d \neq 0$. Then
\begin{equation}\label{eqn:dlaest}
d \le d_k(\la).
\end{equation}
\end{thm}

\begin{ex}
In the setting of Example~\ref{ex:7schur}, we have $k=7$ and
\[
\la=(6,5^2,3,1^2), \qquad (\la^{\omega_k})'=(12,6,3),
\]
and in this case $d_7(\la)=8$.
\end{ex}

\begin{proof}[Proof of Theorem~\ref{thm:dla}]
By Theorem~\ref{thm:kcore}, there is a $k+1$-core $\kappa$ corresponding to $\la$, denoted $\mathfrak c(\la)$.

The $k$-Pieri rule~\cite[Theorem~2.3]{BMPS} states that
\begin{equation}\label{eqn:t1e}
\we_{1}^{\perp}\, \ws_\la^{(k)} = \sum_{(\mu,r)} q^{\mathrm{spin}(\mu,r)} \, \ws_\mu^{(k)},
\end{equation}
where $\mu$ ranges over $k$-bounded $\mu\subset\la$ with $|\la|=|\mu|+1$, and $\mu^{\omega_k}\subset \la^{\omega_k}$. The index $r$ denotes a marking on $\mathfrak c(\la)$. Note that $\mathfrak c(\la)/\mathfrak c(\mu)$ is a union of $c$-ribbons of common height $h$, and
\[
\mathrm{spin}(\mu,r)=c(h-1)+N \le c(h-1)+(c-1)=ch-1,
\]
with $N$ the number of ribbons in rows strictly above the marked row $r$, hence $N\le c-1$ (see~\cite[p.~924]{BMPS}). Consequently,
\[
\mathrm{spin}(\mu,r) \le |\mathfrak c(\la)|-|\mathfrak c(\mu)|-1.
\]
Since $\we_1^{\perp}$ preserves the $q$-grading and acts nontrivially on each Schur component, iterating~\eqref{eqn:t1e} a total of $|\la|$ times yields that all $q$-degrees appearing in $\ws_\la^{(k)}$ are bounded above by
\[
|\mathfrak c(\la)|-|\la| = d_k(\la),
\]
and thus we obtain~\eqref{eqn:dlaest}.
\end{proof}

\begin{thm}\label{thm:CHid}
Fix $k>0$. For every $\la \in \Par_m^{(k)}$, the module $\ssc^{(k)}_\la$ is isomorphic to $M_{\la,(\la^{\omega_k})'}$ as defined in~\cite[4.2.1,~4.4.2]{Che10}. In particular, $M_{\la,(\la^{\omega_k})'}$ is the unique finite-dimensional graded $A_m$-module with head $L_\la$ and socle $\mathsf q^{d_k(\la)} L_{(\la^{\omega_k})'}$. Moreover, $L_{(\la^{\omega_k})'}$ occurs in $P_\la$ first in degree $d_k(\la)$, and with multiplicity one in that degree.
\end{thm}

\begin{proof}
Note that $\Psi [\la,k]$ is shallow by inspection. By Theorem~\ref{thm:Cat-sh} and Theorem~\ref{thm:HLss}, the module $\sk_\la$ has simple head and socle.

By Theorems~\ref{thm:ksch} and~\ref{thm:FKM}, the module $\ssc^{(k)}_\la$ has head $L_\la$ in degree $0$ and socle $\mathsf q^{d} L_{\mu}$ with $d \le d_k(\la)$. The quantity $d_k(\la)$ agrees with $d(\la,(\la^{\omega_k})')$ of~\cite[\S3.6]{Che10} by comparing Theorem~\ref{thm:dla},~\cite[Proposition~4.2.8(2)]{Che10}, and~\cite[Corollary~4.4.2]{Che10} (see also~\cite[Proposition~4.3.3(1)]{Che10}). Here $d(\la,(\la^{\omega_k})')$ is the minimal degree in which $L_{(\la^{\omega_k})'}$ can occur in an $A_m$-module generated by $L_\la$~\cite[Proposition~4.3.3(5)]{Che10}.

Applying the Schur--Weyl functor $\SW$ to the multiplicity statement in Theorem~\ref{thm:ksch} shows that the bound~\eqref{eqn:dlaest} is sharp, and that $[\ssc^{(k)}_\la : L_\mu]_q$ has top degree $d_k(\la)$ when $\mu=(\la^{\omega_k})'$. Uniqueness then follows from~\cite[Proposition~4.2.1]{Che10}.

The last assertion follows from~\cite[Theorem~4.2.2(3) and Lemma~4.2.4]{Che10}, since $d_k(\la)=d(\la,(\la^{\omega_k})')$ as observed above.
\end{proof}

\begin{cor}
Let $k \ge 1$. For every $\la \in \Par_m^{(k)}$, we have
\[
\ws^{(k)}_\la = \gch M_{\la,(\la^{\omega_k})'}.
\]
In particular, this verifies~\cite[Conjecture~4.4.3]{Che10}.
\end{cor}

\begin{proof}
Combining Theorem~\ref{thm:ksch}, Theorem~\ref{thm:FKM}, and Theorem~\ref{thm:CHid}, we obtain
\[
\ws^{(k)}_\la
= \gch \sk_\la
= \gch \ssc^{(k)}_\la
= \gch M_{\la,(\la^{\omega_k})'}.
\qedhere
\]
\end{proof}

\begin{ex}[Case $m=3$]\label{ex:n3}
When $m=3$,
\[
\Par_3^{(1)}=\{(1^3)\},\qquad
\Par_{3}^{(2)}=\{(21),(1^3)\},\qquad
\Par_{3}^{(3)}=\{(3),(21),(1^3)\}.
\]

We have
\[
\ws_{1^3}^{(1)}=\ws_{1^3}+(q+q^2)\ws_{21}+q^3\ws_3,\quad
\ws_{1^3}^{(2)}=\ws_{1^3}+q\ws_{21},\quad
\ws_{21}^{(2)}=\ws_{21}+q\ws_{3},
\]
while $\ws_\la^{(3)}=\ws_\la$ for $\la\in\{(3),(21),(1^3)\}$. Moreover,
\[
(1^3)^{\omega_2}=(21)=(21)',\qquad
(21)^{\omega_2}=(1^3)=(3)'.
\]

The preceding discussion implies
\[
\gch M_{(1^3),(21)}=\gch \ssc^{(2)}_{1^3}=\gch \mathtt{s}_{1^3}^{(2)},\qquad
\gch M_{(21),(3)}=\gch \ssc^{(2)}_{21}=\gch \mathtt{s}_{21}^{(2)}.
\]
\end{ex}

\begin{cor}\label{cor:kdual}
Let $k\in\Z_{>0}$. For every $\la \in \Par_m^{(k)}$,
\[
(\ssc^{(k)}_\la)^{\circledast} \cong \mathsf q^{d_k(\la)}\, (\ssc^{(k)}_{\la^{\omega_k}})^{\vee}.
\]
\end{cor}

\begin{proof}
We have $d_k(\la)=d_k(\la^{\omega_k})$ by Corollary~\ref{cor:dktr}. Applying $\circledast$ to the asserted isomorphism, it suffices to identify $\ssc^{(k)}_\la$ with $\mathsf q^{d_k(\la)}\bigl((\ssc^{(k)}_{\la^{\omega_k}})^{\vee}\bigr)^{\circledast}$. By Theorem~\ref{thm:CHid}, the module $\ssc^{(k)}_{\la^{\omega_k}}$ has head $L_{\la^{\omega_k}}$ and socle $\mathsf q^{d_k(\la)}L_{\la'}$, so that the latter module has head $L_\la$ and socle $\mathsf q^{d_k(\la)}L_{(\la^{\omega_k})'}$. Since $\la \in \Par^{(k)}_m$, the uniqueness statement in Theorem~\ref{thm:CHid} identifies it with $\ssc^{(k)}_\la$.
\end{proof}

\begin{cor}\label{cor:rigid}
Assume $n>m$ and let $k \in \Z_{>0}$. For each $\la \in \Par^{(k)}_m$, any finite-dimensional module $M \in \g[z]\gmod_m$ whose head is $V_\la^*$ and whose socle is $\mathsf q^{d_k(\la)} V^*_{(\la^{\omega_k})'}$ is necessarily isomorphic to $\sk_\la$.
\end{cor}

\begin{proof}
Transport the characterization of $M_{\la,(\la^{\omega_k})'}$ from~\cite[4.2.1,~4.4.2]{Che10} along $\SW \colon \g[z]\gmod_m \rightleftarrows A_m\gmod$ to conclude $M \cong \sk_\la$.
\end{proof}

\begin{cor}\label{cor:CHindep}
Let $n>m$ and view $\la \in \Par_m^{(k)}$ as an element of $\Par_m^{(k)}(n)$. Then
\[
\ssc_\la^{(k)} \cong \SW(\sk_\la),
\]
where $\SW$ denotes the Schur--Weyl duality functor for $\g=\mathfrak{gl}(n)$.
\end{cor}

\begin{proof}
By~\cite{BMPS}, we have $\mathsf{tr}_n ( \widetilde{s}^{(k)}_\la ) = \underline{\gch} \, \sk_\la$. 
Theorem~\ref{thm:ksch} yields the Schur expansion
\[
\widetilde{s}^{(k)}_\la = \sum_{\mu \in \Par_m,\ \ell(\mu) \le \ell(\la) \le n} c_{\la\mu}(q)\, \ws_\mu,
\]
so only partitions of size $m$ with at most $n$ parts appear. 
It follows that $\gch \sk_\la$ is independent of the choice of $n > \ell(\la)$.

Since the coefficients $\{c_{\la\mu}(q)\}_{\mu}$ are independent of $n > \ell(\la)$, the module $\sk_\la$ contains $V_\la^*$ in degree~$0$ and $V_{(\la^{\omega_k})'}^*$ in degree $d_k(\la)$. Here $\sk_\la$ has simple head $V_\la^*$ by Theorem~\ref{thm:Cat-sh}. 
Moreover, this degree calculation shows that $\mathsf q^{d_k(\la)} V_{(\la^{\omega_k})'}^*$ coincides with the simple socle of $\sk_\la$, as established in Theorem~\ref{thm:HLss}. Thus, Corollary~\ref{cor:rigid} implies the result.
\end{proof}

\section{A twisted vector bundle}\label{sec:tvb}

We continue to work in the setting of the previous section. Let $\{e_i^*\}_{i=1}^n$ and $\{e_i\}_{i=1}^n$ denote the standard bases of $V_1^* = ( \C^n )^*$ and $V_1$, respectively.  We set
\[
F_i V_1^* := \sum_{j > i} \C e_j^*, \qquad F_i V_1^*[z] := F_i V_1^* + V_1^* \otimes \C[z] z \subset V_1^* \otimes \C [z] \qquad (0 \le i \le n).
\]
It is easy to see that each $F_i V_1^*$ is a $\gb$-submodule of $V_1^*$, and each $F_i V_1^*[z]$ is a $\tb$-submodule of $V_1^* \otimes \C [z]$ on which $\g[z]$ acts by the loop extension of the $\g$-action on $V_1^*$. Dually, we set $F^i V_1 := \left( \coker ( F_i V_1^* \to V_1^* ) \right)^* = \sum_{j \le i} \C e_j \subset V_1$ for $0 \le i \le n$. Note that $F^i V_1 \subset V_1$ is the unique $B$-submodule of dimension $i$.

The main result of this section is the following:

\begin{thm}\label{thm:tvb}
Let $\Psi$ be a root ideal. There exists a $\widehat{G}[\![z]\!]$-equivariant vector bundle $\mathcal V$ on $\sX_{\Psi}$ with the following properties.
\begin{enumerate}
\item There exists a $\widehat{G}[\![z]\!]$-equivariant filtration
\[
\mathcal V = F^n \mathcal V \supset F^{n-1} \mathcal V \supset \cdots \supset F^0 \mathcal V = 0
\]
such that each successive quotient admits a $\widehat{G}[\![z]\!]$-equivariant isomorphism
\[
\cO_{\sX_{\Psi}} ( - \epsilon_i ) \cong F^{i} \mathcal V / F^{i-1} \mathcal V.
\]
Here $\cO_{\sX_{\Psi}}(-\epsilon_i)$ carries the equivariant structure specified in Theorem~\ref{thm:Cat}\textup{(3)}.
\item The restriction of $F^i \mathcal V$ to $T^*_\Psi X$ is isomorphic to $\pi^* ( G \times^B F^i V_1)$.
\end{enumerate}
\end{thm}

Below, we may denote $F^p\mathcal V \otimes_{\cO_{\sX_\Psi}} \cO_{\sX_\Psi} ( \la )$ by $F^p\mathcal V(\la)$ for simplicity.
The proof of Theorem~\ref{thm:tvb} is given at the end of this section. Recall that we have functors $\mathscr C_j ^{\Psi} ( m )$ for a root ideal $\Psi$ and $m \in \Z$ defined in~\eqref{eqn:Cjpdef}.

\begin{lem}\label{lem:twisting}
For a root ideal $\Psi$ and $m \in \Z$, we have an isomorphism of functors
\[
\mathscr C_j ^{\Psi} ( m ) \circ ( \C_{\La_{j-1}} \otimes \bullet ) \cong \C_{\La_{j-1}} \otimes \mathscr C_j ^{\Psi} ( m ) (\bullet) \qquad (1 \le j \le n),
\]
where the functor $\mathscr C_j ^{\Psi} ( m )$ is defined in~\eqref{eqn:Cjpdef}.
\end{lem}

\begin{proof}
This follows from the fact that the Demazure functors
 \[
\mathscr D_{j}, \ldots, \mathscr D_{n-1}, \mathscr D_{0}, \mathscr D_1, \ldots, \mathscr D_{\hs_j ( \Psi )-1},
\]
which constitute $\mathscr C_j ^{\Psi} ( m )$ up to character twists, commute with $\C_{\La_{j-1}}$ by Theorem~\ref{thm:Jos2}\textup{(2)} and $\hs_j (\Psi) < j$.
\end{proof}

\begin{lem}\label{lem:FV-pres}
For $0 \le i \le n$, the subspace
\[
F_i V_1^* + V_1^* \otimes \C[z] z \subset V_1^* \otimes \C[z,z^{-1}]
\]
is stable under $\tp_j$ for $j \not\equiv i \mod n$.
\end{lem}

\begin{proof}
The assertion for $j \neq 0$ is straightforward. The assertion for $j = 0$ follows as $E_{1n} \otimes z^{-1}$ acts nontrivially on $e_j^* \otimes z^r$ only when $j = 1$, and its image is $- e_n^* \otimes z^{r-1}$.
\end{proof}

For a root ideal $\Psi$, the homogeneous coordinate ring of $\sX_\Psi$ is given as
\begin{align}\nonumber
R_\Psi &:= \bigoplus_{\la \in \Par(n)} H^0 ( \sX_\Psi, \cO_{\sX_\Psi} ( \la ) )\\
& \cong \bigoplus_{m_1,\ldots,m_n \ge 0}  \left( \mathscr C_1^\Psi ( m_1 ; \bullet ) \circ \mathscr C_2^\Psi (m_2; \bullet) \circ \cdots \circ \mathscr C_n^\Psi (m_n; \bullet) \right) ( \C )^*,\label{eqn:Rpsi}
\end{align}
where $\la = \sum_i m_i \varpi_i$ and the last isomorphism is from Theorem~\ref{thm:catalan}. It follows that a $\Par(n)$-graded $R_\Psi$-module gives rise to a quasicoherent sheaf on $\sX_\Psi$, and it admits a $\widehat{G}[\![z]\!]$-action if the corresponding module admits a $\widehat{G}[\![z]\!]$-action.

For $0 \le i \le n$, we define
\begin{align}\nonumber
R_\Psi^{(i)} [\mathbf m] & := \left( \mathscr C_1^\Psi ( m_1 ; \bullet ) \circ \cdots \circ \mathscr C_i^\Psi (m_i; \bullet)  \circ \bigl( F_i V_1^* [z] \otimes \bullet \bigr) \circ \mathscr C_{i+1}^\Psi (m_{i+1}; \bullet)  \circ\cdots \circ \mathscr C_n^\Psi (m_n; \bullet) \right) ( \C )^{\vee},\\
R_\Psi^{(i)} & := \bigoplus_{\mathbf m = (m_j)_j \in ( \Z_{\ge 0} )^n}  R_\Psi^{(i)} [\mathbf m],\label{eqn:bigRpsi}
\end{align}
where $R_\Psi^{(0)}$ is understood to be $V_1 \otimes \C[z^{-1}] \otimes R_\Psi$ in which we understand that
\[
( V_1^*[z] )^{\vee} = V_1 \otimes \C[z^{-1}] \cong V_1 \otimes \C[z,z^{-1}] /  V_1 \otimes \C[z]z
\]
as graded $\g[z]$-modules.

\begin{prop}\label{prop:Rdeg}
For each $0 \le i \le n-1$ and $\mathbf m = (m_j)_j \in \Z_{\ge 0}^n$ such that $m_1 > m_2 > \cdots > m_n > 0$, we set $\la := \sum_i m_i \varpi_i$. We have a short exact sequence
\[
0 \rightarrow H^0 ( \sX_\Psi, \cO_{\sX_\Psi} ( \la - \epsilon_{i+1} ) ) \rightarrow R_\Psi^{(i)} [\mathbf m] \rightarrow R_\Psi^{(i+1)} [\mathbf m] \rightarrow 0
\]
of $\widehat{G}[\![z]\!]$-modules.
\end{prop}

\begin{proof}
We consider a short exact sequence
\[
0 \rightarrow F_{i+1} V_1^*[z] \rightarrow F_{i} V_1^*[z] \rightarrow \C_{-\epsilon_{i+1}} \rightarrow 0
\]
of graded $\tb$-modules. By Lemma~\ref{lem:FV-pres} and Theorem~\ref{thm:Jos2}\textup{(2)}, $\mathscr D_{i+1},\ldots, \mathscr D_{n-1}, \mathscr D_0, \ldots, \mathscr D_{\hs_{i+1} ( \Psi)-1}$, which constitute $\mathscr C_{i+1}^\Psi (m_{i+1} ; \bullet )$ up to character twists, preserve the space $F_{i} V_1^*[z]$. It follows that
\[
R_\Psi^{(i)} [\mathbf m] \cong \left( \mathscr C_1^\Psi ( m_1 ; \bullet ) \circ \cdots \circ \mathscr C_{i}^\Psi (m_{i}; \bullet)  \circ \mathscr C_{i+1}^\Psi (m_{i+1}; \bullet)  \circ \bigl( F_i V_1^* [z] \otimes \bullet \bigr) \circ \cdots \circ \mathscr C_n^\Psi (m_n; \bullet) \right) ( \C )^{\vee}.
\]
By the right exactness of the Demazure functors and the exactness of the character twist, we have a right exact sequence of $\tb$-modules
\begin{equation}
( R_\Psi^{(i+1)} [\mathbf m] ) ^{\vee} \rightarrow ( R_\Psi^{(i)} [\mathbf m] )^{\vee} \rightarrow \mathbb L^0 M \rightarrow 0,\label{eqn:Vfilt-D}
\end{equation}
where
\[
\mathbb L^{p} M := \mathbb L^{p} \left( \mathscr C_1^\Psi ( m_1 ; \bullet ) \circ \cdots \circ \mathscr C_{i+1}^\Psi (m_{i+1}; \bullet) \circ ( \C_{- \epsilon_{i+1}} \otimes \bullet ) \circ \cdots \circ \mathscr C_n^\Psi (m_n; \bullet) \right) ( \C ) \qquad (p \in \Z).
\]
The long exact sequence of derived functors shows that~\eqref{eqn:Vfilt-D} is short exact provided that $\mathbb L^{-1}M=0$.
Applying Lemma~\ref{lem:twisting} with $j=i+1$ puts the expression into a form to which Theorem~\ref{thm:catalan} applies. Thus Theorem~\ref{thm:catalan} gives
\[
\mathbb L^{p} M\cong \begin{cases} H^{-p} ( \sX_\Psi, \cO_{\sX_\Psi} ( \la - \epsilon_{i+1} ) )^* & (1 \le i < n)\\
H^{-p} ( \sX_\Psi, \cO_{\sX_\Psi} ( \la - \epsilon_{i+1} ) )^* \otimes \C_{\La_0} & (i=0)
\end{cases}.
\]
The factor $\C_{\La_0}$ in the case $i=0$ is trivial after restriction to the current algebra and loop rotation. Theorem~\ref{thm:Cat}\textup{(4)} gives $\mathbb L^{-1}M=0$, since $\la-\epsilon_{i+1}\in\Par(n)$ under the stated assumptions. Dualizing~\eqref{eqn:Vfilt-D} therefore gives the asserted short exact sequence.

The $\tb$-action on $R_\Psi^{(0)} [\mathbf m]$ and each $H^0 ( \sX_\Psi, \cO_{\sX_\Psi} ( \la - \epsilon_{i+1} ) )^*$ naturally extends to a $\widehat{G}[\![z]\!]$-action. A $\wgt$-semisimple $\tb$-action on $M$ extends to a $\widehat{G}[\![z]\!]$-action if the natural map $M \rightarrow \mathscr D_{w_0} ( M )$ is an isomorphism by Theorem~\ref{thm:kum}. In particular, we have
\[
\mathbb L^{\bullet} \mathscr D_{w_0}( H^0 ( \sX_\Psi, \cO_{\sX_\Psi} ( \la - \epsilon_{i+1} ) )^* ) \cong H^0 ( \sX_\Psi, \cO_{\sX_\Psi} ( \la - \epsilon_{i+1} ) )^*.
\]
From this, we deduce that
\[
0 \rightarrow \mathscr D_{w_0} \bigl( ( R_\Psi^{(i+1)} [\mathbf m] )^{\vee} \bigr) \rightarrow \mathscr D_{w_0} \bigl( ( R_\Psi^{(i)} [\mathbf m] )^{\vee} \bigr) \rightarrow H^0 ( \sX_\Psi, \cO_{\sX_\Psi} ( \la - \epsilon_{i+1} ) )^* \rightarrow 0.
\]
Thus, the $\tb$-action on $R_\Psi^{(i)} [\mathbf m]$ ($1 \le i \le n$) extends to a natural $\widehat{G}[\![z]\!]$-action by induction from the case $i=0$.
\end{proof}

Let us record one specific vanishing result which is useful in the sequel.

\begin{lem}\label{lem:shift-vanish}
Let $\Psi$ be a root ideal and let $\la = \sum_{j=1}^n m_j \varpi_j \in \sP$. Then we have
\begin{equation}
H^\bullet ( \sX_{\Psi}, \cO_{\sX_{\Psi}} (\la) ) = 0\label{eqn:totalXP}
\end{equation}
in each of the following cases:
\begin{enumerate}
\item $\la_i + 1 = \la_{i+1}$ for some $1 \le i < n$ such that $\Psi$ has no corner of shape $(i,\bullet)$ or $(\bullet,i+1)$;
\item $\la_n = -1$ and $\Psi \neq \emptyset$.
\end{enumerate}
\end{lem}

\begin{proof}
We regard the Demazure functors as total left derived functors, for which Theorem~\ref{thm:Jos}, Theorem~\ref{thm:Jos2}, Theorem~\ref{thm:kum}, and Lemma~\ref{lem:Dstab} remain valid (cf.~\cite[\S2]{Kat23a}), and compute~\eqref{eqn:totalXP} through Theorem~\ref{thm:catalan}. For $j \in \tI_\af$ and $\La \in \sP_\af$, the one-dimensional $\tb$-module $\C_\La$ extends to a $\tp_j$-module if $s_j \La = \La$, and $\mathbb L^\bullet \mathscr D_j ( \C_\La ) = 0$ if $s_j \La = \La + \al_j$ (cf.~\cite[\S 8.1]{Kum02}). Hence it suffices to find one stage of the composition at which the derived functor vanishes.

In case~\textup{(2)}, we have $m_n = -1$ and $\hs_n (\Psi) \ge 1$, so that $\mathscr C_n^\Psi (m_n; \C) = \mathscr D_{[0,\hs_n(\Psi))} ( \C_{-\La_n} )$ starts with $\mathscr D_0 ( \C_{-\La_n} ) = 0$, as $s_0 ( - \La_n ) = - \La_n + \al_0$ by $\La_n = \La_0 + \varpi_n$.

In case~\textup{(1)}, we have $m_i = -1$, and $h := \hs_i (\Psi) = \hs_{i+1} (\Psi) \le i-1$ since $\Psi$ has no corner of shape $(\bullet,i+1)$. Since $\Psi$ has no corner of shape $(i,\bullet)$, we have $\hs_r(\Psi) \neq i$ for all $r$, and hence there exists $i+1 \le \ell \le n$ such that $\hs_r (\Psi) \le i-1$ for $r \le \ell$ and $\hs_r (\Psi) \ge i+1$ for $r > \ell$. Then
\[
M := \mathbb L^\bullet \bigl( \mathscr C_{i+2}^\Psi (m_{i+2}; \bullet) \circ \cdots \circ \mathscr C_n^\Psi (m_n; \bullet) \bigr) ( \C )
\]
is a complex of $\tp_i$-modules: if $\ell<n$, the part $r > \ell$ is a complex of $\tp_{[\ell+1,\hs_{\ell+1}(\Psi))}$-modules by Lemma~\ref{lem:Dstab} and $i \in [\ell+1,\hs_{\ell+1}(\Psi))$; if $\ell=n$, this part is $\C$. For $i+1 < r \le \ell$ the functor $\mathscr C_r^\Psi (m_r;\bullet)$ involves only $\C_{m_r \La_r} \otimes \bullet$ and $\mathscr D_l$ with $l \notin \{i-1,i,i+1\}$, which preserve $\tp_i$-modules by Theorem~\ref{thm:Jos2}\textup{(2)} and Theorem~\ref{thm:kum}. Since $s_i \La_{i+1} = \La_{i+1}$, the complex $M' := \C_{m_{i+1} \La_{i+1}} \otimes M$ is again a complex of $\tp_i$-modules. As $s_l \La_i = \La_i$ for $l \neq i$ and $i \notin [i+1,h)$, we obtain
\begin{align*}
\mathbb L^\bullet \mathscr C_{i}^\Psi \bigl( m_{i}; \mathbb L^\bullet \mathscr C_{i+1}^\Psi (m_{i+1}; M ) \bigr)
&= \mathbb L^\bullet \mathscr D_{[i,h)} \circ \mathbb L^\bullet \mathscr D_{[i+1,h)} ( \C_{-\La_i} \otimes M' )\\
&\cong \mathbb L^\bullet \mathscr D_{[i,h-1)} \circ \mathbb L^\bullet \mathscr D_{[i+1,h)} \circ \mathbb L^\bullet \mathscr D_{i} ( \C_{-\La_i} \otimes M' ),
\end{align*}
where the second isomorphism is Theorem~\ref{thm:Jos} (cf.~\cite[Corollary~2.9]{Kat23a}). Finally, $\mathbb L^\bullet \mathscr D_i ( \C_{-\La_i} \otimes M' ) \cong M' \otimes \mathbb L^\bullet \mathscr D_i ( \C_{-\La_i} ) = 0$ by Theorem~\ref{thm:Jos2}\textup{(2)}, since $s_i ( - \La_i ) = - \La_i + \al_i$. This proves~\eqref{eqn:totalXP}.
\end{proof}

\begin{proof}[Proof of Theorem~\ref{thm:tvb}]
By Corollary~\ref{cor:Cat-sg} and Corollary~\ref{cor:inclCat}, it is enough to prove the assertion for $\Psi = \Delta^+$ and henceforth we omit $\Psi$ from notation (during this proof).

We first prove~\textup{(1)}.
Let $f:\widetilde{\sX}\to\sX_{\Delta^+}$ be the Bott--Samelson--Demazure--Hansen resolution associated with the sequence of Demazure functors in~\eqref{eqn:bigRpsi}. The coefficient module $F_iV_1^*[z]$ defines a graded coefficient sheaf on this resolution, and the intervening character twists correspond to line bundle twists~\cite[Proposition~8.1.17]{Kum02}. Thus $R^{(i)}$ is the multigraded section module of the corresponding pushforward. Multiplication of sections makes it an $R$-module, compatibly with the graded $\g[z]$-action.

The maps in Proposition~\ref{prop:Rdeg} are induced by the short exact
sequences of coefficient modules used in its proof. Tensoring these
sequences with line bundles and multiplying sections commute.
For $m_1>\cdots>m_n>0$, we have identifications
\[
\ker\bigl(R^{(i-1)}[\mathbf m]\to R^{(i)}[\mathbf m]\bigr)
\cong H^0\bigl(\sX_{\Delta^+},
\cO_{\sX_{\Delta^+}}(\la-\epsilon_i)\bigr),
\qquad \la=\sum_jm_j\varpi_j,
\]
as graded $\g[z]$-modules.
By Theorem~\ref{thm:HLss}, \eqref{eqn:stdmult}, and
\cite[Chap.~III,~(6.5)]{Mac95}, the highest-degree part of
$\HL_\mu$ ($\mu \in \Par(n)$) is its simple socle
$\mathsf q^{\mathsf n(\mu)}V_{|\mu|}^*$.
Since $\mathsf n(\mu+\nu)=\mathsf n(\mu)+\mathsf n(\nu)$
for $\mu,\nu\in\Par(n)$, the Littlewood--Richardson rule implies
that a nonzero graded $\g[z]$-module map
\[
\HL_\mu\otimes\HL_\nu\to\HL_{\mu+\nu}
\]
is unique up to scalar.
Thus the multiplication maps from the two constructions are proportional.

To pass to sheaves, choose integers $a_1>\cdots>a_n>0$ and set
$\eta=\sum_ja_j\varpi_j$.
The line bundle $\cO(\eta)$ is ample, and every fixed multidegree,
after adding $r\eta$ for all sufficiently large $r$, lies in the range above.
Along this ray, normalize the successive identifications to respect
multiplication by the degree-$\eta$ component of $R$; associativity
then gives compatibility with all components of degree $p\eta$.
The Serre correspondence therefore identifies the sheaf associated with
$\ker(R^{(i-1)}\to R^{(i)})$ with
$\cO_{\sX_{\Delta^+}}(-\epsilon_i)$, equivariantly.

Applying Proposition~\ref{prop:Rdeg}, we conclude that the $\widehat{G}[\![z]\!]$-equivariant quasicoherent sheaf $\mathcal V$ on $\sX_{\Delta^+}$ associated to $\ker ( R^{(0)} \to R^{(n)} )$ admits the desired filtration by setting $F^i \mathcal V$ to the sheaf corresponding to $\ker ( R^{(0)} \to R^{(i)} )$. In particular, $\mathcal V$ is a vector bundle of rank $n$ since it admits an $n$-step filtration by line bundles $\{\cO_{\sX_{\Delta^+}}(-\epsilon_i)\}_i$.

We prove~\textup{(2)}.
The fibers over $B/B \in X$ of the sheaves associated to $R^{(0)}$ and $R^{(n)}$ are $V_1[z^{-1}]$ and $V_1[z^{-1}]z^{-1}$ as graded $B$-modules, respectively. The map between these fibers is the natural projection that discards the constant term. It follows that the fiber of $\mathcal V \subset V_1[z^{-1}] \otimes \cO_{\sX_{\Delta^+}}$ over $B/B$ surjects onto $V_1 = \ker ( V_1[z^{-1}] \rightarrow V_1[z^{-1}]z^{-1} )$ by the right exactness of the fiber functor. As $\rk \mathcal V = n$, the fiber is isomorphic to $V_1$. 

The fiber of $F^i \mathcal V \subset \mathcal V$ at $B/B$ defines an $i$-dimensional $B$-submodule of $V_1$, which must be $F^i V_1$. Since $\C[\pi^{-1}(B/B)]$ is nonnegatively graded with degree zero part $\C$, and the fiber of $F^i \mathcal V$ at $B/B$ is concentrated in degree zero, the graded Nakayama lemma yields $\Gamma ( \pi^{-1}(B/B), F^i \mathcal V ) \cong \C[\pi^{-1}(B/B)] \otimes F^i V_1$ as graded $B$-equivariant modules. Applying the $G$-action identifies $F^i \mathcal V$ with $\pi^* ( G \times ^B F^i V_1 )$.
\end{proof}

\section{Geometric Pieri rule}\label{sec:gpr}

We retain the setting of the previous section.

This section provides geometric counterparts of several statements from~\cite[\S7.2]{BMPS}. In particular, familiarity with the examples given there and in~\cite{BMPS2} will help the reader appreciate the material more readily.

We begin by recalling a basic observation that will be used repeatedly.

\begin{thm}[{\cite[Proposition~5.12]{Kat23a}}]\label{thm:Psi-rest}
Let $\Psi$ and $\Psi'$ be root ideals with $\gn(\Psi)\subset \gn(\Psi')$, and let $\la \in \sP$. For each $i\in\Z$, there is a natural restriction map
\[
H^i(\sX_{\Psi'}, \cO_{\sX_{\Psi'}}(\la)) \longrightarrow H^i(\sX_{\Psi}, \cO_{\sX_{\Psi}}(\la)),
\]
and this map is a homomorphism of graded $\g[z]$-modules.
\end{thm}

\begin{prop}
Let $\Psi$ be a root ideal, and let $\Psi'$ be obtained from $\Psi$ by removing a corner $(i,j)$. Then the closed embedding $\sX_{\Psi'} \hookrightarrow \sX_{\Psi}$ is a $\widehat{G}[\![z]\!]$-equivariant smooth Cartier divisor.
\end{prop}

\begin{proof}
The construction in~\cite[\S3]{Kat23a} presents $\sX_\Psi$ as an iterated equivariant projective-space bundle. Removing the corner $(i,j)$ reduces $\hs_j(\Psi)$ by one and leaves the other $\hs_r(\Psi)$ unchanged. At the column-$j$ stage of the tower~\cite[Proposition~3.7]{Kat23a}, the base $\P^{\hs_j(\Psi)+n-j}$ is replaced by its $\tB$-stable hyperplane, while the subsequent fiber is unchanged. Pulling back this hyperplane and carrying it through the preceding bundle stages gives $\sX_{\Psi'}$. Hence the given equivariant closed embedding is a smooth Cartier divisor.
\end{proof}

\begin{thm}\label{thm:gmaps}
Let $\Psi$ be a root ideal and $\la \in \sP$. If $(i,j)$ is a corner of $\Psi$, then there is a natural left exact sequence of graded $\g[z]$-modules
\[
0 \rightarrow \mathsf q H^0 \bigl(T^*_\Psi X, \cO_{T^*_\Psi X}(\la+\epsilon_i-\epsilon_j)\bigr) \longrightarrow H^0 \bigl(T^*_\Psi X, \cO_{T^*_\Psi X}(\la)\bigr)\longrightarrow H^0\bigl(T^*_{\Psi\setminus(i,j)} X, \cO_{T^*_{\Psi\setminus(i,j)} X}(\la)\bigr),
\]
which restricts to a left exact sequence of graded $\g[z]$-modules
\[
0 \rightarrow \mathsf q H^0 \bigl(\sX_\Psi, \cO_{\sX_\Psi}(\la+\epsilon_i-\epsilon_j)\bigr) \longrightarrow H^0 \bigl(\sX_\Psi, \cO_{\sX_\Psi}(\la)\bigr) \longrightarrow\ H^0 \bigl(\sX_{\Psi\setminus(i,j)}, \cO_{\sX_{\Psi\setminus(i,j)}}(\la)\bigr).
\]
\end{thm}

\begin{proof}
We set $\Psi':=\Psi\setminus(i,j)$. By Corollary~\ref{cor:inclCat}, the natural surjection $\cO_{\sX_\Psi}(\la) \twoheadrightarrow \cO_{\sX_{\Psi'}}(\la)$ restricts along $T^*_{\Psi'}X\subset T^*_{\Psi}X$ to give
\[
\cO_{T^*_\Psi X}(\la) \longrightarrow \cO_{T^*_{\Psi'} X}(\la).
\]
This fits into a short exact sequence on $T^*_\Psi X$
\begin{equation}\label{eqn:SESTX}
0 \to \mathsf{q} \cO_{T^*_\Psi X}(\la+\epsilon_i-\epsilon_j)\to \cO_{T^*_\Psi X}(\la) \to \cO_{T^*_{\Psi'} X}(\la) \to 0,
\end{equation}
coming from the relative hyperplane section corresponding to the corner $(i,j)$.

Let $D$ denote the $\widehat{G}[\![z]\!]$-equivariant Cartier divisor on $\sX_\Psi$ cut out by $\sX_{\Psi'}$ with multiplicity one. Then~\eqref{eqn:SESTX} extends to
\[
0 \to \cO_{\sX_\Psi}(\la-D)\to \cO_{\sX_\Psi}(\la)\to \cO_{\sX_{\Psi'}}(\la) \to 0.
\]
By the Picard group computation~\cite[Corollary~4.9]{Kat23a}, line bundles on $\sX_\Psi$ are classified by $\sP$ and this classification is compatible with the restriction to $X\subset T^*_\Psi X\subset \sX_\Psi$ (Corollary~\ref{cor:inclCat}). Since $T^*_{\Psi'} X \subset T^*_\Psi X$ defines a $G$-equivariant divisor corresponding to $\cO_{T^*_\Psi X}(\epsilon_i-\epsilon_j)$, we have a canonical isomorphism (here $\mathsf{q}$ denotes the degree one $\Gm$-character twist)
\[
\cO_{\sX_\Psi}(\la-D)\simeq \mathsf q \cO_{\sX_\Psi}(\la+\epsilon_i-\epsilon_j + a \varpi_n ) \otimes \C_{a \varpi_n}
\]
for some $a \in \Z$. Here $a \in \Z$ is detected by its restriction to $\P^1 \cong Y\subset \sX_\Psi$ used in the proof of~\cite[Corollary~4.9]{Kat23a}, that is a subspace of the last column stage in the construction of $\sX_\Psi$ (\cite[\S3]{Kat23a}). Here $D$ is pulled back from the base space of the column-$j$ stage of the tower.
Thus $\deg\bigl(\cO(-D)|_Y\bigr)=-\delta_{jn}$. Thus, the preceding isomorphism gives $a=0$.

Passing to global sections yields the claimed left exact sequences. They are graded $\g[z]$-module maps by the derivations of the $\widehat{G}[\![z]\!]$-action on $\sX_{\Psi'}\subset \sX_\Psi$, that descend to the derivation action on functions on $T^*_{\Psi'} X\subset T^*_{\Psi} X$.
\end{proof}

Here we give a geometric interpretation of~\cite[Lemma~6.2]{BMPS}.

\begin{lem}\label{lem:triv}
Let $\Psi$ be a root ideal and $\la \in \sP$ with $\la_i=\la_{i+1}$. Assume that $\Psi$ has no corner of the form $(\bullet,i+1)$ or $(i,\bullet)$. Then
\[
H^{\bullet}\bigl(T^*_\Psi X, \cO_{T^*_\Psi X}(\la-\epsilon_i)\bigr)=0.
\]
\end{lem}

\begin{proof}
The assumption ensures that $\gn(\Psi)$ is $\SL^{+}(2,i)$-stable, or equivalently, stable under the action of $\mathfrak{sl}(2,i)$. Hence, there is a $G$-equivariant fibration
\[
\eta: T^*_\Psi X=G\times^B \gn(\Psi) \longrightarrow G\times^{P_i}\gn(\Psi),
\]
with fibers $\P^1\simeq P_i/B$. The line bundle $\cO_{T^*_\Psi X}(\la-\epsilon_i)$ restricts on each fiber to
\[
\cO_{\P^1}\bigl(\langle \la-\epsilon_i,\alpha_i^\vee\rangle\bigr)
=\cO_{\P^1}(-1),
\]
since $\langle \la,\alpha_i^\vee\rangle=\la_i-\la_{i+1}=0$. Hence $\R^s\eta_*\cO_{T^*_\Psi X}(\la-\epsilon_i)=0$ for all $s$ by cohomology and base change for projective bundles. The Leray spectral sequence then yields the desired vanishing.
\end{proof}

\begin{thm}[Geometric Pieri rule]\label{thm:assgr}
Let $\Psi$ be a root ideal and $\la \in \sP^+$. Then the graded $\g[z]$-module $H^{0}\bigl(T^*_\Psi X, \mathcal V (\la)\bigr)$ admits a
filtration whose associated graded is a direct sum of graded $\g[z]$-modules
\[
H^{0}\bigl(T^*_\Psi X, \cO_{T^*_\Psi X}(\la-\epsilon_i)\bigr)\qquad (1\le i\le n),
\]
provided that
\begin{equation}\label{eqn:hcvan}
H^{r}\bigl(T^*_\Psi X, \cO_{T^*_\Psi X}(\la-\epsilon_i)\bigr)=0\qquad\text{for all}\quad r>0\text{ and }1\le i\le n.
\end{equation}
\end{thm}

\begin{proof}
By Theorem~\ref{thm:tvb}\textup{(1)}, for $0 \le i < n$ we have a short exact sequence
\[
0 \rightarrow \cO_{\sX_\Psi} (\la - \epsilon_{i+1}) \rightarrow ( \mathcal V / F^i \mathcal V ) (\la) \rightarrow ( \mathcal V / F^{i+1} \mathcal V ) (\la) \rightarrow 0
\]
of $\widehat{G}[\![z]\!]$-equivariant sheaves on $\sX_\Psi$. Restricting to $T^*_\Psi X$ and taking global sections, we obtain a short exact sequence of graded $\g[z]$-modules by~\eqref{eqn:hcvan}. Since $\mathcal V / F^0 \mathcal V = \mathcal V$ and $\mathcal V / F^n \mathcal V = 0$, d\'evissage yields the assertion.
\end{proof}

\begin{prop}\label{prop:SW-res}
Assume $n > m$. Let $\Psi$ be a root ideal and $\la \in \Par_m(n)$. Assume that the left-hand side of
\[
H^{0}\bigl(\sX_\Psi, \mathcal V (\la)\bigr) \hookrightarrow V_1 \otimes H^{0}\bigl(T^*_\Psi X, \cO_{T^*_\Psi X}(\la)\bigr)
\] 
is precisely the span of all polynomial $G$-submodules of the right-hand side. Then, we have
\begin{equation}
H^{0}\bigl(\sX_\Psi, \mathcal V (\la)\bigr) \cong \WS \Bigl( \mathrm{Res}_{A_{m-1}}^{A_m} \bigl( \SW ( H^{0}\bigl(\sX_\Psi, \cO_{\sX_\Psi}(\la)\bigr)  ) \bigr)\Bigr).\label{eqn:Sm-res}
\end{equation}
\end{prop}

\begin{proof}
By the Pieri rule, we have
\[
V_1 \otimes V^{*}_\mu \cong \bigoplus_{\nu = \mu - \epsilon_i \in \sP^+} V_\nu^* \qquad\text{and}\qquad \mathrm{Res}_{\mathfrak S_{m-1}}^{\mathfrak S_m} L _\mu \cong \bigoplus _{\nu = \mu - \epsilon_i \in \Par_{m-1}} L_\nu.
\]
The indices common to the two sums are precisely those $\nu=\mu-\epsilon_i$ that lie in $\Par_{m-1}(n)=\Par_{m-1}$. By Theorem~\ref{thm:FKM} and assumption, we deduce that both sides of~\eqref{eqn:Sm-res} are abstractly isomorphic as graded $G$-modules.

By construction, $\mathcal V$ is a $\widehat{G}[\![z]\!]$-equivariant vector subbundle of $V_1[z^{-1}] \otimes \cO_{\sX_\Psi}$. This implies
\begin{equation}
H^{0}\bigl(\sX_\Psi, \mathcal V (\la)\bigr) \subset V_1 [z^{-1}] \otimes H^{0}\bigl(\sX_\Psi, \cO_{\sX_\Psi}(\la)\bigr)\label{eqn:max-Res}
\end{equation}
as graded $\g[z]$-modules. The left-hand side of~\eqref{eqn:max-Res} belongs to $\g[z]\gmod_{m-1}$. Here we have
\begin{align*}
\SW ( H^{0}\bigl(\sX_\Psi, \cO_{\sX_\Psi}(\la)\bigr) ) & = \hom_{\g[z]} ( ( V_1^* [z] )^{\otimes m}, H^{0}\bigl(\sX_\Psi, \cO_{\sX_\Psi}(\la)\bigr)) \\
& \cong \hom_{\g[z]} ( ( V_1^* [z] )^{\otimes (m-1)}, V_1 [z^{-1}] \otimes H^{0}\bigl(\sX_\Psi, \cO_{\sX_\Psi}(\la)\bigr)).
\end{align*}
Restriction to $A_{m-1}$ corresponds to its action on the first $m-1$ tensor factors, namely on $(V_1^*[z])^{\otimes(m-1)}$ in the last term.

Applying $\WS$, the right-hand side of~\eqref{eqn:Sm-res} is the sum of graded $\g[z]$-submodules in $V_1 [z^{-1}] \otimes H^{0}\bigl(\sX_\Psi, \cO_{\sX_\Psi}(\la)\bigr)$ belonging to $\g[z]\gmod_{m-1}$. Hence,~\eqref{eqn:Sm-res} is an inclusion. Therefore, both sides of~\eqref{eqn:Sm-res} define the same subspace in~\eqref{eqn:max-Res} as they have the same character.
\end{proof}

\begin{prop}\label{prop:pol}
Suppose that $M$ is a graded $\g[z]$-module whose underlying $G$-module is rational, equipped with a finite filtration
\[
0 = F_0 M  \subset F_1 M \subset \cdots \subset F_r M = M
\]
such that, for each $1 \le i \le r$, we have
\[
F_i M / F_{i-1} M \cong \mathsf{q}^{s_i} H^0 ( T^*_{\Psi[\mu^{(i)}, k]} X, \cO_{T^*_{\Psi[\mu^{(i)}, k]} X} ( \mu^{(i)})) \qquad s_i \in \Z, \mu^{(i)} \in \sP^+
\]
whose $\g[z]$-module structure is inherited from $\cO_{\sX_{\Psi[\mu^{(i)}, k]}} ( \mu^{(i)} )$. If $N \subset M$ is the span of all irreducible polynomial $G$-submodules which is also a graded $\g[z]$-submodule, then $N$ admits a finite filtration as graded $\g[z]$-modules whose associated graded is $\bigoplus_{\mu^{(i)} \in \Par(n)} \mathsf{q}^{s_i} \sk_{\mu^{(i)}}$.
\end{prop}

\begin{proof}
Set $F_iN:=F_iM\cap N$. Since $G$ is reductive, taking the sum of the polynomial isotypic components is exact on rational $G$-modules. Thus the natural injection
\[
F_iN/F_{i-1}N\hookrightarrow F_iM/F_{i-1}M
\]
identifies its source with the polynomial part of its target. If $\mu^{(i)}\in\Par(n)$, this polynomial part is $\mathsf q^{s_i}\sk_{\mu^{(i)}}$ by Theorem~\ref{thm:Cat}\textup{(4)}. If $\mu^{(i)}\notin\Par(n)$, it is zero by Theorem~\ref{thm:tri}.
Since $N$ and every $F_iM$ are graded $\g[z]$-submodules, these identifications respect the graded $\g[z]$-actions. Omitting the zero quotients gives the required finite filtration.
\end{proof}

\section{$k$-Schur straightening law}\label{sec:kst}

We retain the setting of the previous section.

\begin{lem}\label{lem:ctwist}
Let $\Psi$ be a root ideal, and let $\la \in \sP$ satisfy $\la_i=\la_{i+1}$. Assume that $\Psi$ has no corner of the form $(\bullet,i+1)$ or $(i,\bullet)$, and that we can delete one of the corners of the form $(\bullet,i)$ and $(i+1,\bullet)$ to obtain smaller root ideals $\Psi_c$ and $\Psi_r$, respectively. Then there are graded $\g[z]$-module isomorphisms
\[
H^{\bullet}\bigl(T^*_{\Psi_c} X, \cO_{T^*_{\Psi_c} X}(\la)\bigr) \cong
H^{\bullet}\bigl(T^*_{\Psi} X, \cO_{T^*_{\Psi} X}(\la)\bigr) \cong
H^{\bullet}\bigl(T^*_{\Psi_r} X, \cO_{T^*_{\Psi_r} X}(\la)\bigr).
\]
\end{lem}

\begin{proof}
Set $\Psi_{0} := \Psi_c \cap \Psi_r$. Then the subspaces $\gn(\Psi_0) \subset \gn(\Psi)$ are $\SL^{+}(2,i)$-stable.

Consider the maps
\[
\pi_c : T^*_{\Psi_c} X \longrightarrow G \times^{P_i} \gn(\Psi),
\qquad
\pi_r : T^*_{\Psi_r} X \longrightarrow G \times^{P_i} \gn(\Psi),
\]
and also $\eta : T^*_{\Psi} X \longrightarrow G \times^{P_i} \gn(\Psi)$. The maps $\pi_{c}$ and $\pi_{r}$ are $G$-inflations of the products of $\gn (\Psi_0)$ and
\begin{equation}
\SL^{+}(2,i) \times^{B_i} \bigl( \begin{psmallmatrix} \C \\ 0 \end{psmallmatrix} \oplus \C^2 \bigr) \longrightarrow \C^2 \oplus \C^2.\label{eqn:SL2map}
\end{equation}
Thus, they are Zariski-locally the product of the identity on $G \times^{P_i} \gn(\Psi_0)$ with~\eqref{eqn:SL2map}. Since the line bundle $\cO_{T^*_\Psi X}(\la)$ is trivial along these fibers (because $\la_i=\la_{i+1}$), we deduce
\[
\R^{\bullet}(\pi_c)_*\cO_{T^*_{\Psi_c} X}(\la) \cong
\R^{\bullet}(\eta)_*\cO_{T^*_{\Psi} X}(\la) \cong
\R^{\bullet}(\pi_r)_*\cO_{T^*_{\Psi_r} X}(\la)
\]
by the standard $\SL(2)$-calculation. This implies the desired isomorphisms after taking cohomology. Since each isomorphism arises from a natural map, their composition yields an isomorphism of graded $\g[z]$-modules.
\end{proof}

\begin{lem}\label{lem:atwist}
Let $\Psi_0$ be a root ideal, and let $\la \in \sP$ satisfy $\la_i+1=\la_{i+1}$. Assume that $\Psi_0$ has no corner of the form $(\bullet,i+1)$ or $(i,\bullet)$, and let
\[
\Psi := \Psi_0 \cup (j,i+1), \qquad \Psi^+ := \Psi \cup (j,i)
\]
be root ideals. Then there is an isomorphism of graded $\g[z]$-modules
\[
H^{\bullet}\bigl(T^*_{\Psi} X, \cO_{T^*_{\Psi} X}(\la)\bigr)
\cong \mathsf q H^{\bullet}\bigl(T^*_{\Psi^+} X, \cO_{T^*_{\Psi^+} X}(\la+\epsilon_j-\epsilon_{i+1})\bigr).
\]
\end{lem}

\begin{proof}
We have a short exact sequence
\[
0 \longrightarrow \mathsf{q} \cO_{T^*_{\Psi} X}(\la+\epsilon_j-\epsilon_{i+1})
\longrightarrow \cO_{T^*_{\Psi} X}(\la)
\longrightarrow \cO_{T^*_{\Psi_0} X}(\la) \longrightarrow 0.
\]
By Lemma~\ref{lem:triv}, we have $H^{\bullet}\bigl(T^*_{\Psi_0} X, \cO_{T^*_{\Psi_0} X}(\la)\bigr)\equiv 0$. It follows that
\[
H^{\bullet}\bigl(T^*_{\Psi} X, \cO_{T^*_{\Psi} X}(\la)\bigr)
\cong \mathsf q H^{\bullet}\bigl(T^*_{\Psi} X, \cO_{T^*_{\Psi} X}(\la+\epsilon_j-\epsilon_{i+1})\bigr),
\]
where the grading shift arises from the construction of the map in the short exact sequence, which corresponds to multiplication by a local section of degree one. Now apply Lemma~\ref{lem:ctwist} to conclude that
\[
H^{\bullet}\bigl(T^*_{\Psi^+} X, \cO_{T^*_{\Psi^+} X}(\la+\epsilon_j-\epsilon_{i+1})\bigr)
\cong
H^{\bullet}\bigl(T^*_{\Psi} X, \cO_{T^*_{\Psi} X}(\la+\epsilon_j-\epsilon_{i+1})\bigr).
\]
These isomorphisms are induced by natural maps and hence are isomorphisms of graded $\g[z]$-modules, as required.
\end{proof}

\begin{defn}
Let $k\in\Z_{\ge 1}$. A pair $(\Psi,\la)$, where $\Psi$ is a root ideal and $\la \in \sP$, is \emph{$(k,i)$-adapted} if
\[
\la_1\ge \la_2\ge\cdots\ge \la_i
\quad\text{and}\quad
E_{jl}\in\gn(\Psi)\ \Longleftrightarrow\ l> \mathsf o^k_\la(j)\ \ \ (1\le j\le i,\ j<l\le n),
\]
where $\mathsf o^k_\la(j)$ is defined in~\eqref{eqn:okdef}. The second condition can equivalently be written as
\[
\gn(\Psi)\cap \bigoplus_{j=1}^i\bigoplus_{l>j}\C E_{jl}
=
\bigoplus_{j=1}^i\bigoplus_{\mathsf o^k_\la(j)<l\le n}\C E_{jl}.
\]
\end{defn}

\noindent
Thus $(k,i)$-adaptedness requires that $\la_1\ge\cdots\ge\la_i$ and that the first $i$ full rows of $\Psi$ agree with those prescribed by $\mathsf o^k_\la$. In particular, restriction to the $i\times i$ upper-left block gives $\Phi=\Psi[\mu,k]$ for $\g=\gl(i)$, where $\mu=(\la_1,\ldots,\la_i)$; this restricted equality alone does not specify the entries to the right of that block.

In the following proposition and its proof, we use the convention $\la_0=+\infty$ in comparisons. Thus a corner chain stops as soon as it reaches row $1$.

\begin{prop}\label{prop:upshifts}
Let $k \in \Z_{\ge 1}$ and let $(\Psi,\la)$ be a $(k,i)$-adapted pair such that $\la_{i+1}=\la_i+1$, $\Psi$ has no corner of the form $(i,\bullet)$, and, whenever $E_{i+1,j}\in \gn(\Psi)$ for some $j>i+1$, the root ideal $\Psi$ has a corner of the form $(i+1,\bullet)$. Then:
\begin{enumerate}
\item If $\Psi$ has no corner of the form $(\bullet,i+1)$, then
\[
H^{\bullet}\bigl(T^*_\Psi X,\cO_{T^*_\Psi X}(\la)\bigr)\equiv 0.
\]
\item If $\Psi$ has a corner of the form $(\bullet,i+1)$, find such a corner $(i_1,i+1)$ and form a maximal sequence of corners
\[
(i_1,i_0):=(i_1,i+1),\ (i_2,i_1),\ \ldots,\ (i_\ell,i_{\ell-1})
\]
with $i_\ell<\cdots<i_2<i_1<i$ and $\la_{i_s-1}=\la_{i_s}$ for all $1\le s<\ell$. If $\la_{i_\ell-1}=\la_{i_\ell}$, then
\[
H^\bullet\bigl(T^*_\Psi X,\cO_{T^*_\Psi X}(\la)\bigr)=0.
\]
Otherwise $\la_{i_\ell-1}>\la_{i_\ell}$, and there exists a root ideal $\Phi$ with the following properties:
\begin{itemize}
\item $\gn(\Phi)$ agrees with $\gn(\Psi)$ in all rows of index $>i+1$;
\item $(\Phi,\la+\epsilon_{i_\ell}-\epsilon_{i+1})$ is $(k,i+1)$-adapted;
\item there is an isomorphism of graded $\g[z]$-modules
\[
H^\bullet\bigl(T^*_\Psi X,\cO_{T^*_\Psi X}(\la)\bigr)
\cong\mathsf q^\ell H^\bullet\bigl(T^*_\Phi X,\cO_{T^*_\Phi X}(\la+\epsilon_{i_\ell}-\epsilon_{i_0})\bigr).
\]
\end{itemize}
\end{enumerate}
\end{prop}

\begin{ex}[$n=5$, $k=2$]\label{n5ex1}
We illustrate how Proposition~\ref{prop:upshifts}\textup{(2)} works in the case $k=2, \la=(11112)$ using the graded characters. Here $i=4$ since $\la_4+1 = \la_5 = 2$ and $\Psi [\la,2]$ has no corner of shape $(4,\bullet)$.

For each $\la \in \sP$, we set
$$\mathrm{HL}^\Psi_\la := \ugch H^0 ( T^*_\Psi X, \cO_{T^*_\Psi X} ( \la )).$$
We set $\Psi = \Psi [ \la, 2 ]$ by putting $\la = ( \la_1,\la_2,\la_3,\la_4,\la_5)$ in the diagonal entries, and indicate the shape of $\Psi$ as red boxes. In other words, we have
\[
\mathrm{HL}^{\Psi[11112,2]}_{11112} = \CatalanInlineRC{1,1,1,1,2}{1}{1/3,1/4,1/5,2/4,2/5,3/5}.
\]
The algorithm in Proposition~\ref{prop:upshifts} proceeds as follows:
\begin{align}
\CatalanInlineRC{1,1,1,1,2}{1}{1/3,1/4,1/5,2/4,2/5,3/5} & = \mathsf q \CatalanInlineRC{1,1,2,1,1}{1}{1/3,1/4,1/5,2/4,2/5,3/4,3/5} = \mathsf q^2 \CatalanInlineRC{2,1,1,1,1}{1}{1/2,1/3,1/4,1/5,2/4,2/5,3/4,3/5} = \mathsf q^2 \CatalanInlineRC{2,1,1,1,1}{1}{1/2,1/3,1/4,1/5,2/4,2/5,3/5}.
\end{align}
Here, the sequence of corners corresponds to indices $i_1=3$ and $i_2=1$ (so $\ell = 2$).
\end{ex}

\begin{proof}[Proof of Proposition~\ref{prop:upshifts}]
We first prove~\textup{(1)}. Since $\Psi$ has no corner of the form $(\bullet,i+1)$, Lemma~\ref{lem:triv} applies, with $\mu_i=\mu_{i+1}$ for the weight $\mu = \la + \epsilon_i$ and the stated absence of corners, and yields
\[
H^{\bullet}\bigl(T^*_\Psi X,\cO_{T^*_\Psi X}(\la)\bigr)\equiv 0.
\]
This proves~\textup{(1)}.

\medskip

We now prove~\textup{(2)}. Suppose $\Psi$ has a corner of the form $(\bullet,i_0)$.
Find such a corner $(i_1,i_0)$ with $i_1<i$.
Since $(\Psi,\la)$ is $(k,i)$-adapted, $\Phi_1 := \Psi \cup (i_1,i_0-1)$ is a root ideal, and $\la_{i_1}<k$.
By assumption, $\Phi_1$ has no corner of the form $(\bullet,i_0)$, and hence $\gn(\Phi_1)$ is $P_i$-stable.

By Lemma~\ref{lem:atwist}, we have
\begin{equation}\label{eqn:step1}
H^{\bullet}\bigl(T^*_\Psi X,\cO_{T^*_\Psi X}(\la)\bigr)
 \cong
\mathsf q H^{\bullet}\bigl(T^*_{\Phi_1}X,\cO_{T^*_{\Phi_1}X}(\la+\epsilon_{i_1}-\epsilon_{i_0})\bigr).
\end{equation}

If $\la_{i_1-1}>\la_{i_1}$, we may take $\ell=1$. If $\Phi_1$ has a corner of the form $(i+1,\bullet)$, let $\Phi:=\Phi_1 \setminus (i+1,j)$ for such a corner $(i+1,j)$. Otherwise, set $\Phi:=\Phi_1$. In the former case, Lemma~\ref{lem:ctwist} shows that $(\Phi,\la+\epsilon_{i_1}-\epsilon_{i_0})$ satisfies the required conclusion. In the latter case, the pair $(\Phi_1,\la+\epsilon_{i_1}-\epsilon_{i_0})$ already satisfies the conclusion (as we have $E_{i+1,j} \notin \gn(\Phi_1)$ for every $j$ by assumption).

Otherwise, if $\la_{i_1-1}=\la_{i_1}$ and $\Phi_1$ has no corner of the form $(\bullet,i_1)$, Lemma~\ref{lem:triv} implies
\[
H^{\bullet}\bigl(T^*_{\Phi_1}X,\cO_{T^*_{\Phi_1}X}(\la+\epsilon_{i_1}-\epsilon_{i_0})\bigr)=0,
\]
and hence the desired vanishing holds.
In the remaining case, we continue this process.

We build a maximal chain of corners
\[
(i_1,i_0),\ (i_2,i_1),\ \ldots,\ (i_\ell,i_{\ell-1}),
\qquad
i_\ell<\cdots<i_2<i_1<i,
\]
such that $\la_{i_s-1}=\la_{i_s}$ for all $1\le s<\ell$, together with root ideals $\Phi_s$ obtained by successively adding a corner of shape $(i_s,\bullet)$ to $\Phi_{s-1}$. By repeated application of Lemma~\ref{lem:atwist}, we obtain
\begin{equation}\label{eqn:qell}
H^{\bullet}\bigl(T^*_\Psi X,\cO_{T^*_\Psi X}(\la)\bigr)
 \cong
\mathsf q^{\ell}
H^{\bullet}\bigl(T^*_{\Phi_\ell}X,\cO_{T^*_{\Phi_\ell}X}(\la+\epsilon_{i_\ell}-\epsilon_{i_0})\bigr).
\end{equation}

At this final stage, there are two possibilities.

\emph{(a)} We have $\la_{i_\ell-1}=\la_{i_\ell}$ and there is no corner of the form $(\bullet,i_\ell)$.
At the $s$-th step ($s > 0$), $\Phi_s$ has no corner of the form $(i_s-1,\bullet)$ as we added a corner of the form $(i_s,\bullet)$ to the root ideal in a $(k,i_s-1)$-adapted pair $(\Phi_{s-1}, \la )$; in particular, for $s=\ell$ there is no corner $(i_\ell-1,\bullet)$.
Thus the hypotheses of Lemma~\ref{lem:triv}, namely, the absence of corners of shapes $(\bullet,i_\ell)$ and $(i_\ell-1,\bullet)$ together with $\la_{i_\ell-1}=\la_{i_\ell}$, are satisfied. Therefore
\[
H^{\bullet}\bigl(T^*_{\Phi_\ell}X,\cO_{T^*_{\Phi_\ell}X}(\la+\epsilon_{i_\ell}-\epsilon_{i_0})\bigr)=0.
\]
Hence~\eqref{eqn:qell} gives $H^{\bullet}\bigl(T^*_\Psi X,\cO_{T^*_\Psi X}(\la)\bigr)\equiv 0$.

\emph{(b)} We have $\la_{i_\ell-1}>\la_{i_\ell}$.
Since $\Phi_\ell$ is obtained from $\Psi$ by adding corners of the shape $(i_{s},\bullet)$ to the root ideal in a $(k,i_s-1)$-adapted pair $(\Phi_{s-1}, \la )$ for $1 \le s \le \ell$, we deduce that $\gn(\Phi_\ell)$ has no corners of the forms
\[
(i_{\ell-1}-1,\bullet),\, (i_{\ell-2}-1,\bullet),\, \ldots,\, (i_1-1,\bullet),\, (i_0-1,\bullet)
\quad
(\bullet,i_{\ell-1}),\, (\bullet,i_{\ell-2}),\, \ldots,\ (\bullet,i_1),\, (\bullet,i_0),
\]
where $\ell-1$ terms overlap between the two sequences, and $\la_{i_a-1}=\la_{i_a}$ for $0< a<\ell$ and $\la_{i_0-1}+1=\la_{i_0}$. Hence, repeated application of Lemma~\ref{lem:ctwist} allows us to remove the corners
\[
(i_0,\bullet),\, (i_1,\bullet),\, \ldots,\, (i_{\ell-2},\bullet),\, (i_{\ell-1},\bullet)
\]
successively, skipping $(i_0,\bullet)$ if no such corner exists. This produces a root ideal $\Phi$ that agrees with $\Psi$ in all rows of index $>i+1$.
The resulting pair $(\Phi,\ \la+\epsilon_{i_\ell}-\epsilon_{i_0})$ is $(k,i+1)$-adapted, and the above discussion yields
\[
H^{\bullet}\bigl(T^*_{\Phi_\ell}X,\cO_{T^*_{\Phi_\ell}X}(\la+\epsilon_{i_\ell}-\epsilon_{i_0})\bigr)
 \cong
H^{\bullet}\bigl(T^*_{\Phi}X,\cO_{T^*_{\Phi}X}(\la+\epsilon_{i_\ell}-\epsilon_{i_0})\bigr).
\]
Substituting into~\eqref{eqn:qell} gives
\[
H^{\bullet}\bigl(T^*_\Psi X,\cO_{T^*_{\Psi}X}(\la)\bigr)
 \cong \mathsf q^{\ell} 
H^{\bullet}\bigl(T^*_{\Phi}X,\cO_{T^*_{\Phi}X}(\la+\epsilon_{i_\ell}-\epsilon_{i_0})\bigr).
\]

This completes the proof of~\textup{(2)}, and hence of the proposition.
\end{proof}

The following Corollary is a geometric counterpart of~\cite[Theorem~7.12]{BMPS}.

\begin{cor}\label{cor:k-straight}
Let $k \in \Z_{\ge 1}$.  
Assume $\la \in \sP$ satisfies
\begin{equation}\label{eqn:qpos}
\la_j \le \la_i + 1 \qquad \text{for all} \quad 1 \le i \le j \le n.
\end{equation}
Then the graded $\g[z]$-module
\[
H^{\bullet}\bigl(T^*_{\Psi[\la,k]} X, \cO_{T^*_{\Psi[\la,k]} X}(\la)\bigr)
\]
is either $0$ or, up to a grading shift, isomorphic to
\[
H^\bullet \bigl(T^*_{\Psi[\mu,k]} X, \cO_{T^*_{\Psi[\mu,k]} X}(\mu)\bigr)
\]
for some $\mu \in \sP^+$ with
\[
\{ \la_i > k \} = \{ \mu_i > k \} \quad \text{as multisets}, \qquad \text{and} \quad \mu \ge \la.
\]
\end{cor}

For $\la \in \sP$ satisfying~\eqref{eqn:qpos} but $\la \notin \sP^+$, we call the minimal index $i$ with $\la_{i+1} = \la_i + 1$ the \emph{ascent} of $\la$; this is the index $i$ at which Corollary~\ref{cor:k-straight} applies Theorem~\ref{thm:gmaps} or Proposition~\ref{prop:upshifts} to $(\Psi[\la,k],\la)$.

\begin{proof}[Proof of Corollary~\ref{cor:k-straight}]
If $\la$ is already dominant, we may take $\mu=\la$. Otherwise, let $i$ be the ascent of $\la$. If $i+1<n$, then $\la_{i+2}\le\la_i+1=\la_{i+1}$. Consequently, whenever row $i+1$ of $\Psi[\la,k]$ is nonempty, it has a corner. If $i+1=n$, the row is empty and this condition is vacuous. We distinguish two cases.

If $\la_{i+1}>k$, then $\Psi[\la,k]$ has a corner $(i,i+1)$. Set $\nu:=\la+\epsilon_i-\epsilon_{i+1}$ and $\Psi:=\Psi[\la,k]\setminus(i,i+1)$. The sheaf exact sequence in the proof of Theorem~\ref{thm:gmaps}, together with
\[
H^\bullet\bigl(T^*_\Psi X,\cO_{T^*_\Psi X}(\la)\bigr)=0
\]
from Lemma~\ref{lem:triv}, yields an isomorphism in every cohomological degree:
\[
H^\bullet\bigl(T^*_{\Psi[\la,k]}X,\cO(\la)\bigr)
\cong\mathsf q H^\bullet\bigl(T^*_{\Psi[\la,k]}X,\cO(\nu)\bigr).
\]
Here $\nu$ swaps $\la_i$ and $\la_{i+1}$, and $\Psi[\nu,k]=\Psi[\la,k]$. Thus the swap contributes a grading shift by $1$.

If we have $\la_{i+1} \le k$, then the pair $(\Psi[\la,k],\la)$ is $(k,i)$-adapted and has no corner of shape $(i,\bullet)$ by~\eqref{eqn:qpos}.
If there is no corner of shape $(\bullet,i+1)$, Proposition~\ref{prop:upshifts}\textup{(1)} yields
\[
H^{\bullet}(T^*_{\Psi[\la,k]} X,\cO(\la)) = 0.
\]
Otherwise, Proposition~\ref{prop:upshifts}\textup{(2)} either gives vanishing or, with $j=i_\ell<i$, provides an ``upshift''
\[
\nu=\la+\epsilon_j-\epsilon_{i+1},
\]
together with an isomorphism (up to a grading shift)
\[
H^{\bullet}\bigl(T^*_{\Psi[\la,k]} X,\cO(\la)\bigr)
 \cong
H^{\bullet}\bigl(T^*_{\Psi[\nu,k]} X,\cO(\nu)\bigr).
\]
In either nonvanishing case, the new pair $(\Psi[\nu,k],\nu)$ is $(k,i+1)$-adapted, still satisfies~\eqref{eqn:qpos}, has $\{ \nu_i > k \} = \{ \la_i > k \}$, and moreover $\nu \ge \la$.  
If we still have $\nu_{j+1}=\nu_{j}+1$ for some $j$, then we replace $(\la,i)$ by $(\nu,j)$ and iterate the same procedure.  
At each step, the ascent index strictly increases, so the process terminates after at most $n$ iterations.

If the procedure ends with vanishing, the claim follows. 
Otherwise, it terminates with a dominant $\mu$ satisfying $\{ \mu_i > k \} = \{ \la_i > k \}$, and $\mu \ge \la$.  
Composing the isomorphisms completes the proof.
\end{proof}

The following is the geometric counterpart of~\cite[Theorem~9.2]{BMPS}:

\begin{cor}[Geometric $k$-Pieri rule]\label{cor:kPieri}
Let $k \in \Z_{\ge 1}$ and let $\la \in \sP^+$. Set $\Psi := \Psi[\la,k]$. For $1 \le i \le n$ we have
\[
H^{>0}\bigl(T^*_\Psi X, \cO_{T^*_\Psi X}(\la - \epsilon_i)\bigr) = 0,
\]
and the space
\begin{equation}\label{eqn:sbbranch}
H^0\bigl(T^*_\Psi X, \cO_{T^*_\Psi X}(\la - \epsilon_i)\bigr)
\end{equation}
admits a finite filtration by
\[
H^0\bigl(T^*_{\Psi[\mu,k]} X, \cO_{T^*_{\Psi[\mu,k]} X}(\mu)\bigr)
\quad\text{with}\quad
\mu \in \sP^+,\quad \la - \epsilon_i \le \mu,
\]
as graded $\g[z]$-modules. In addition, we have
\[
\{ \la_j -\delta_{ij} > k \} = \{\mu_j > k \} \qquad \text{(as multisets)}. 
\]
\end{cor}

\begin{proof}
We prove the assertion by descending induction on $i$. For $i=n$, the weight $\la-\epsilon_n$ is dominant and $\Psi[\la-\epsilon_n,k]=\Psi[\la,k]$, so higher cohomology vanishes and the filtration has a single step. Assume the assertion for all indices greater than $i$, and consider $i<n$.

If there is no corner of shape $(i,\bullet)$ in $\Psi$ or $\la_i > k$, then we have $\Psi = \Psi[\la-\epsilon_i,k]$. If $\la_i > \la_{i+1}$ in addition, then $\la - \epsilon_i \in \sP^+$ and there is nothing to prove. If $\la_i = \la_{i+1}$ in addition, then $\la - \epsilon_i$ satisfies~\eqref{eqn:qpos}, and Corollary~\ref{cor:k-straight} applied to $(\Psi[\la-\epsilon_i,k],\la-\epsilon_i)$, together with Theorem~\ref{thm:Cat}\textup{(4)}, yields the assertion.

If $\la_i \le k$ and there is a corner of shape $(i,j)$ in $\Psi$ (so that $j = i+k+1-\la_i$), then we have a short exact sequence
$$0 \to \mathsf{q} \cO_{T^*_\Psi X} ( \la - \epsilon_j) \to \cO_{T^*_\Psi X} ( \la - \epsilon_i) \to \cO_{T^*_{\Psi'} X} ( \la - \epsilon_i) \to 0,$$
where $\Psi' = \Psi \setminus(i,j) = \Psi [\la-\epsilon_i,k]$. We apply Corollary~\ref{cor:k-straight} and Theorem~\ref{thm:Cat}\textup{(4)} to deduce that $H^{\bullet} ( T^*_{\Psi'} X, \cO_{T^*_{\Psi'} X} ( \la - \epsilon_i)  )$ satisfies $H^{>0} = 0$ and $H^0$ is of the desired form. Since $i < j$, the induction hypothesis implies that $H^{\bullet} ( T^*_\Psi X, \cO_{T^*_\Psi X} ( \la - \epsilon_j)  )$ also satisfies $H^{>0} = 0$ and $H^0$ admits the desired filtration. Thus, the induction proceeds.

This proves the assertion.
\end{proof}

\begin{cor}\label{cor:Bsk}
For $\la \in \Par_m(n)$ such that $n > k - \la_1+1$, the $\g[z]$-module $H^0 ( \sX_{\Psi[\la,k]}, \mathcal V (\la) )$ admits a finite filtration whose associated graded is the direct sum of grading shifts of $\{\sk_\mu\}_{\mu \in \Par_{m-1}(n)}$.

We assume $n > m$. For each $\la \in \Par_m(n)$ $($resp. $\la \in \Par_m^{(k)}(n))$, the module
\[
\WS ( \Res^{A_m}_{A_{m-1}} ( \SW (  \sk_\la )))
\]
admits a finite graded $\g[z]$-module filtration whose associated graded is the direct sum of grading shifts of $\{\sk_\mu\}_{\mu \in \Par_{m-1}(n)}$ $($resp. $\{\sk_\mu\}_{\mu \in  \Par_{m-1}^{(k)}(n)})$.
\end{cor}

\begin{proof}
By Corollary~\ref{cor:kPieri}, the hypothesis~\eqref{eqn:hcvan} of Theorem~\ref{thm:assgr} holds for $\Psi := \Psi[\la,k]$. Refining the resulting filtration of $H^0 ( T^*_{\Psi} X, \mathcal V (\la) )$ by Corollary~\ref{cor:kPieri}, we obtain a finite filtration by graded $\g[z]$-submodules whose pieces are $\mathsf q^{s} H^0 ( T^*_{\Psi[\mu,k]} X, \cO_{T^*_{\Psi[\mu,k]} X} ( \mu ) )$ with $\mu \in \sP^+$, $\mu \ge \la - \epsilon_i$, and with the $\g[z]$-structure inherited from $\cO_{\sX_{\Psi[\mu,k]}} ( \mu )$: all maps involved are restrictions of $\widehat{G}[\![z]\!]$-equivariant maps of line bundles on $\sX_\Phi$'s (Theorem~\ref{thm:Psi-rest} and Theorem~\ref{thm:gmaps}).

The same d\'evissage can be carried out on $\sX_\Phi$: the short exact sequences of Theorem~\ref{thm:gmaps} exist on $\sX_\Phi$, Lemma~\ref{lem:triv} is replaced by Lemma~\ref{lem:shift-vanish}\textup{(1)}, and the isomorphisms of Lemma~\ref{lem:ctwist} follow from these two applied to the corner $(\hs_i(\Psi),i)$ (resp.\ $(i+1,j)$), as the weight $\la + \epsilon_{\hs_i(\Psi)} - \epsilon_i$ (resp.\ $\la + \epsilon_{i+1} - \epsilon_j$) satisfies the hypothesis of Lemma~\ref{lem:shift-vanish}\textup{(1)}. Each coordinate of a weight is decreased at most once in the procedure, so every piece satisfies $\mu_n \ge -1$; moreover $\mu \notin \Par(n)$ occurs only for $i \ge 2$, in which case $\mu_1 \ge \la_1 > k+1-n$ and hence $\al_{1n} \in \Psi[\mu,k]$. Therefore $H^{>0} ( \sX_{\Psi[\mu,k]}, \cO_{\sX_{\Psi[\mu,k]}} ( \mu ) ) = 0$ by Theorem~\ref{thm:Cat}\textup{(4)} if $\mu \in \Par(n)$, and $H^{\bullet} = 0$ by Lemma~\ref{lem:shift-vanish}\textup{(2)} otherwise. By d\'evissage, $H^{>0} ( \sX_\Psi, \cO_{\sX_\Psi} ( \la - \epsilon_i ) ) = 0$ for all $i$, so that $F^{\bullet} \mathcal V (\la)$ induces a filtration of $H^0 ( \sX_\Psi, \mathcal V ( \la ) )$ with pieces $H^0 ( \sX_\Psi, \cO_{\sX_\Psi} ( \la - \epsilon_i ) )$, each of which is filtered by $\{ \mathsf q^s \sk_\mu \}_{\mu \in \Par_{m-1}(n)}$. This proves the first assertion.

Comparing with Theorem~\ref{thm:Cat}\textup{(4)} and Theorem~\ref{thm:tri}, we also find that $H^0 ( \sX_\Psi, \mathcal V ( \la ) ) \subset H^0 ( T^*_\Psi X, \mathcal V ( \la ) )$ is the span of all irreducible polynomial $G$-submodules. Hence Proposition~\ref{prop:SW-res} yields the second assertion when $n > k - \la_1 + 1$. If $n \le k - \la_1 + 1$, then $m < n$ gives $|\la| < k$, so $\sk_\la \cong V_\la^*$ by Lemma~\ref{lem:gchcomp}, and the assertion follows from the branching rule of $\Sym_m$ and $\sk_\nu \cong V_\nu^*$ for $\nu \in \Par_{m-1}(n)$. The assertion for $\la \in \Par_m^{(k)}(n)$ follows from the last assertion of Corollary~\ref{cor:kPieri}.
\end{proof}

\section{$k$-Schur branching law and its applications}\label{sec:kbr}

We continue with the setting of the previous section. Fix $k \in \Z_{>0}$ and $n > m$.

For each root ideal $\Phi$, by Theorem~\ref{thm:tvb}\textup{(1)}, we have a $\widehat{G}[\![z]\!]$-equivariant vector bundle $\mathcal V$ on $\sX_\Phi$ equipped with a filtration $F^{\bullet} \mathcal V$ such that $F^{i} \mathcal V/F^{i-1} \mathcal V \cong \cO_{\sX_\Phi} (-\epsilon_i)$ ($1 \le i \le n$). We set
\[
\mathcal F^{d,p}_\Phi ( \gamma ) := \bigl(  \wedge^d F^p \mathcal V \bigr) \otimes_{\cO_{\sX_\Phi}} \cO_{\sX_\Phi} ( \gamma ) \qquad  0 \le d, p \le n,\gamma \in \sP.
\]
We have $\mathcal F^{d,p}_\Phi ( \gamma ) = 0$ for $d > p$ and $\mathcal F^{d,d}_\Phi (\gamma) \cong \cO_{\sX_{\Phi}} ( \gamma - \varpi_d )$ for $0 \le d \le n$ by Theorem~\ref{thm:tvb}\textup{(1)}.

If $p < n$, then we have a short exact sequence of $\widehat{G}[\![z]\!]$-equivariant sheaves on $\sX_\Phi$
\begin{equation}
0 \rightarrow \mathcal F^{d+1,p}_\Phi \rightarrow \mathcal F^{d+1,p+1}_\Phi \rightarrow \mathcal F^{d,p}_\Phi ( - \epsilon_{p+1}) \rightarrow 0,\label{eqn:dpSES}
\end{equation}
whose fiber over $B/B$ gives rise to a short exact sequence of graded $B$-modules by Theorem~\ref{thm:tvb}\textup{(2)}
\[
0 \rightarrow \wedge^{d+1} ( F^{p} V_1 )  \rightarrow \wedge^{d+1} ( F^{p+1} V_1 ) \rightarrow ( \wedge^d ( F^p V_1 ) ) \otimes \C_{\epsilon_{p+1}} \rightarrow 0.
\]
Since $\mathcal F^{d,p}_{\Phi} ( \gamma )$ is obtained by restricting $\mathcal F^{d,p}_{\Delta^+} ( \gamma )$ along $\sX_\Phi \hookrightarrow \sX_{\Delta^+}$, we henceforth omit $\Phi$ from the notation.

The following result is a geometric counterpart of~\cite[Theorem~9.3]{BMPS}.

\begin{thm}\label{thm:Fdp}
Let $k \in \Z_{> 0}$. Consider a family of modules
\begin{equation}
H^0 ( T^*_{\Psi[\mu,k]} X, \cO_{T^*_{\Psi[\mu,k]} X} ( \mu ) ) \qquad (\mu \in \sP^+),
\label{eqn:Fdp}
\end{equation}
where the $\g[z]$-action is inherited from $\cO_{\sX_{\Psi[\mu,k]}}(\mu)$.
For each $0 \le d \le p \le n$ and $\la \in \sP^+$,
\begin{enumerate}
\item The module $H^0 ( T^*_{\Psi[\la,k]} X, \mathcal F^{d,p} ( \la ))$ admits a finite filtration by~\eqref{eqn:Fdp} and
\[
H^{>0} ( T^*_{\Psi[\la,k]} X, \mathcal F^{d,p} ( \la )) = 0.
\]
\item If $p < n$, then the module $H^0 ( T^*_{\Psi[\la,k]} X, \mathcal F^{d,p} ( \la -\epsilon_{p+1} ))$ admits a finite filtration by~\eqref{eqn:Fdp} and
\[
H^{>0} ( T^*_{\Psi[\la,k]} X, \mathcal F^{d,p} ( \la - \epsilon_{p+1})) = 0.
\]
\end{enumerate}
\end{thm}

\begin{proof}
Since $\mathcal F^{d,p}_\Phi (\xi)$ is the restriction of the $\widehat{G}[\![z]\!]$-equivariant bundle $\wedge^d F^p \mathcal V \otimes_{\cO_{\sX_\Phi}} \cO_{\sX_\Phi}(\xi)$ for every root ideal $\Phi$ and $\xi \in \sP$, and every map in the d\'evissage is the restriction of a $\widehat{G}[\![z]\!]$-equivariant map of sheaves on the $\sX_\Phi$'s (cf.\ the proof of Lemma~\ref{lem:ctwist}), the resulting long exact sequences consist of graded $\g[z]$-modules, each term $H^\bullet ( T^*_\Phi X, \mathcal F^{d,p}_\Phi (\xi) )$ carrying the structure inherited from $\sX_\Phi$ up to grading shifts. In particular, the pieces~\eqref{eqn:Fdp} carry the structure inherited from $\cO_{\sX_{\Psi[\mu,k]}}(\mu)$.

We first reduce~\textup{(2)} to~\textup{(1)} for fixed $(d,p)$. The case $d=0$ is given by Corollary~\ref{cor:kPieri}, and we now generalize the discussion there to $d>0$.

We show that in the d\'evissage of $\cO_{T^*_{\Psi[\la,k]} X} ( \la - \epsilon_{p+1} )$, no application of Lemma~\ref{lem:triv}, Lemma~\ref{lem:ctwist}, Lemma~\ref{lem:atwist}, or of the first case of Corollary~\ref{cor:k-straight} has index $p$. The indices used in an application of Proposition~\ref{prop:upshifts} with the sequence $i_1,\ldots,i_\ell$ are $i_0 - 1 = i$ and $i_a - 1$ ($1 \le a \le \ell$), and every ascent is $\ge p+1$. Hence it suffices to show $p + 1 \notin \{i_1,\ldots,i_\ell\}$ for every such application.

If $\la_{p+1} > k$, then $\Psi[\la,k] = \Psi[\la-\epsilon_{p+1},k]$, and the straightening of $\la - \epsilon_{p+1}$ consists only of the swaps in the first case of Corollary~\ref{cor:k-straight}, at the indices $p+1, p+2, \ldots$. Hence we may assume $\la_{p+1} \le k$.

Set $r_0:=p+1$ and define
\[
r_{a+1}:=r_a+k+1-\la_{r_a}
\]
while $r_a\le n$. Let $c$ be determined by $r_c\le n<r_{c+1}$. Since $\la_{r_a}\le\la_{p+1}\le k$, this sequence is strictly increasing. For $a<c$, the pair $(r_a,r_{a+1})$ is a corner of $\Psi[\la,k]$. By the proof of Corollary~\ref{cor:kPieri}, the d\'evissage reduces to straightening $(\Phi_a,\la-\epsilon_{r_a})$ for $0\le a\le c$, where
\[
\Phi_a:=\Psi[\la-\epsilon_{r_a},k]
=\begin{cases}
\Psi[\la,k]\setminus(r_a,r_{a+1})&a<c,\\
\Psi[\la,k]&a=c.
\end{cases}
\]
For a fixed $a$, row $r_a$ of $\Phi_a$ starts at column $r_{a+1} +1$, or is empty if $r_{a+1}>n$.

We claim that throughout this straightening every ascent $i$ satisfies $i+1\le r_{a+1}$, while the current weight $\xi$ satisfies
\[
\xi_{r_e}=\la_{r_e}\quad(e<a),\qquad \xi_{r_a}=\la_{r_a}-1.
\]
If the initial weight is dominant, there is nothing to check. Otherwise its ascent is $r_a$. On a nonvanishing branch, an upshift at ascent $i$ changes only positions $i_\ell$ (increased) and $i+1$ (decreased), and the next ascent, if any, is $i+1$. In particular, entries with index less than $r_a$ never decrease.

Suppose that the stated conditions hold before an upshift with $i_0=i+1\le r_{a+1}$. The rows $r_0,\ldots,r_{a-1},r_a$ of $\Psi[\xi,k]$ start at columns $r_1,\ldots,r_a,r_{a+1}+1$, respectively. If $i_b=r_e$ in its corner chain, then earlier corner additions have changed only rows with the larger indices $i_1,\ldots,i_{b-1}$. Row $r_e$ is therefore unchanged, and its corner $(r_e,i_{b-1})$ forces $i_{b-1}=r_{e+1}$, read as $r_{a+1}+1$ when $e=a$. For $e<a$, we have $r_{e+1}\le r_a<i_0$, so this inverse tracing cannot stop at $i_0$. Repeating it forces $r_{a+1}+1$ into the decreasing chain starting at $i_0\le r_{a+1}$, a contradiction. Thus the chain avoids $r_0,\ldots,r_a$, and the displayed entries of $\xi$ remain unchanged.

Finally, if $r_{a+1}\le n$ and the ascent reaches $r_{a+1}-1$, there is no corner of shape $(r_{a+1}-1, \bullet)$. Row $r_a$ starts at $r_{a+1}+1$. For $j<r_a$, we have $\xi_j\ge\la_j\ge\la_{r_a}$, so row $j$ starts at or before column $r_{a+1}-(r_a-j)<r_{a+1}$. For $r_a<j<r_{a+1}$, condition~\eqref{eqn:qpos} gives $\xi_j\le\xi_{r_a}+1=\la_{r_a}$, so row $j$ starts at or after column $r_{a+1}+(j-r_a)>r_{a+1}$. Hence there is no corner of shape $(\bullet, r_{a+1})$. Thus Proposition~\ref{prop:upshifts}\textup{(1)} gives vanishing at index $r_{a+1}-1\ge r_a >p$. If $r_{a+1} >n$, this boundary cannot be reached. This proves the claim and, in particular, $p+1=r_0\notin\{i_1,\ldots,i_\ell\}$ in every upshift.

Consequently, every application of Lemma~\ref{lem:triv}, Lemma~\ref{lem:ctwist}, Lemma~\ref{lem:atwist}, and of the first case of Corollary~\ref{cor:k-straight} in the d\'evissage has index $i \neq p$. For such $i$, the $B$-module $\wedge^d ( F^p V_1 )$ is $P_i$-stable, so that $\mathcal F^{d,p} |_{T^*_\Phi X}$ is the pullback of a vector bundle along $\eta : T^*_\Phi X \to G \times^{P_i} \gn ( \Phi )$ whenever $\gn(\Phi)$ is $P_i$-stable. By the projection formula, the vanishing in Lemma~\ref{lem:triv} and the isomorphisms in Lemma~\ref{lem:ctwist} and Lemma~\ref{lem:atwist} remain valid after tensoring with $\mathcal F^{d,p}$, while the short exact sequences in Theorem~\ref{thm:gmaps} remain exact after tensoring with the vector bundle $\mathcal F^{d,p}$. Hence the d\'evissage of $\cO ( \la - \epsilon_{p+1} )$ in Corollary~\ref{cor:kPieri} yields a d\'evissage of $\mathcal F^{d,p} ( \la - \epsilon_{p+1} )$ whose terms are $H^\bullet ( T^*_{\Psi[\kappa,k]} X, \mathcal F^{d,p} (\kappa) )$ for $\kappa \in \sP^+$ (up to grading shifts) and vanishing terms. Assuming~\textup{(1)} for $(d,p)$, the resulting d\'evissage yields
\begin{equation}
H^{>0} ( T^*_{\Psi[\la,k]} X, \mathcal F^{d,p} (\la-\epsilon_{p+1}) ) = 0,\label{eqn:Fdpvan}
\end{equation}
and $H^0 ( T^*_{\Psi[\la,k]} X, \mathcal F^{d,p} (\la-\epsilon_{p+1}) )$ admits a filtration by~\eqref{eqn:Fdp}. Thus~\textup{(1)} implies~\textup{(2)} for $(d,p)$.

We now prove the first assertion by induction on $d$. The case $d=0$ follows from Theorem~\ref{thm:Cat}\textup{(4)}. Assume that the assertions~\textup{(1)(2)} for $d-1$ have been established for all values of $p$ in their respective ranges. For this fixed $d$, we prove~\textup{(1)} by induction on $p$, beginning with $p=d$.

For the base case, the isomorphism
\[
\wedge^d(F^d \mathcal V)
\cong
\wedge^{d-1}(F^{d-1}\mathcal V)\otimes_{\cO_{\sX_{\Delta^+}}}\cO_{\sX_{\Delta^+}} (-\epsilon_d)
\]
induces an isomorphism
\[
\mathcal F^{d,d}(\la)
\cong
\mathcal F^{d-1,d-1}(\la-\epsilon_d).
\]
Hence~\textup{(1)} for $(d,d)$ follows from~\textup{(2)}  for $(d-1,d-1)$, which is available by the induction hypothesis on $d$.

Suppose now that~\textup{(1)} holds for $(d,p)$, where $d\le p<n$. Replacing $d$ by $d-1$ in~\eqref{eqn:dpSES} and twisting by $\mathcal O(\la)$, we obtain
\[
0
\longrightarrow
\mathcal F^{d,p}(\la)
\longrightarrow
\mathcal F^{d,p+1}(\la)
\longrightarrow
\mathcal F^{d-1,p}(\la-\epsilon_{p+1})
\longrightarrow
0.
\]
By the induction hypothesis on $p$, the first term has vanishing higher cohomology and filtered zeroth cohomology. By the induction hypothesis on $d$, the same is true of the third term. The associated long exact sequence therefore yields
\[
H^{>0}\bigl(
T^*_{\Psi[\la,k]}X,
\mathcal F^{d,p+1}(\la)
\bigr)=0.
\]
Moreover, its zeroth cohomology is an extension of modules admitting filtrations by~\eqref{eqn:Fdp}; hence it also admits such a filtration. This proves~\textup{(1)} for $(d,p+1)$.

The induction on $p$ proves~\textup{(1)} for the fixed value of $d$. The straightening argument given above then proves~\textup{(2)} for every $d\le p<n$.

This completes the induction on $d$ and proves the theorem.
\end{proof}

\begin{thm}\label{thm:skfilt}
Let $k \in \Z_{> 1}$. For each $\la \in \Par_m(n)$, the module $\mathtt{s}^{(k-1)}_\la$ admits a finite filtration by
\begin{equation}
\sk_\mu = H^0 ( \sX_{\Psi[\mu,k]}, \cO_{\sX_{\Psi[\mu,k]}} ( \mu ) ) \qquad (\mu \in \Par_m(n)).
\label{eqn:piece-general}
\end{equation}
\end{thm}

\begin{proof}
We consider $\la \in \Par_m(n)$, and set $\nu := \la + \varpi_n$. Since $\Psi[\nu,k] = \Psi[\nu-\varpi_n,k{-}1] = \Psi[\la,k{-}1] $, we have a natural isomorphism of graded $\g[z]$-modules
\begin{equation}
H^{0} \bigl(T^*_{\Psi[\nu,k]} X, \wedge^n \mathcal V \otimes_{\cO_{T^*_{\Psi[\nu,k]} X}} \cO_{T^*_{\Psi[\nu,k]} X}( \nu )\bigr) \cong H^{0} \bigl(T^*_{\Psi[\nu,k]} X, \cO_{T^*_{\Psi[\nu,k]} X}( \la )\bigr)\label{eqn:ktok-1}.
\end{equation}
By Theorem~\ref{thm:Fdp}\text{(1)} for $(d,p)=(n,n)$, this admits a finite filtration by~\eqref{eqn:Fdp}.

The left-hand side of
\[
\mathtt{s}^{(k-1)}_\la = H^0 ( \sX_{\Psi[\la,k-1]}, \cO_{\sX_{\Psi[\la,k-1]}}(\la)) \subset H^0 ( T^*_{\Psi[\la,k-1]} X, \cO_{T^*_{\Psi[\la,k-1]} X}(\la))
\]
is a graded $\g[z]$-submodule spanned by all irreducible polynomial $G$-modules in the right-hand side by Theorem~\ref{thm:Cat}\textup{(4)}. Applying Proposition~\ref{prop:pol}, we deduce that $\mathtt{s}^{(k-1)}_\la$ admits a finite filtration by~\eqref{eqn:piece-general}.
\end{proof}

\begin{rem}
Theorem~\ref{thm:skfilt} is a counterpart of~\cite[\S9.2]{BMPS}. Note that we cannot prove Theorem~\ref{thm:skfilt} as a direct generalization of Corollary~\ref{cor:kPieri}.

The following exhibits a remnant~\eqref{eqn:Hone} of the obstruction discussed in~\cite[Remark~9.5]{BMPS} in the case $n=6$, $k=3$, and $\Psi = \Psi[222222,3]$.

In Corollary~\ref{cor:kPieri}, we appealed to certain sheaf-theoretic calculations, which eliminate the term
\[
H^{\bullet} ( T^*_{\Psi} X, \cO_{T^*_{\Psi} X} (122222)),
\]
since it vanishes as
\[
H^{\bullet} ( T^*_{\Psi} X, \cO_{T^*_{\Psi} X} (221222)) 
\equiv 
H^{\bullet} ( T^*_{\Psi [232322,3]} X, \cO_{T^*_{\Psi [232322,3]} X} (221222)) \equiv 0
\]
by Lemma~\ref{lem:ctwist} (applied to the 1st and 2nd rows to add a corner $(2,3)$, and to the 5th and 6th rows to add a corner $(4,5)$) and Lemma~\ref{lem:triv} (applied to the 3rd and 4th rows), and
\[
H^{\bullet} ( T^*_{\Psi [122222,3]} X, \cO_{T^*_{\Psi [122222,3]} X} (122222)) \equiv 0
\]
by Lemma~\ref{lem:triv} (applied to the 1st and 2nd rows).  
However, we still have
\begin{equation}
H^{1} ( T^*_{\Psi} X, \cO_{T^*_{\Psi} X} (122212)) \neq 0.\label{eqn:Hone}
\end{equation}
This nonvanishing prevents a direct termwise extension of the argument in Corollary~\ref{cor:kPieri} to the line-bundle filtration of $\wedge^2\mathcal V(222222)$.
\end{rem}

\begin{cor}\label{cor:sskbranch}
Let $k \in \Z_{> 1}$. For each $\la \in \Par_m^{(k-1)} ( n )$, the module $\mathtt{s}^{(k-1)}_\la$ admits an $\sk$-filtration.
\end{cor}

\begin{proof}
By Theorem~\ref{thm:skfilt}, $\mathtt{s}^{(k-1)}_\la$ admits a filtration by $\{\sk_\mu\}_{\mu \in \Par_m (n)}$. We have
\begin{equation}
\gch \mathtt{s}^{(k-1)}_\la = \sum_{\mu \in \Par_{m}^{(k)}(n)} c_\mu(q) \gch \mathtt{s}^{(k)}_\mu \qquad (c_\mu(q) \in \Z_{\ge 0}[q])\label{eqn:kbnumpos}
\end{equation}
by~\cite[(1.1)]{BMPS}.

By Theorem~\ref{thm:tri},~\eqref{eqn:kbnumpos} is the only possible way to write
\[
\gch \mathtt{s}^{(k-1)}_\la = \sum_{\mu \in \Par_{m}(n)} c' _\mu(q) \gch \mathtt{s}^{(k)}_\mu \qquad (c'_\mu(q) \in \Z_{\ge 0}[q]),
\]
and hence $\mathtt{s}^{(k-1)}_\la$ admits a filtration by $\{\sk_\mu\}_{\mu \in \Par_m^{(k)} (n)}$ as required.
\end{proof}

\begin{cor}\label{cor:HLfilt}
Let $k \in \Z_{> 0}$. For each $\la \in \Par_m^{(k)}(n)$, the module $\HL_\la$ admits an $\sk$-filtration.
\end{cor}

\begin{proof}
For each $\la \in \Par_m(n)$, we have $\HL_\la = \mathtt{s}^{(1)}_\la$.  
Thus, by inductively applying Theorem~\ref{thm:skfilt} to $\mathtt{s}^{(r)}_\mu$ ($\mu \in \Par_m(n)$) from $r=1$ to $r=k-1$, we obtain a filtration by~\eqref{eqn:piece-general}.

In view of Theorem~\ref{thm:tri}, it follows that the characters $\gch \sk_\la$ ($\la \in \Par_m(n)$) are linearly independent over $\Z[q,q^{-1}]$.  
Hence we obtain a unique expression
\[
\gch \HL_\la = \sum_{\mu \in \Par_m(n)} c_\mu(q)\, \gch \sk_\mu
\qquad (c_\mu(q) \in \Z_{\ge 0}[q]).
\]
Here $c_\mu(q)$ records the graded occurrence of $\sk_\mu$ in $\HL_\la$; in particular, $c_\mu(q)=0$ means that $\sk_\mu$ does not appear in $\HL_\la$.  

By~\cite[Theorem~4.1]{BMPS2}, we have $c_\mu(q)=0$ for all $\mu \in \Par_m(n) \setminus \Par_m^{(k)}(n)$.  
This establishes the claim.
\end{proof}

\begin{cor}\label{cor:Kfilt}
Let $k \in \Z_{> 0}$. For each $\la \in \Par_m^{(k)}$, the module $K_\la$ admits a $\ssc^{(k)}$-filtration.
\end{cor}

\begin{proof}
Apply $\SW$ to Corollary~\ref{cor:HLfilt} by choosing $n = m+1$.
\end{proof}

\begin{cor}\label{cor:Pfilt}
Let $k \in \Z_{> 0}$. For each $\la \in \Par_m^{(k)}$, the module $P_\la$ admits a $\ssc^{(k)}$-filtration.
\end{cor}

\begin{proof}
By Theorem~\ref{thm:Kostka}\textup{(3)(4)}, $P_\la$ has a decreasing separable filtration by grading shifts of $K_\mu$ with $\mu\le\la$. Since $\mu_1\le\la_1\le k$, each such $K_\mu$ has a finite $\ssc^{(k)}$-filtration by Corollary~\ref{cor:Kfilt}. Refine by these finite filtrations. The refinement is finite in each internal degree, so it remains separable and gives the required filtration.
\end{proof}

\begin{cor}\label{cor:maxco}
Let $k \in \Z_{> 0}$. Let $\la_{\max} \in \Par_m^{(k)}(n)$ denote the unique maximal partition with respect to the dominance order, namely, $(k^{\lfloor m/k\rfloor}(m-k\lfloor m/k\rfloor))$. Then
\[
\HL_{\la_{\max}} \cong \sk_{\la_{\max}}.
\]
\end{cor}

\begin{proof}
We recall the unique expansion
\[
\gch \HL_{\la_{\max}} = \sum_{\mu \in \Par_m^{(k)}(n)} c_\mu (q)\, \gch \sk_\mu 
\qquad (c_{\mu}(q) \in \Z_{\ge 0}[q]),
\]
as established in the proof of Corollary~\ref{cor:HLfilt}.  
By~\cite[\S3(16)]{LM07} we have
\[
\ch \HL_{\la_{\max}} = \ch \sk_{\la_{\max}}.
\]
It follows that $c_{\la_{\max}}(q) = 1$ and $c_\mu(q)=0$ for $\mu \neq \la_{\max}$, compare Theorem~\ref{thm:ksch}.

The natural map
\[
\HL_{\la_{\max}}\longrightarrow\sk_{\la_{\max}}
\]
is surjective by Corollary~\ref{cor:Cat-sg}. Since the graded characters are equal, this map is an isomorphism.
\end{proof}

\begin{cor}\label{cor:korth}
Let $k \in \Z_{> 0}$. For each $\la \in \Par_m^{(k)}$, we have
\[
\ext_{A_m}^{\bullet}\bigl(\ssc_\la^{(k)},\,L_{\mu}\bigr)\equiv 0
\qquad \bigl(\mu \in \Par_m \setminus \Par_m^{(k)}\bigr).
\]
In particular, the minimal projective resolution of $\ssc_\la^{(k)}$ in $A_m\gmod$ involves only projectives from $\{P_\mu\}_{\mu \in \Par_m^{(k)}}$ (up to grading shifts).
\end{cor}

\begin{proof}
We argue by downward induction on $\la$ with respect to dominance.
Assume the statement holds for all $\mu \in \Par_m^{(k)}$ with $\mu>\la$, and we prove it for $\la$.

For the base case, let $\la_{\max}\in \Par_m^{(k)}$ be the maximal element.
Applying $\SW$ to Corollary~\ref{cor:maxco} by choosing $n > m$ gives $K_{\la_{\max}} \cong \ssc^{(k)}_{\la_{\max}}$, and hence
\[
\ext_{A_m}^{\bullet} \bigl(\ssc^{(k)}_{\la_{\max}},\,L_{\mu}\bigr)\equiv 0
\qquad \bigl(\mu \in \Par_m \setminus \Par_m^{(k)}\bigr)
\]
by Corollary~\ref{cor:Kostka} and Lemma~\ref{lem:interval}. This constitutes the base step.

For general $\la$, Corollary~\ref{cor:Kfilt} and Theorem~\ref{thm:CHid} imply that
$\ker \bigl(K_\la \twoheadrightarrow \ssc_\la^{(k)}\bigr)$ admits an $\ssc^{(k)}$-filtration.
By~\eqref{eqn:Ktri} and Theorem~\ref{thm:ksch}, all subquotients are of the form $\ssc_{\mu}^{(k)}$ with $\mu>\la$.
Applying d\'evissage with respect to this filtration using the induction hypothesis, we obtain
\[
\ext_{A_m}^{\bullet} \bigl(\ker(K_\la \to \ssc_\la^{(k)}),\,L_{\gamma}\bigr)\equiv 0
\qquad \bigl(\gamma \in \Par_m \setminus \Par_m^{(k)}\bigr).
\]

Together with Corollary~\ref{cor:Kostka} for $K_\la$ and Lemma~\ref{lem:interval}, the long exact sequence attached to
\[
0 \longrightarrow \ker(K_\la \to \ssc_\la^{(k)}) \longrightarrow K_\la \longrightarrow \ssc_\la^{(k)} \longrightarrow 0
\]
yields
\[
\ext_{A_m}^{\bullet} \bigl(\ssc_{\la}^{(k)},\,L_{\mu}\bigr)\equiv 0
\qquad \bigl(\mu \in \Par_m \setminus \Par_m^{(k)}\bigr).
\]
This proves the vanishing.

The vanishing implies that in a minimal projective resolution of $\ssc_\la^{(k)}$ only projective covers of simples $L_\mu$ with $\mu\in \Par_m^{(k)}$ can occur, up to grading shifts.
This completes the proof.
\end{proof}

\begin{cor}
Let $k \in \Z_{> 0}$. Assume $n>m$. For each $\la \in \Par_m^{(k)}(n)$ and each $\mu \in \Par_m(n) \setminus \Par_m^{(k)}(n)$, we have
\[
\ext_{\g[z]\gmod_m}^{\bullet} \bigl(\sk_\la,\,V_{\mu}^*\bigr)\equiv 0.
\]
\end{cor}

\begin{proof}
Apply $\WS$ to Corollary~\ref{cor:korth}.
\end{proof}

\section{Filtration criteria}\label{sec:main}

We retain the notation from the previous section. We fix
$k,m \in \mathbb{Z}_{>0}$.

\begin{lem}\label{lem:fres}
Let $M$ be a finite-dimensional graded $A_m$-module such that
\[
  \ext^{\bullet}_{A_m}\bigl(M,L_{\mu}\bigr)=0 \qquad\text{for every }\mu \in \Par_m \setminus \Par^{(k)}_m.
\]
Then $M$ admits a finite resolution
\[
0 \to Q_r \to Q_{r-1} \to \cdots \to Q_1 \to Q_0 \to M \to 0
\]
whose terms $Q_0,\ldots,Q_r$ are finite-dimensional $\ssc^{(k)}$-filtered modules.
\end{lem}

\begin{proof}
By Theorem~\ref{thm:gdimA}, there is a finite projective resolution $(P^{\bullet},d_{\bullet})$ of $M$ by graded maps of internal degree $0$.
The hypothesis $\ext^{\bullet}_{A_m} (M,L_\mu)=0$ for $\mu\notin\Par^{(k)}_m$ allows us to replace
$(P^{\bullet},d_{\bullet})$ by a minimal projective resolution in which only graded shifts of $\{P_{\la}\}_{\la\in \Par^{(k)}_m}$ appear among the direct summands of $\bigoplus_a P^a$. Indeed, in a minimal graded projective resolution, the multiplicity of a grading shift $\mathsf q^sP_\mu$ in $P^a$ is detected by the corresponding
graded piece of $\hom_{A_m}(P^a,L_\mu)\cong \ext_{A_m}^a(M,L_\mu)$.

Let $s_{-1}$ be the maximal internal degree with $M_{s_{-1}}\neq 0$. Since each $d_i$ has internal degree $0$, the complex $(P^{\bullet},d_{\bullet})$ is exact in all internal degrees $> s_{-1}$. 

Enumerate
\[
  \Par^{(k)}_m =\{\la^{(1)},\la^{(2)},\ldots,\la^{(l)}\}
\]
so that $\la^{(i)}\le \la^{(j)}$ implies $i\le j$. For each $a \in \Z_{\ge 0}$ we have a finite filtration
\[
  0=F_0 P^a \subset F_1 P^a \subset \cdots \subset F_{l}P^a=P^a
\]
such that $F_j P^a /F_{j-1}P^a$ is a direct sum of grading shifts of $\tK_{\la^{(j)}}$ by Corollary~\ref{cor:filtproj2}.

We now truncate in internal degree and construct successive quotients. Let $s_0$ be the maximal degree in which (the head of) $\tK_{\la^{(j)}}$ ($1 \le j \le l$) appears in $\bigoplus_a P^a$, which is finite since $\dim M < \infty$ and $A_m$ has a finite global dimension. Set
\[
s_i  :=1+ \max\{ s_{-1}, s_0 \} +\sum_{j=1}^{i-1} \Bigl(1 + \max\{e\mid (K_{\la^{(j)}})_e\neq 0\}\Bigr)
  \qquad (1\le i\le l).
\]
For $0\le i \le l$, define $(P^{[i],\bullet},d_{\bullet})$ as the termwise quotient of $(P^{\bullet},d_{\bullet})$
by the $A_m$-submodules generated by all simple graded $\Sym_m$-submodules of the form
$\mathsf q^{s}L_{\la^{(j)}}$ with $1\le j\le i$ and $s>s_j$.

Because the differentials are $A_m$-linear, $\Sym_m$-equivariant, and
homogeneous of internal degree $0$, the $A_m$-submodules generated by the
$\Sym_m$-direct summands of the form
$\mathsf q^sL_{\la^{(i+1)}}$ with $s>s_{i+1}$ are stable under the
differential. Hence we obtain a quotient of complexes
\[
  (P^{[i],\bullet},d_{\bullet}) \longrightarrow
  (P^{[i+1],\bullet},d_{\bullet}).
\]

\begin{claim}\label{clm:pik-filt}
For each $0\le i<l$, each term of the subcomplex
\[
\ker ( P^{[i],\bullet} \to P^{[i+1],\bullet} ) \subset P^{[i],\bullet}
\]
admits a decreasing separable filtration whose associated graded is a direct
sum of modules $\mathsf q^{s}K_{\la^{(i+1)}}$ with $s>s_{i+1}$.
Moreover, for each $0\le i\le l$ and each $a\ge0$, the term $P^{[i],a}$
admits a decreasing separable filtration by
$\{K_{\la}\}_{\la \in \Par^{(k)}_m}$.
\end{claim}

\begin{proof}
We prove the assertion by induction on $i$ from $i=0$. More precisely, we prove that
each term $P^{[i],a}$ admits a filtration whose successive quotients are of
the form
\begin{equation}\label{eqn:Rmult}
  R_\la \otimes_{\cend_{A_m} (\tK_{\la})}\tK_{\la}
  \qquad (\la\in\Par^{(k)}_m)
\end{equation}
for some finitely generated graded $\cend_{A_m}(\tK_\la)$-module $R_\la$,
where $R_\la$ is free when $\la=\la^{(b)}$ with $b>i$, and
\[
  (R_{\la^{(b)}})_e=0
  \qquad\text{for }e>s_b \text{ and } b\le i.
\]
The case $i=0$ follows from Corollary~\ref{cor:filtproj2} and
Lemma~\ref{lem:Ksplit}.

Assume that the assertion has been proved for $i$. The definition of
$s_{i+1}$ implies that, for every $b\le i$,
\[
  \bigl(R_{\la^{(b)}} \otimes_{\cend_{A_m}(\tK_{\la^{(b)}})}
  \tK_{\la^{(b)}}\bigr)_e=0
  \qquad\text{for }e>s_{i+1}.
\]
Thus all contributions in internal degree $>s_{i+1}$ with label
$\la^{(b)}$ for $b\le i$ have already been eliminated in
$P^{[i],a}$.

Together with Lemma~\ref{lem:Ksplit} and the triangularity
\eqref{eqn:Ktri}, this shows that every map from
$\tK_{\la^{(i+1)}}$ to $P^{[i],a}$ of degree $>s_{i+1}$ lands in the
$\la^{(i+1)}$-part
\[
R_{\la^{(i+1)}}\otimes_{\cend_{A_m}(\tK_{\la^{(i+1)}})}
\tK_{\la^{(i+1)}}.
\]

In addition, triangularity~\eqref{eqn:Ktri} gives
\[
[\tK_{\la^{(b)}}:L_{\la^{(c)}}]=0
\qquad
(b>i+1\geq c).
\]
Hence, by Lemma~\ref{lem:Ksplit}, the passage from
$P^{[i],a}$ to $P^{[i+1],a}$ does not alter the subquotients
of the form~\eqref{eqn:Rmult} indexed by
$\la^{(b)}$ with $b>i+1$.

Therefore, the passage from
$P^{[i],a}$ to $P^{[i+1],a}$ only replaces the coefficient module
$R_{\la^{(i+1)}}$ with its truncation in degrees 
$\le s_{i+1}$. In particular,
$P^{[i+1],a}$ and $\ker ( P^{[i],a} \to P^{[i+1],a} )$ ($a \in \Z_{\ge 0}$) admit decreasing separable filtrations by $\{K_{\la}\}_{\la \in \Par^{(k)}_m}$ and by $K_{\la^{(i+1)}}$, respectively. In addition, only $s > s_{i+1}$ is allowed for $\mathsf{q}^s K_{\la^{(i+1)}}$ to appear in $\ker ( P^{[i],a} \to P^{[i+1],a} )$ by the construction of truncations. This advances the induction and proves the claim.
\end{proof}

\begin{claim}\label{clm:pik-filt-diff}
For each $0\leq i<l$, the filtration of
\[
\ker\bigl(P^{[i],\bullet}\to P^{[i+1],\bullet}\bigr)
\]
given by Claim~\ref{clm:pik-filt} can be chosen to be compatible with
the differential.
\end{claim}

\begin{proof}
By~\eqref{eqn:Rmult}, each term of the kernel complex is of the form
\[
R_{\la^{(i+1)}}\otimes_{\cend_{A_m}(\tK_{\la^{(i+1)}})}
\tK_{\la^{(i+1)}}.
\]
The differential induces homogeneous maps between the coefficient
modules $R_{\la^{(i+1)}}$. By Theorem~\ref{thm:Kostka}\textup{(4)(5)}, the
$K_{\la^{(i+1)}}$-filtration of $\tK_{\la^{(i+1)}}$ is induced by
the grading of $\cend_{A_m} ( \tK_{\la^{(i+1)}})$. Hence the induced
filtrations in~\eqref{eqn:Rmult} are functorial for homogeneous maps
of coefficient modules. Therefore, they may be chosen compatibly with
the differential.
\end{proof}

Returning to the proof of the lemma, we prove that $P^{[i],\bullet}$ is acyclic in internal degree $>s_0$ by induction on $i$ from $i=0$. The case $i=0$ is clear by the definition of the projective resolution.
Let
\[
\ker_i^\bullet:=\ker(P^{[i],\bullet}\to P^{[i+1],\bullet}).
\]
For each $s>s_{i+1}$, let $U_s^\bullet$ denote the multiplicity
complex of $\mathsf q^sL_{\la^{(i+1)}}$ in
$P^{[i],\bullet}$. Since $P^{[i],\bullet}$ is exact in internal
degrees $> s_0$, and since taking an
$L_{\la^{(i+1)}}$-isotypic component is exact over
$\Sym_m$, each $U_s^\bullet$ is acyclic.

By Claim~\ref{clm:pik-filt} and Claim~\ref{clm:pik-filt-diff}, the complex $\ker_i^\bullet$ admits a decreasing separable filtration by subcomplexes
whose associated graded terms are direct sums of
$\mathsf q^sK_{\la^{(i+1)}}$, with multiplicity spaces
$U_s^\bullet$. Thus, the acyclicity of the complexes $U_s^\bullet$ implies that the associated graded complex is acyclic. Since the filtration
is finite in each internal degree and separable,
$\ker_i^\bullet$ is acyclic. Hence the quotient complex $P^{[i+1],\bullet}$ has
the same homology as $P^{[i],\bullet}$. This proves the induction.

At the final stage $i=l$, the coefficient modules corresponding to all
$\la^{(1)}, \ldots, \la^{(l)}$ have been truncated above the bounds
$s_1, \ldots, s_l$. Hence each term $P^{[l],a}$ is supported in only finitely
many internal degrees. Since each graded piece is finite-dimensional, each
$P^{[l],a}$ is finite-dimensional. By Claim~\ref{clm:pik-filt}, the terms
$P^{[l],a}$ are filtered by $\{K_\la\}_{\la \in \Par^{(k)}_m}$. Therefore, by Corollary~\ref{cor:Kfilt}, each
$P^{[l],a}$ is a finite-dimensional $\ssc^{(k)}$-filtered module. Setting
$Q_a:=P^{[l],a}$ for $0\le a\le r$ yields the desired finite resolution of
$M$.
\end{proof}

\begin{thm}\label{thm:kfilt}
Suppose that we have a family of graded $A_m$-modules $\{ T_\la \}_{\la \in S}$ indexed by $S \subset \Par_m$ with the following properties:
\begin{itemize}
\item For each $\la \in \Par^{(k)}_m$, the module $\ssc^{(k)}_\la$ admits a finite filtration by $\{ T_\mu \}_{\mu \in S}$.
\item For an inclusion $f : N \hookrightarrow M$ of two graded $A_m$-modules with finite filtrations by $\{ T_\la \}_{\la \in S}$, $\coker f$ also admits a finite filtration by $\{ T_\la \}_{\la \in S}$.
\end{itemize}
Let $M$ be a finite-dimensional graded $A_m$-module such that
\[
\ext^{\bullet}_{A_m} ( M, L_\la ) = 0 \qquad (\la \in \Par_m \setminus \Par_m^{(k)}).
\]
Then $M$ admits a finite filtration by $\{ T_\la \}_{\la \in S}$.
\end{thm}

\begin{proof}
By Lemma~\ref{lem:fres}, $M$ admits a finite
graded $A_m$-module resolution
\[
0 \to Q_r \to Q_{r-1} \to \cdots \to Q_1 \to Q_0 \to M \to 0
\]
in which each $Q_i$ admits a finite $\ssc^{(k)}$-filtration.

The first assumption implies that each $Q_i$ also admits a finite filtration by $\{ T_\la \}_{\la \in S}$.

Applying the second assumption to the inclusion
\[
Q_r \hookrightarrow Q_{r-1},
\]
we may replace the last two terms of the resolution with their cokernel, which admits a finite filtration by $\{ T_\la \}_{\la \in S}$. Thus we obtain a shorter exact sequence of the same
form for $M$. Iterating this process $r$ times, we
deduce that $M$ itself admits a finite filtration by $\{ T_\la \}_{\la \in S}$.
\end{proof}

\section{A combinatorial criterion for $\ssc^{(k)}$-filtrations}\label{sec:comb}

We retain the setting of the previous section and fix $k,m \in \mathbb{Z}_{>0}$. The formulation of Hypothesis~\ref{hypo:cycle} uses only the material from Subsections~\ref{subsec:part},~\ref{subsec:roots},~\ref{subsec:Amod}, and Section~\ref{sec:ksch}, together with standard facts on $k$-Schur functions; see, for example,~\cite{LLM03,BMPS}.

The goal of this section is to give a combinatorial criterion which implies the following filtration property, used in the proof of Theorem~\ref{fthm:cmain} and Corollary~\ref{fcor:rmp}.

\begin{conjecture}\label{conj:ker}
Let $M,N \in A_m\gmod$ be finite-dimensional $\ssc^{(k)}$-filtered modules, and let
\[
f : N \hookrightarrow M
\]
be an injective homomorphism of graded $A_m$-modules. Then $\coker f$ also admits a $\ssc^{(k)}$-filtration.
\end{conjecture}

For each pair $\la,\mu \in \Par^{(k)}_m$, we set
\[
e_k(\la,\mu)
:=
\min \left\{
d \mid
\Hom_{\Sym_m}
\left(
L_{(\mu^{\omega_k})'},
\bigl(P_\la\bigr)_d
\right)
\neq 0
\right\}
-d_k(\mu).
\]
Since $\gch P_\la$ is Macdonald's big Schur function~\cite[Chap.~III~(4.5) and~(4.9)]{Mac95}, we have
\[
\gch P_\la=\ws_\la\!\left[\frac{X}{1-q}\right]=\sum_{\nu \in \Par_m}c_k(\la,\nu;q)\,\ws_\nu
\qquad
\bigl(c_k(\la,\nu;q) \in \Z[\![q]\!]\bigr),
\]
and hence $e_k(\la,\mu)+d_k(\mu)=\operatorname{ord}_q c_k(\la,(\mu^{\omega_k})';q)$. In particular, $e_k$ is computable in \texttt{SageMath}~\cite{sage}. Theorem~\ref{thm:CHid} gives $e_k(\la,\la)=0$.

\begin{hypo}\label{hypo:cycle}
For each sequence $(\la^{(1)},\ldots,\la^{(l)})$ of distinct elements of $\Par^{(k)}_m$ with $l \ge 2$, we have
\begin{equation}
\sum_{i=1}^{l} e_k\bigl(\la^{(i)},\la^{(i+1)}\bigr)>0,
\label{eqn:loop}
\end{equation}
where we set $\la^{(l+1)}:=\la^{(1)}$ by convention.
\end{hypo}

\begin{rem}\label{rem:min}
The bound $m\leq 2k$ in Theorem~\ref{thm:app-cycle-m2k} is not optimal. Since extending this range requires substantially more combinatorial analysis, we record only the simpler proof needed here. The smallest counterexamples to Hypothesis~\ref{hypo:cycle} occur for $m=20$, with $k=4,5$.
\end{rem}

\begin{cor}\label{cor:pot}
Assume Hypothesis~\ref{hypo:cycle}. Then there exists a function $F_k : \Par^{(k)}_m \rightarrow \mathbb R$ such that
\begin{equation}
e_k(\la,\mu)>F_k(\la)-F_k(\mu)
\qquad
\bigl(\la,\mu \in \Par^{(k)}_m,\ \la \neq \mu\bigr).
\label{eqn:Fpot}
\end{equation}
\end{cor}

\begin{proof}
We equip the complete directed graph without loops on $\Par^{(k)}_m$ with the arc weights $e_k(\la,\mu)-\epsilon$, where $\epsilon>0$ is small enough that every directed cycle keeps positive total weight. Fixing $\la_0 \in \Par^{(k)}_m$ and letting $F_k(\la)$ be the minimal total weight of a directed path from $\la$ to $\la_0$, we obtain $F_k(\la) \le e_k(\la,\mu)-\epsilon+F_k(\mu)$ by prefixing the arc from $\la$ to $\mu$; see~\cite[Theorem~7.7]{KV}.
\end{proof}

For $(\la,d),(\mu,d') \in \Par^{(k)}_m \times \Z$, we write
\[
(\mu,d') \preceq (\la,d)
\qquad \Longleftrightarrow \qquad
d' \le d+F_k(\la)-F_k(\mu).
\]
In the sequel, we freely use the fact that $\ssc^{(k)}_\la$ has head $L_\la$ in degree $0$ and simple socle $\mathsf q^{d_k(\la)}L_{(\la^{\omega_k})'}$, and that it is the unique finite-dimensional such module (Theorem~\ref{thm:ksch} and Theorem~\ref{thm:CHid}).

\begin{conjecture}\label{conj:ntext}
There exists a function $F_k : \Par^{(k)}_m \rightarrow \mathbb R$ with the following property.
Let $(\la,d),(\mu,d') \in \Par^{(k)}_m \times \Z$.
\begin{enumerate}
\item If $(\mu,d') \preceq (\la,d)$, then every short exact sequence
\[
0 \rightarrow \mathsf q^{d'}\ssc^{(k)}_\mu \rightarrow E \rightarrow \mathsf q^d\ssc^{(k)}_\la \rightarrow 0
\]
of graded $A_m$-modules splits;
\item If $(\la,d) \preceq (\mu,d')$, then every nonzero graded $A_m$-module map
$\mathsf q^{d'}\ssc^{(k)}_\mu \rightarrow \mathsf q^d\ssc^{(k)}_\la$ is an isomorphism,
and $(\la,d)=(\mu,d')$ in that case.
\end{enumerate}
Equivalently, if $(\mu,d') \preceq (\la,d)$, then
$\Ext^1_{A_m\gmod}(\mathsf q^{d}\ssc^{(k)}_\la,\mathsf q^{d'}\ssc^{(k)}_\mu)=0$, and
\[
\Hom_{A_m\gmod}(\mathsf q^{d}\ssc^{(k)}_\la,\mathsf q^{d'}\ssc^{(k)}_\mu) \cong \C^{\delta_{\la,\mu}\delta_{d,d'}}.
\]
\end{conjecture}

\begin{lem}\label{lem:ntext}
Hypothesis~\ref{hypo:cycle} implies Conjecture~\ref{conj:ntext}, with the function $F_k$ specified in Corollary~\ref{cor:pot}.
\end{lem}

\begin{proof}
Assume $(\mu,d') \preceq (\la,d)$, and let $E$ be an extension as in~\textup{(1)} with morphisms $\psi$ and $\phi$. Let $\widetilde\phi : \mathsf q^dP_\la \rightarrow E$ be a lift of the projective cover of $\mathsf q^d\ssc^{(k)}_\la$ through $\phi$, and set $N:=\mathrm{Im}\,\widetilde\phi$. Since $\phi(N)=\mathsf q^d\ssc^{(k)}_\la$ and $\ker\phi=\psi(\mathsf q^{d'}\ssc^{(k)}_\mu)$, we have $N+\psi(\mathsf q^{d'}\ssc^{(k)}_\mu)=E$. If $N \cap \psi(\mathsf q^{d'}\ssc^{(k)}_\mu)$ does not contain $\psi(\mathsf q^{d'}\soc\ssc^{(k)}_\mu)$, then this intersection is a submodule of $\mathsf q^{d'}\ssc^{(k)}_\mu$ missing its simple socle, and hence is zero. In this case the sequence splits, as required.

It remains to rule out the case $\psi(\mathsf q^{d'+d_k(\mu)}L_{(\mu^{\omega_k})'}) \subset N$. Since $N$ is a quotient of $\mathsf q^dP_\la$, this forces $d' \ge d+e_k(\la,\mu)$. For $\la \neq \mu$, this contradicts~\eqref{eqn:Fpot} and $(\mu,d') \preceq (\la,d)$. For $\la=\mu$, it gives $d' \ge d$, and hence $d'=d$ by $(\la,d') \preceq (\la,d)$. Then $N$ contains $\psi(\mathsf q^{d+d_k(\la)}L_{(\la^{\omega_k})'}) \subset \ker\phi$, while $\phi(N)=\mathsf q^d\ssc^{(k)}_\la$ contains its socle $\mathsf q^{d+d_k(\la)}L_{(\la^{\omega_k})'}$. Thus $L_{(\la^{\omega_k})'}$ occurs at least twice in degree $d+d_k(\la)$ in the quotient $N$ of $\mathsf q^dP_\la$, which contradicts the multiplicity one in Theorem~\ref{thm:CHid}. This proves~\textup{(1)}.

We next prove~\textup{(2)}. Assume $(\la,d) \preceq (\mu,d')$ and let
$\theta : \mathsf q^{d'}\ssc^{(k)}_\mu \rightarrow \mathsf q^d\ssc^{(k)}_\la$ be nonzero. Then
$\mathrm{Im}\,\theta$ contains the socle $\mathsf q^{d+d_k(\la)}L_{(\la^{\omega_k})'}$, so that
$L_{(\la^{\omega_k})'}$ occurs in $\mathsf q^{d'}P_\mu$ in degree $d+d_k(\la)$. This yields
$d-d' \ge e_k(\mu,\la)$. If $\mu \neq \la$, then~\eqref{eqn:Fpot} gives
$d'<d+F_k(\la)-F_k(\mu)$, which is incompatible with $(\la,d) \preceq (\mu,d')$.
Hence $\mu=\la$, and $e_k(\la,\la)=0$ combined with $(\la,d) \preceq (\la,d')$ gives $d=d'$. Thus $\theta$ is an isomorphism by Theorem~\ref{thm:CHid}.
\end{proof}

For a finite-dimensional $\ssc^{(k)}$-filtered module $M$, we denote by $G_M$ the multiset of labels $(\la,d)$ of the filtration factors $\mathsf q^d\ssc^{(k)}_\la$, and we call $(\la,d) \in G_M$ \emph{minimal} if $(\la,d) \preceq (\mu,d')$ for each $(\mu,d') \in G_M$.

\begin{lem}\label{lem:label}
Let $M$ be a finite-dimensional $\ssc^{(k)}$-filtered module and let $\pi : M \twoheadrightarrow \mathsf q^dL_\la$ be a surjection with $\la \in \Par^{(k)}_m$. Then we have $(\la,d) \in G_M$.
\end{lem}

\begin{proof}
Fix a filtration $0=M_0 \subset \cdots \subset M_r=M$ with $M_i/M_{i-1} \cong \mathsf q^{d_i}\ssc^{(k)}_{\la^{(i)}}$, and let $j$ be minimal with $M_j \not\subset \ker\pi$. Then $\pi$ induces a nonzero map $M_j/M_{j-1} \rightarrow \mathsf q^dL_\la$. Since $\hd \ssc^{(k)}_{\la^{(j)}}=L_{\la^{(j)}}$ sits in degree $0$, this implies $(\la^{(j)},d_j)=(\la,d)$.
\end{proof}

\begin{lem}\label{lem:peel}
Assume Conjecture~\ref{conj:ntext}. Let $M$ be a finite-dimensional $\ssc^{(k)}$-filtered module and let $(\la,d) \in G_M$ be minimal. Then a surjection $M \twoheadrightarrow \mathsf q^d\ssc^{(k)}_\la$ exists, and the kernel of every such surjection $p$ is $\ssc^{(k)}$-filtered with $G_{\ker p}=G_M \setminus \{(\la,d)\}$.
\end{lem}

\begin{proof}
We first construct such a surjection. Fix a filtration with factors $\mathsf q^{d_i}\ssc^{(k)}_{\la^{(i)}}$ $(1 \le i \le r)$, where $r=|G_M|$, and let $(\la^{(i)},d_i)=(\la,d)$. If $i<r$, then $(\la^{(i)},d_i) \preceq (\la^{(i+1)},d_{i+1})$, and hence
\[
0 \rightarrow \mathsf q^{d_i}\ssc^{(k)}_{\la^{(i)}} \rightarrow M_{i+1}/M_{i-1} \rightarrow \mathsf q^{d_{i+1}}\ssc^{(k)}_{\la^{(i+1)}} \rightarrow 0
\]
splits by Conjecture~\ref{conj:ntext}\textup{(1)}. Replacing $M_i$ by the preimage in $M_{i+1}$ of the second summand moves $(\la,d)$ to the position $i+1$ without changing the labels. Iterating this, we may assume $i=r$, which yields the required surjection.

We prove the second assertion by induction on $r$, the case $r=1$ being trivial. Let $p : M \twoheadrightarrow \mathsf q^d\ssc^{(k)}_\la$ be a surjection, let $M_1 \cong \mathsf q^{a}\ssc^{(k)}_\gamma$ be the smallest nonzero term of a filtration of $M$, and set $\theta:=p|_{M_1}$. Since $(\la,d) \preceq (\gamma,a)$, Conjecture~\ref{conj:ntext}\textup{(2)} leaves two cases.

If $\theta$ is an isomorphism, then $M=M_1 \oplus \ker p$, and $\ker p \cong M/M_1$ is $\ssc^{(k)}$-filtered with labels $G_M \setminus \{(\gamma,a)\}=G_M \setminus \{(\la,d)\}$.

If $\theta=0$, then $p$ factors through a surjection $\bar p : M/M_1 \twoheadrightarrow \mathsf q^d\ssc^{(k)}_\la$, and Lemma~\ref{lem:label} gives $(\la,d) \in G_{M/M_1}$, where it is again minimal. The induction hypothesis shows that $\ker\bar p$ is $\ssc^{(k)}$-filtered with labels $G_M \setminus \{(\gamma,a),(\la,d)\}$, and $\ker p$ is an extension of $\ker\bar p$ by $M_1$.
\end{proof}

\begin{thm}\label{thm:ker}
For fixed $(k,m)$, Conjecture~\ref{conj:ntext} implies Conjecture~\ref{conj:ker}. In particular, Hypothesis~\ref{hypo:cycle} implies Conjecture~\ref{conj:ker}.
\end{thm}

\begin{proof}
Applying $\vee$ and $\circledast$ and using Corollary~\ref{cor:kdual}, it suffices to prove that $\ker f$ is $\ssc^{(k)}$-filtered for a surjection $f : M \twoheadrightarrow N$ between finite-dimensional $\ssc^{(k)}$-filtered modules. We argue by induction on the sum of the $\ssc^{(k)}$-filtration lengths of $M$ and $N$, the case $N=0$ being trivial.

Choose a minimal $(\la,d) \in G_M$. By Lemma~\ref{lem:peel}, we have a surjection $g : M \twoheadrightarrow \mathsf q^d\ssc^{(k)}_\la$ such that $M':=\ker g$ is $\ssc^{(k)}$-filtered with $G_{M'}=G_M \setminus \{(\la,d)\}$.

Assume first $f(M')=N$. Then $M=\ker f+M'$, and hence
\[
\ker f/\ker(f|_{M'}) \cong M/M' \cong \mathsf q^d\ssc^{(k)}_\la .
\]
Here $\ker (f|_{M'})$ is $\ssc^{(k)}$-filtered by the induction hypothesis applied to $f|_{M'} : M' \twoheadrightarrow N$. Therefore $\ker f$ is $\ssc^{(k)}$-filtered.

Assume next $f(M') \subsetneq N$. Then $N/f(M')$ is a nonzero quotient of $M/M' \cong \mathsf q^d\ssc^{(k)}_\la$, so that $\mathsf q^dL_\la$ occurs in $\hd N$, and Lemma~\ref{lem:label} gives $(\la,d) \in G_N$.

We claim that $(\la,d)$ is minimal in $G_N$. Indeed, let $(\mu,d') \in G_N$ be minimal. By Lemma~\ref{lem:peel}, we have a surjection $N \twoheadrightarrow \mathsf q^{d'}\ssc^{(k)}_\mu$, and composing it with $f$ and with the projection to the head yields a surjection $M \twoheadrightarrow \mathsf q^{d'}L_\mu$. Hence $(\mu,d') \in G_M$ by Lemma~\ref{lem:label}, and the minimality of $(\la,d)$ in $G_M$ gives $(\la,d) \preceq (\mu,d')$.

By Lemma~\ref{lem:peel}, we thus obtain a surjection $p_N : N \twoheadrightarrow \mathsf q^d\ssc^{(k)}_\la$ such that $\ker p_N$ is $\ssc^{(k)}$-filtered, and the same lemma applied to $p_N \circ f$ shows that $\ker (p_N \circ f)$ is $\ssc^{(k)}$-filtered. Since $f$ restricts to a surjection
\[
\ker (p_N \circ f) \twoheadrightarrow \ker p_N
\]
with kernel $\ker f$, and both sides have smaller filtration length, the induction hypothesis shows that $\ker f$ is $\ssc^{(k)}$-filtered.

This proves Conjecture~\ref{conj:ker}. The final assertion follows from Lemma~\ref{lem:ntext}.
\end{proof}

\section{Consequences}\label{sec:CH}

In this section, we work in the setting of Subsection~\ref{subsec:Amod} and use several results from Sections~\ref{sec:ksch}, \ref{sec:kbr}, and \ref{sec:main}.

Define $\mathcal M^{(k)}_m$ to be the full subcategory of $A_m\gmod$ consisting of finite-dimensional modules $M$ such that
\begin{equation}\label{eqn:kbddinA}
\ext^{\bullet}_{A_m}(M,L_{\la}) = 0 \quad \text{for all } \la \in \Par_m \setminus \Par^{(k)}_m.
\end{equation}

\begin{thm}\label{thm:crange}
Conjecture~\ref{conj:ker} holds for $m \le \max \{ 2 k, 19 \}$.
\end{thm}

\begin{proof}
Hypothesis~\ref{hypo:cycle} is proved for $m\leq 2k$ in Theorem~\ref{thm:app-cycle-m2k}. We have also verified it by computer for $m\leq 19$; see Remark~\ref{rem:min}. Hence Conjecture~\ref{conj:ker} holds in this range by Theorem~\ref{thm:ker}.
\end{proof}

The next theorem reduces~\cite[Conjecture~5.1.2]{Che10} to Conjecture~\ref{conj:ker}.

\begin{thm}\label{thm:str}
Let $k,m \in \Z_{>0}$. Then
\begin{enumerate}
\item Every finite-dimensional $\ssc^{(k)}$-filtered graded $A_m$-module $M$ belongs to $\mathcal M^{(k)}_m$.
\item Assume Conjecture~\ref{conj:ker} for $(k,m)$. Then every $M \in \mathcal M^{(k)}_m$ admits a finite filtration whose associated graded is a direct sum of grading shifts of $\{\ssk_\la\}_{\la \in \Par^{(k)}_m}$. Consequently, we have
\begin{equation}
\gch M \in \sum_{\la \in \Par_m^{(k)}} \Z_{\ge 0} [q^{\pm 1}] \cdot \ws^{(k)}_\la.\label{eqn:gchpos}
\end{equation}
\end{enumerate}
\end{thm}

\begin{proof}
The assertion~\textup{(1)} for $M = \ssk_\la$ ($\la \in \Par^{(k)}_m$) follows from Corollary~\ref{cor:korth}. Since the condition~\eqref{eqn:kbddinA} is stable under extensions, the general case follows.

In order to prove the first assertion in~\textup{(2)}, it suffices to verify the assumptions of Theorem~\ref{thm:kfilt} by choosing $S = \Par_m^{(k)}$ and $T_\la = \ssc^{(k)}_\la$ for each $\la \in S$. The first assumption of Theorem~\ref{thm:kfilt} is immediate. The second assumption of Theorem~\ref{thm:kfilt} is exactly the assumed Conjecture~\ref{conj:ker}.

Taking the graded character proves the second assertion in~\textup{(2)}.
\end{proof}

\begin{ex}
For $k=2$ and $m=3$, we have
$\mathsf{s}_{1^3}^{(1)}=K_{1^3}\in \mathcal M^{(2)}_3$. Moreover, there is
a short exact sequence
\[
0 \longrightarrow \mathsf{q}^2 \ssc_{21}^{(2)} \longrightarrow \ssc_{1^3}^{(1)} \longrightarrow \ssc_{1^3}^{(2)} \longrightarrow 0.
\]
\end{ex}

\begin{rem}
As noted in~\cite[\S5.1]{Che10}, the subcategory $\mathcal M^{(k)}_m$ is, in general, not abelian. For instance, when $k=2$ and $m=3$,
\[
\gch \ssc_{21}^{(2)} = \ws_{21} + q\,\ws_{3},
\qquad
\gch \ssc_{1^3}^{(2)} = \ws_{1^3} + q\,\ws_{21}.
\]
It follows that $\hd \ssc_{21}^{(2)} = L_{21}$ and $\soc \ssc_{1^3}^{(2)} = \mathsf{q} L_{21}$, and hence
\[
\hom_{A_3}\bigl(\ssc_{21}^{(2)},\ssc_{1^3}^{(2)}\bigr)\neq 0.
\]
This is incompatible with $\mathcal M^{(2)}_3$ being an abelian category in
which the objects $\ssc_{21}^{(2)}$ and $\ssc_{1^3}^{(2)}$ are
nonisomorphic simple objects: in an abelian category, any nonzero morphism
between simple objects is an isomorphism.
\end{rem}

For each $k\in\Z_{>0}$ and $\la\in\Par_m$, we set
\[
\mathsf{W}_\la^{(k)}
:=
\SW\bigl((W_\la^{(k)})^*\bigr)
\in A_m\gmod.
\]
This is well-defined by Proposition~\ref{prop:D-dual}, applied with
$n=m+1>m$.

\begin{thm}\label{thm:Dstr}
Let $k,m \in \Z_{>0}$. Every $M \in \mathcal M^{(k)}_m$ admits a finite filtration whose associated graded is a direct sum of grading shifts of $\{\mathsf{W}_\la^{(k)}\}_{\la \in \Par_m}$. In particular,
\begin{equation}
\gch M \in \sum_{\la \in \Par_m} \Z_{\ge 0} [q^{\pm 1}] \cdot \gch \mathsf{W}_\la^{(k)}. \label{eqn:dgchpos}
\end{equation}
\end{thm}

\begin{proof}
In order to prove the first assertion, it suffices to verify the assumptions of Theorem~\ref{thm:kfilt} by choosing $S = \Par_m$ and $T_\la = \mathsf{W}^{(k)}_\la$ for each $\la \in S$. Both assertions follow after applying $\SW$ (Theorem~\ref{thm:FKM}): the first from Theorem~\ref{thm:D-filt}, and the second from Proposition~\ref{prop:D-ker}.

Taking the graded character proves the second assertion.
\end{proof}

\begin{rem}
Theorem~\ref{thm:Dstr} is relatively weak when $m$ is small compared with $k$, since $\mathsf{W}^{(k)}_\la \cong L_\la$ for every $\la \in \Par_m^{(k)}$.

By contrast, when $m > k^2$, we have $(\la^{\omega_k})' \notin \Par_m^{(k)}$ for every $\la \in \Par_m^{(k)}$.
It follows that every $\ssc_\la^{(k)}$ has a filtration containing a factor
$\mathsf{W}^{(k)}_\mu$ with $\mathsf{W}^{(k)}_\mu \neq L_\mu$.
Thus the affine Demazure filtrations in Theorem~\ref{thm:Dstr} are genuinely nontrivial when $m$ is large relative to $k$.
\end{rem}

Recall that
\[
A_m = \C \Sym_m \ltimes \C[X] \subset B_m = \C \Sym_m \ltimes \C[X,Y].
\]
We denote by $\psi_X: A_m \hookrightarrow B_m$ the embedding induced by the inclusion $\C[X]\hookrightarrow \C[X,Y]$, and by $\psi_Y: A_m \hookrightarrow B_m$ the embedding sending $X_i\mapsto Y_i$ for $1\le i\le m$. We regard $B_m$ as a bigraded algebra by setting
\[
\deg w=(0,0)\ (w\in \Sym_m),\qquad \deg X_i=(1,0),\qquad \deg Y_i=(0,1)\quad (1\le i\le m).
\]
A bigraded $B_m$-module is a decomposition $M=\bigoplus_{p,q\in\Z} M_{p,q}$ such that $\Sym_m$ preserves each $M_{p,q}$, the action of $X_i$ sends $M_{p,q}$ to $M_{p+1,q}$, and the action of $Y_i$ sends $M_{p,q}$ to $M_{p,q+1}$ for all $1\le i\le m$.

\begin{thm}[Garsia--Haiman~\cite{GH93}, Haiman~\cite{Hai01,Hai03}]\label{thm:GH}
For each $\la \in \Par_{m}$, there exists a finite-dimensional bigraded $B_m$-module $\mathsf{GH}_{\la}$ with the following properties:
\begin{enumerate}
\item $\mathsf{GH}_\la \cong \C \Sym_m$ as an (ungraded) $\Sym_m$-module.
\item $\mathsf{GH}_\la$ has a unique simple $B_m$-quotient isomorphic to $\mathsf{triv}$, and a unique simple $B_m$-submodule isomorphic to $\mathsf{sgn}$, up to bidegree shifts.
\item $\ext_{A_m}^{\bullet}\bigl(\psi_X^* \mathsf{GH}_\la, L_{\mu}\bigr) = 0  \qquad \text{for all} \quad \mu \not\ge \la$.
\item $\ext_{A_m}^{\bullet}\bigl(\psi_Y^* \mathsf{GH}_\la, L_{\mu'}\bigr) = \ext_{A_m}^{\bullet}\bigl( \psi_Y^* \mathsf{GH}_\la^{\circledast}, L_{\mu} \bigr) =0 \qquad \text{for all} \quad \mu \not\le \la$.
\end{enumerate}
\end{thm}

\begin{cor}[Macdonald positivity {\cite{Hai01}}]
We have
\[
\widetilde K_{\la,\mu}(q,t)=\sum_{i,j\ge0}q^it^j\,
\dim\Hom_{\Sym_m}\bigl(L_\la,(\mathsf{GH}_\mu)_{i,j}\bigr),
\qquad
\bgch\mathsf{GH}_\mu=\widetilde H_\mu(q,t),
\]
where the convention follows~\cite[Chap.~VI,\S8]{Mac95}.
\end{cor}

We understand that the two gradings of a bigraded $B_m$-module are swapped by $\psi_Y^*$. In particular, we have
\[
 \bgch\bigl(\psi_Y^*\mathsf{GH}_\la^{\circledast}\bigr)
 =\omega\widetilde H_\la(t,q)
 =\sum_{\mu\in\Par_m}\widetilde K_{\mu,\la}(t,q)\ws_{\mu'},
 \qquad \omega(\ws_\mu)=\ws_{\mu'}.
\]

By Theorem~\ref{thm:Dstr}, we have the following:

\begin{cor}[affine Demazure positivity of Macdonald polynomials]\label{cor:dMac}
Let $k\in\Z_{>0}$. For each $\la\in\Par^{(k)}_m$, the module $\psi_Y^*\mathsf{GH}_\la^{\circledast}$ admits a finite filtration by $\{\mathsf{W}^{(k)}_\mu\}_{\mu \in \Par_m}$. In particular,
\begin{equation}\label{eqn:drefpos}
\sum_{\mu\in\Par_m} \widetilde K_{\mu,\la}(t,q)\cdot \ws_{\mu'} \in \sum_{\mu\in\Par_m} \Z_{\ge0}[t,q^{\pm 1}]\cdot \gch \mathsf{W}^{(k)}_\mu.
\end{equation}
\end{cor}

\begin{proof}
Let $\la\in\Par^{(k)}_m$. For $\nu\in\Par_m$ with $\nu\le \la$, we have $\nu\in\Par^{(k)}_m$. In particular,
\[
\Par_m\setminus\Par^{(k)}_m \subset \{\mu\in\Par_m\mid \mu\not\le \la\}.
\]
By Theorem~\ref{thm:GH}\textup{(4)}, it follows that
\[
\ext_{A_m}^{\bullet}\bigl(\psi_Y^*\mathsf{GH}_\la^{\circledast},L_{\mu}\bigr)=0
\qquad \text{for all} \quad \mu\in \Par_m\setminus\Par^{(k)}_m.
\]
Hence $\psi_Y^*\mathsf{GH}_\la^{\circledast}$ and each of its $X$-homogeneous components belong to $ \mathcal M_m^{(k)}$. Theorem~\ref{thm:Dstr} therefore implies that the graded character of $\psi_Y^*\mathsf{GH}_\la^{\circledast}$ has the form~\eqref{eqn:drefpos}.
\end{proof}

\begin{cor}[$k$-Schur positivity of Macdonald polynomials]\label{cor:rMac}
Assume Conjecture~\ref{conj:ker} for $(k,m)$. For each $\la\in\Par^{(k)}_m$, the module $\psi_Y^*\mathsf{GH}_\la^{\circledast}$ admits an $\ssc^{(k)}$-filtration. In particular,
\begin{equation}\label{eqn:refpos}
\sum_{\mu\in\Par_m} \widetilde{K}_{\mu,\la}(t,q)\cdot \ws_{\mu'} \in \sum_{\mu\in\Par^{(k)}_m} \Z_{\ge0}[q,t]\cdot \ws^{(k)}_{\mu}.
\end{equation}
\end{cor}

\begin{proof}
As in the proof of Corollary~\ref{cor:dMac},
Theorem~\ref{thm:GH}\textup{(4)} implies that
\[
\psi_Y^*\mathsf{GH}_\la^{\circledast}\in\mathcal M_m^{(k)}.
\]
Applying Theorem~\ref{thm:str}\textup{(2)} gives an
$\ssc^{(k)}$-filtration compatible with the $X$-grading. Taking bigraded
characters yields~\eqref{eqn:refpos} with coefficients a priori in
$\Z_{\ge0}[q^{\pm1},t]$. Since both
$\psi_Y^*\mathsf{GH}_\la^{\circledast}$ and $\ssc_\mu^{(k)}$ are concentrated in
nonnegative degrees, and the degree-zero part of
$\ssc_\mu^{(k)}$ is $L_\mu$, no negative grading shift can occur.
Hence all coefficients belong to $\Z_{\ge0}[q,t]$.
\end{proof}

\begin{cor}\label{cor:nMac}
Assume that $(k,m)$ satisfies $m \le \max \{2k, 19\}$. For each $\la\in\Par^{(k)}_m$, 
\[
\sum_{\mu\in\Par_m} \widetilde{K}_{\mu,\la}(t,q)\cdot \ws_{\mu'} \in \sum_{\mu\in\Par^{(k)}_m} \Z_{\ge0}[q,t]\cdot \ws^{(k)}_{\mu}.
\]
\end{cor}

\begin{proof}
Combine Theorem~\ref{thm:crange} and Corollary~\ref{cor:rMac}.
\end{proof}

\begin{ex}
Writing $\gch_{q,t}$ for the bigraded Frobenius characteristic, we have
\[
{\textstyle\gch_{q,t} ( \psi_Y^* \mathsf{GH}_{1^3}^{\circledast} )
= \gch\ssc^{(1)}_{1^3}
= \gch\ssc^{(2)}_{1^3} + q^{2}\,\gch\ssc^{(2)}_{21},
\qquad
\gch_{q,t} ( \psi_Y^* \mathsf{GH}_{21}^{\circledast} )
= \gch \ssc^{(2)}_{1^3} + t\,\gch\ssc^{(2)}_{21}.}
\]
\end{ex}

\begin{rem}\label{rem:ksch}
\textbf{(1)} Corollary~\ref{cor:rMac} is the refined Macdonald positivity presented in~\cite{Che10}. However, the $k$-Schur functions used here and those in the original conjecture~\cite[Conjecture~8]{LLM03} are not proved to be identical. Despite recent progress in understanding the various definitions of $k$-Schur functions~\cite{BMPS,BMPS2}, it is still unclear whether Corollary~\ref{cor:rMac} implies the original conjecture. Nevertheless, the notion of $k$-Schur functions used here matches the convention in~\cite{LM03a,LM03b,Che10,LLMS} that has been standard since~\cite{Lam06}.
\textbf{(2)} In view of Haglund--Haiman--Loehr~\cite{HHL05}, Corollary~\ref{cor:rMac} would follow from the claim that a certain class of LLT polynomials expands positively in certain families of $k$-Schur polynomials. Such positivity results were established for $k=2$ and $k=3$ in~\cite{Lee20,Mil19}.
\end{rem}

\appendix
\section{A proof of Hypothesis~\ref{hypo:cycle} in the range $m\leq2k$}
\label{app:cycle-m2k}

In this appendix, we fix $k,m\in\Z_{>0}$ with $m\leq2k$ and prove
Hypothesis~\ref{hypo:cycle}. We retain the notation for partitions,
$k$-conjugation, $d_k$, and the statistics $\mathsf{m},\mathsf{n}$
from Subsection~\ref{subsec:part}, the modules $L_\la,P_\la$ from
Subsection~\ref{subsec:Amod}, and the definition of $e_k$ from
Section~\ref{sec:comb}. The proof uses only this material and
standard facts about representations of $\Sym_m$.
We first prove a cycle inequality for a transportation statistic
and then compare this statistic with $e_k$.

Throughout this appendix, $\la,\mu\in\Par_m^{(k)}$, unless otherwise
specified. Auxiliary partitions need not be $k$-bounded.
By a nonconstant cycle, we mean a finite
cyclic list with at least two distinct vertices.

For $\la\in\Par_m^{(k)}$, set
\begin{equation}\label{eqn:app-Ik}
 I_k(\la)=\{i:i+k-\la_i+1\le\ell(\la)\}.
\end{equation}
Thus, $I_k(\la)$ consists of the rows in which $\Psi[\la,k]$ has
a root in one of the first $\ell(\la)$ columns; see~\eqref{eqn:okdef}
and the definition of $\Psi[\la,k]$ in Section~\ref{sec:ksch}.

\begin{prop}\label{prop:app-low-mass}
The set $I_k(\la)$ is an initial interval with
$|I_k(\la)|<k-\la_1+2$.
If $|I_k(\la)|\ge2$, then $\la_{k-\la_1+2}=1$.
The $(k+1)$-core associated to $\la$ by Theorem~\ref{thm:kcore}
has row lengths
\begin{equation}\label{eqn:app-core-rows}
 \rho_i=\la_i+\la_{i+k-\la_i+1}.
\end{equation}
Here $\rho_j=\la_j$ for $j=i+k-\la_i+1$, so that~\eqref{eqn:app-core-rows}
also reads $\rho_i=\la_i+\rho_{i+k-\la_i+1}$.
\end{prop}
\begin{proof}
Since $i+k-\la_i+1$ is strictly increasing in $i$,
the set $I_k(\la)$ is an initial interval.
There cannot be two consecutive steps
$i\mapsto j=i+k-\la_i+1\mapsto j+k-\la_j+1$
with all three indices at most $\ell(\la)$.
Indeed, counting boxes in the rows up to the last index
(the rows $1,\ldots,i$ contain at least $\la_i$ boxes each, the
$k+1-\la_i$ rows $i+1,\ldots,j$ at least $\la_j$ boxes each, and the
$k+1-\la_j$ rows $j+1,\ldots,j+k-\la_j+1$ at least one box each) would give
\[
 m\ge i\la_i+(k+1-\la_i)\la_j+k+1-\la_j
   =2k+1+(i-1)\la_i+(k-\la_i)(\la_j-1)\ge2k+1.
\]
Put $p=k-\la_1+2$.
The preceding observation gives $|I_k(\la)|<p$.
If $|I_k(\la)|\ge2$ and $\la_p\ge2$,
counting boxes in the rows up to $k-\la_2+3$
(which are nonempty since $2\in I_k(\la)$, using the rows $1,2$,
the $p-2$ rows $3,\ldots,p$, and the remaining $k+3-\la_2-p$ rows) gives
\[
 m\ge\la_1+\la_2+(p-2)\la_p+k+3-\la_2-p
   =2k+1+(k-\la_1)(\la_p-2)\ge2k+1.
\]
Thus $\la_p=1$ whenever $|I_k(\la)|\ge2$.

The sequence $\rho$ in~\eqref{eqn:app-core-rows} is a partition.
For $j=i+k-\la_i+1$, the preceding observation gives
$\rho_j=\la_j$ and $\rho_i=\la_i+\rho_j$.
In row $i$ of $\rho$, the boxes in columns $\rho_j+1,\ldots,\rho_i$
have hook length at most $k$, and all preceding boxes have hook
length at least $k+2$. Theorem~\ref{thm:kcore} therefore identifies
$\rho$ with the core associated to $\la$.
%
\end{proof}

\subsection{A transportation statistic}

We recall the matrix statistic of Chen's thesis~
\cite[Definition~3.1.1]{Che10}.
For $\la,\nu\in\Par_m$, a transportation matrix with row sums $\nu$
and column sums $\la'$ is a finitely supported matrix
$Y=(y_{ij})_{i,j\ge1}$ with entries in $\Z_{\ge0}$ such that
\[
 \sum_j y_{ij}=\nu_i,\qquad \sum_i y_{ij}=\la'_j.
\]
Here, as usual, partitions are extended by zero parts.
Following~\cite[Definitions~3.5.1 and~3.6.1]{Che10}, we define the degree of $Y$ and its minimum by
\begin{equation}\label{eqn:app-distance}
 d(Y)=\sum_{i,j}\binom{y_{ij}}2,\qquad
 d(\la,\nu)=\min_Y d(Y),
\end{equation}
where the minimum is taken over all matrices with these row and
column sums. Such matrices exist since $|\la|=|\nu|$.
For $\la,\mu\in\Par_m^{(k)}$, put
\begin{equation}\label{eqn:app-b}
 b_k(\la,\mu)=d(\la,(\mu^{\omega_k})')-d_k(\mu).
\end{equation}

\begin{lem}\label{lem:app-weighted-bound}
For each pair $\la,\mu\in\Par_m^{(k)}$, we have
\begin{equation}\label{eqn:app-weighted-bound}
 (k+1)b_k(\la,\mu)\ge
 \sum_i\bigl(\mu_i+\mu_{i+k-\mu_i+1}\bigr)(\mu_i-\la_i)
       +\mathsf{n}(\la)-\mathsf{n}(\mu).
\end{equation}
\end{lem}
\begin{proof}
Let $\rho$ be the core associated to $\mu$, also extended by zero parts.
Consider the steps $i\mapsto i+k-\mu_i+1$ for all positive indices.
They form $k+1$ chains: their initial indices are precisely the
elements of $[1,\ell(\mu)+k+1]$ other than the $\ell(\mu)$ distinct
images of $1,\ldots,\ell(\mu)$.
By~\eqref{eqn:defkinv} and Proposition~\ref{prop:app-low-mass}, the values $\rho_s>0$ at these
initial indices $s$ are the parts of $(\mu^{\omega_k})'$.
We may therefore index the rows of a matrix occurring in
$d(\la,(\mu^{\omega_k})')$ by these indices $s$, so that
\[
 \sum_j y_{sj}=\rho_s,\qquad \sum_s y_{sj}=\la'_j.
\]

For such a matrix $Y$, consider the auxiliary weighted sum
\begin{equation}
 (k+1)d(Y)+\sum_{s,j}(s-\rho_s)y_{sj}.\label{eqn:auxsum}
\end{equation}
Its second term is $\sum_s(s-\rho_s)\rho_s$, independent of $Y$.
The contribution of the entry $y_{sj}$ in~\eqref{eqn:auxsum} is the sum of the first
$y_{sj}$ terms of the sequence
\begin{equation}
 s-\rho_s,\quad s-\rho_s+(k+1),\quad
 s-\rho_s+2(k+1),\quad\ldots.\label{eqn:srseq}
\end{equation}
We have
\[
i - \rho_i + k + 1 = i - (\mu_i+\rho_{i+k-\mu_i+1}) + k + 1 = ( i + k -\mu_i + 1 ) - \rho_{i+k-\mu_i+1}
\]
by Proposition~\ref{prop:app-low-mass}. Since $i-\rho_i$ is strictly increasing in $i$,
the union of the sequences~\eqref{eqn:srseq}, arranged in increasing order, is
$1-\rho_1,2-\rho_2,\ldots$.
Each column of $Y$ selects $\la'_j$ terms from this sequence, so its
contribution is at least the sum of the first $\la'_j$ terms.
Summing over columns gives
\begin{equation}\label{eqn:app-weighted-table}
 (k+1)d(Y)+\sum_s(s-\rho_s)\rho_s
 \ge\sum_i(i-\rho_i)\la_i.
\end{equation}

When $\la=\mu$, equality is attained by the following matrix.
Start with the zero matrix. For each $i$ with $\mu_i\ge j$,
add one to the entry in column $j$
and the row indexed by the initial index of the chain containing $i$.
The column sums are $\mu'_j$, and the row sums are $\rho_s$, since
$\rho_s$ is the sum of the parts $\mu_i$ along that chain.
In each column this selects the first $\mu'_j$ terms of the combined sequence.
The term $\binom{y_{sj}}2$ counts pairs of indices in the same
chain and column. Since the parts weakly decrease along each
chain, each $\mu_i$ is counted once for every earlier index in
that chain. By the recurrence for $\rho$, the total is $|\rho|-|\mu|$.
Thus $d(Y)=|\rho|-|\mu|=d_k(\mu)$.
Subtracting the resulting equality from~\eqref{eqn:app-weighted-table}
and minimizing $d(Y)$ over the prescribed row and column sums gives
\[
 (k+1)b_k(\la,\mu)\ge\sum_i(i-\rho_i)(\la_i-\mu_i).
\]
Since $|\la|=|\mu|$, the right-hand side equals
$\sum_i\rho_i(\mu_i-\la_i)+\mathsf{n}(\la)-\mathsf{n}(\mu)$.
Now~\eqref{eqn:app-core-rows} proves the assertion.
\end{proof}

\subsection{A positive lower bound along cycles}

Set
\begin{equation}\label{eqn:app-potential}
 \Phi(\la)=\mathsf{m}(\la)
       +\sum_{i\ge k-\la_1+2}\la_i-\ell(\la).
\end{equation}

\begin{prop}\label{prop:app-edge}
For each pair $\la,\mu\in\Par_m^{(k)}$, we have
\begin{equation}\label{eqn:app-edge-one}
 \sum_i\bigl(\mu_i+\mu_{i+k-\mu_i+1}\bigr)(\mu_i-\la_i)
       +\Phi(\la)-\Phi(\mu)\ge0.
\end{equation}
Equality cannot hold on every edge of a nonconstant cycle.
\end{prop}

\begin{proof}
Put $\Delta_i=\mu_i-\la_i$, so that $\sum_i\Delta_i=0$.
Expanding the definition of $\Phi$, we see that the left-hand side of
\eqref{eqn:app-edge-one} is
\begin{equation}\label{eqn:app-expanded-edge}
 \frac12\sum_i\Delta_i^2+\sum_i\mu_{i+k-\mu_i+1}\Delta_i
 +\sum_{i\ge k-\la_1+2}\la_i
 -\sum_{i\ge k-\mu_1+2}\mu_i-\ell(\la)+\ell(\mu).
\end{equation}
We have an elementary inequality
\begin{equation}\label{eqn:app-tail}
 \sum_{i\ge k-\la_1+2}\la_i \ge \ell(\la)+\la_1-k-1
\end{equation}
obtained by summing the inequalities $\la_i  \ge 1$ over $k - \la_1 + 2 \le i \le \ell ( \la )$ if $\ell ( \la ) > k+1-\la_1$.

Suppose first that $I_k(\mu)=\varnothing$.
Then~\eqref{eqn:app-expanded-edge} becomes
\[
 \frac12\sum_i\Delta_i^2+
 \sum_{i\ge k-\la_1+2}\la_i-\ell(\la)+\ell(\mu).
\]
This is nonnegative if $\ell(\la)\le\ell(\mu)$.
Otherwise, $\frac12\sum_i\Delta_i^2$ is at least the number of boxes
in the rows of $\la$ below row $\ell(\mu)$ (by $z^2\ge|z|$ for integral $z$ and $\sum_i\Delta_i=0$), and hence at least
$\ell(\la)-\ell(\mu)$. Thus, this is nonnegative in all cases.

Next, suppose that $I_k(\mu)\ne\varnothing$ and
$\mu_{k-\mu_1+2}=1$. Put $e=|I_k(\mu)|$. Then $\mu_i=1$ for every $k-\mu_1+2 \le i \le \ell(\mu)$, and hence
$\sum_{i\ge k-\mu_1+2}\mu_i=\ell(\mu)-k+\mu_1-1$. We have
\[
\Delta_1 - \sum_{i \ge k - \mu_1 + 2} \mu_i + \ell ( \mu ) = \Delta_1 - \ell(\mu) +  k - \mu_1 + 1 + \ell ( \mu ) = k - \la_1 + 1.
\]
Thus, the left-hand side of~\eqref{eqn:app-edge-one} is
\begin{equation}\label{eqn:app-column-edge}
 \frac12\sum_i\Delta_i^2+\sum_{i=2}^e\Delta_i
 +\sum_{i\ge k-\la_1+2}\la_i-\ell(\la)-\la_1+k+1.
\end{equation}
The sum of the first two terms is nonnegative by the same argument,
and the sum of the remaining terms is nonnegative
by~\eqref{eqn:app-tail}.

It remains to consider $\mu_{k-\mu_1+2}\ge2$.
Proposition~\ref{prop:app-low-mass} then gives $I_k(\mu)=\{1\}$.
Put $p=k-\mu_1+2$.
The sum in the definition of $\Phi(\la)$ starts at $p+\Delta_1$. After adding $\ell(\la)-\ell(\mu)$ to the left-hand side
of~\eqref{eqn:app-edge-one}, and using $\frac12\sum_i \Delta_i = 0$, we obtain
\begin{equation}\label{eqn:app-rowplus}
 \sum_{i\ge p+\Delta_1}\binom{\Delta_i}{2}
 +\sum_{i<p+\Delta_1}\binom{\Delta_i+1}{2}
 +\sum_{p\le i<p+\Delta_1}(\mu_p-\mu_i),
 \qquad \Delta_1\ge0,
\end{equation}
or
\begin{equation}\label{eqn:app-rowminus}
 \sum_{i\ge p}\binom{\Delta_i}{2}
 +\sum_{p+\Delta_1\le i<p}
       \left[\binom{\Delta_i}{2}+\mu_i-\mu_p\right]
 +\sum_{i<p+\Delta_1}\binom{\Delta_i+1}{2},
 \qquad \Delta_1<0.
\end{equation}
Here $\binom{z}{2}=z(z-1)/2$ for every integer $z$.
All displayed summands are nonnegative.
Each row of $\la$ below row $\ell(\mu)$ contributes at least one.
Indeed, it contributes $\binom{-\la_i}{2}\ge1$, unless its index
is less than $p+\Delta_1$ in~\eqref{eqn:app-rowplus}; in that case,
the summand $\mu_p-\mu_i=\mu_p\ge2$ suffices.
Thus these expressions are at least
$\max\{\ell(\la)-\ell(\mu),0\}\ge\ell(\la)-\ell(\mu)$, so that the
left-hand side of~\eqref{eqn:app-edge-one} is nonnegative.

Suppose now that equality holds along every edge of a nonconstant
cycle. By removing consecutive repetitions, we may assume that
adjacent vertices are distinct.
If $I_k(\mu)=\varnothing$, equality along $\la\to\mu$ forces
$\ell(\la)>\ell(\mu)$ and
$\sum_{i\ge k-\la_1+2}\la_i=0$.
Hence $I_k(\la)=\varnothing$ as well.
Tracing predecessors around the cycle would then give a strict
decrease in length along every edge, which is impossible.
Thus no vertex of the cycle has $I_k=\varnothing$.

Suppose next that $\mu_{k-\mu_1+2}=1$.
Equality in~\eqref{eqn:app-column-edge} implies equality
in~\eqref{eqn:app-tail} and
\[
 \Delta_i\in\{0,-1\}\quad(2\le i\le e),\qquad
 \Delta_i\in\{0,1\}\quad(i=1\text{ or }i>e).
\]
There are equally many positive and negative entries.
Since $\la_1\le\mu_1$,
$k-\la_1+2\ge k-\mu_1+2>e$.
We have already shown that $I_k(\la)\ne\varnothing$.
Hence equality in~\eqref{eqn:app-tail} forces all parts of $\la$
from row $k-\la_1+2$ to row $\ell(\la)$ to equal one.
Moreover, $\ell(\la)\le\ell(\mu)$ and, for $i>e$,
\[
 i+k-\la_i+1\ge i+k-\mu_i+1>\ell(\mu)\ge\ell(\la).
\]
Thus $|I_k(\la)|\le e$ and $\la_{k-\la_1+2}=1$.
Along every such edge, the triple
\[
 \bigl(|I_k(\la)|,\la_1,\mathsf{n}(\la)\bigr)
\]
strictly increases in lexicographic order. Indeed, if its first
two entries do not change, the signs of $\Delta_i$ show that boxes
move from rows of index at most $e$ to rows of index greater than $e$.
Tracing predecessors again, we see that these vertices cannot
occur on the cycle either.

It remains to consider $\mu_{k-\mu_1+2}\ge2$.
Equality in~\eqref{eqn:app-rowplus}--\eqref{eqn:app-rowminus} gives
\begin{equation}\label{eqn:app-thick-tail-equality}
 \ell(\la)\ge\ell(\mu),\qquad
 \Delta_1\in\{0,-1\},\qquad
 \la_{p+\Delta_1}\in\{\mu_p-1,\mu_p\}.
\end{equation}
Indeed, the rows of $\la$ below row $\ell(\mu)$ already account
for the length difference, so the summand
$\binom{\Delta_1+1}{2}$ corresponding to row one must vanish.
The summand at $p$ when $\Delta_1=0$, or at $p-1$ when
$\Delta_1=-1$, gives the last condition.
Consequently, the tuple
\[
 \bigl(\la_{k-\la_1+2},
             -\ell(\la),-\la_1,\mathsf{n}(\la)\bigr)
\]
strictly increases in lexicographic order along every remaining edge.
If the first three entries do not change, the equality conditions give
$\Delta_i\in\{-1,0\}$ above $p$ and $\Delta_i\in\{0,1\}$
from $p$ onward, so boxes again move to rows of larger index.
This is incompatible with a cycle and completes the proof.
\end{proof}

\begin{thm}\label{thm:app-b-cycle}
For every nonconstant directed cycle $\Gamma$ in $\Par_m^{(k)}$,
\begin{equation}\label{eqn:app-b-cycle}
 \sum_{(\la\to\mu)\in\Gamma}b_k(\la,\mu)\ge1.
\end{equation}
\end{thm}
\begin{proof}
When we sum~\eqref{eqn:app-weighted-bound} around the cycle, the
$\mathsf{n}$ terms cancel. By~\eqref{eqn:app-edge-one}, we obtain
\[
 (k+1)\sum_{(\la\to\mu)\in\Gamma}b_k(\la,\mu)
 \ge\sum_{(\la\to\mu)\in\Gamma}
       \bigl[\Phi(\mu)-\Phi(\la)\bigr]=0.
\]
The $\Phi$ terms also cancel. By Proposition~\ref{prop:app-edge},
at least one of the edge inequalities is strict.
Hence the sum of the $b_k$ values is positive.
Since these values are integers, the assertion follows.
\end{proof}

\subsection{The degree estimate and Hypothesis~\ref{hypo:cycle}}

We now compare $b_k$ with the function $e_k$ defined in
Section~\ref{sec:comb}, using the projective modules of
Subsection~\ref{subsec:Amod}.

\begin{prop}\label{prop:app-energy}
For each pair $\la,\mu\in\Par_m^{(k)}$, we have
\begin{equation}\label{eqn:app-energy}
 e_k(\la,\mu)\ge b_k(\la,\mu).
\end{equation}
\end{prop}
\begin{proof}
For $\la,\nu\in\Par_m$, the representations $L_\la$ and $L_\nu$
are direct summands of
$\C \Sym_m \otimes_{\Sym_{\la'}} \mathsf{sgn}$ and
$\C \Sym_m \otimes_{\Sym_{\nu}}\mathsf{triv}$, respectively.
Since $P_\la\cong\C[X_1,\ldots,X_m]\otimes L_\la$, it follows
from~\cite[Proposition~4.2.8\textup{(1)}]{Che10} that
\begin{equation}\label{eqn:app-big-schur-bound}
 \operatorname{ord}_q[P_\la:L_\nu]_q\ge d(\la,\nu).
\end{equation}
Taking $\nu=(\mu^{\omega_k})'$ and subtracting $d_k(\mu)$ gives
the assertion by the definition of $e_k$ in Section~\ref{sec:comb}.
\end{proof}

\begin{thm}\label{thm:app-cycle-m2k}
Hypothesis~\ref{hypo:cycle} holds for $m\leq2k$.
\end{thm}
\begin{proof}
For a sequence as in Hypothesis~\ref{hypo:cycle},
Proposition~\ref{prop:app-energy} and Theorem~\ref{thm:app-b-cycle} give
\begin{equation}\label{eqn:app-cycle}
 \sum_{i=1}^{l} e_k(\la^{(i)},\la^{(i+1)})
 \ge\sum_{i=1}^{l} b_k(\la^{(i)},\la^{(i+1)})\ge1.
\end{equation}
This proves~\eqref{eqn:loop}. The same argument applies to any
nonconstant cycle, even if some vertices are repeated.
\end{proof}

\subsection*{Acknowledgments}
This paper arose from discussions with Thomas Lam and owes its existence to him; the author is deeply grateful for his insights. The author would also like to thank Anton Khoroshkin for helpful discussions. The example computations in this paper were performed using SageMath~10.4. This work was supported in part by JSPS KAKENHI Grant Number JP24K21192. Part of this work was carried out during the author's visit to the Sydney Mathematical Research Institute in May--June 2025, whose hospitality is gratefully acknowledged.

\subsection*{Disclosure of AI assistance}
ChatGPT (OpenAI; GPT-5.6 and GPT-6) was used to draft proofs and exposition for Appendix~A. ChatGPT (OpenAI; GPT-5.6, GPT-6, and GPT-6.1) and Claude Fable 5.1 (Anthropic) were also used for expository revisions throughout the manuscript, including rewriting proofs for clarity. The author edited the resulting material, independently checked all mathematical arguments, and takes full responsibility for the final manuscript.

{\footnotesize
\bibliography{ref}
\bibliographystyle{plain}}
\end{document}